\documentclass{article}

\usepackage{amssymb}
\usepackage{latexsym}

\usepackage{doi}

\usepackage[numbers]{natbib}

\usepackage{anyfontsize}
\usepackage{stackengine}
\usepackage{siunitx} 
\usepackage{multirow}
\usepackage{xspace}
\usepackage{manfnt}
\usepackage{graphicx}
\usepackage{amsmath}
\usepackage{amsfonts}
\usepackage{url}

\usepackage{pifont}
\usepackage{tikz}
\usepackage{pgffor}
\usetikzlibrary{shapes.geometric, arrows,positioning}
\usetikzlibrary{patterns}
\tikzstyle{codeflowbox} = [draw, thick, text width=0.6\textwidth, inner sep=6mm, align=left]
\tikzstyle{boxRK} = [rectangle, rounded corners, minimum width=2cm, minimum height=2cm, text centered, text width=2cm, draw=black, fill=white!30]
\tikzstyle{loop} = [diamond, minimum width=4cm, aspect=2.5, text centered, draw=black, fill=green!30]
\tikzstyle{stextr} = [rectangle, minimum width=4cm, text centered, draw=red, text=red, fill=white,opacity=0,text opacity=1]
\tikzstyle{stext} = [rectangle, minimum width=4cm, text centered, draw=black, fill=white,opacity=0,text opacity=1]
\tikzstyle{arrow} = [thick,->,>=stealth]

\newcommand{\Rey}{\mathrm{Re}}
\newcommand{\gscdc}{0.75}
\newcommand{\eps}{\epsilon}
\newcommand{\deps}{d\epsilon}

\newcommand{\dts}{\Delta t^*}

\newcommand{\length}{L}
\newcommand{\lx}{\length_x}
\newcommand{\ly}{\length_y}
\newcommand{\lz}{\length_z}

\newcommand{\res}{N}
\newcommand{\resx}{\res_x}
\newcommand{\resy}{\res_y}
\newcommand{\resz}{\res_z}

\newcommand{\bigO}{{\mathcal{O}}}

\newcommand{\presup}{\Sigma}

\newcommand{\subT}{ijk}
\newcommand{\supTu}{\presup100}
\newcommand{\supTuu}{\presup200}
\newcommand{\supTuuu}{\presup300}
\newcommand{\supTuuuu}{\presup400}
\newcommand{\supTuv}{\presup110}
\newcommand{\supTuuv}{\presup210}
\newcommand{\supTuuuv}{\presup310}
\newcommand{\supTuuvv}{\presup220}
\newcommand{\supTuvw}{\presup111}
\newcommand{\supTuuvw}{\presup211}

\newcommand{\subW}{ij}
\newcommand{\supWu}{\presup10}
\newcommand{\supWuu}{\presup20}
\newcommand{\supWuuu}{\presup30}
\newcommand{\supWuuuu}{\presup40}
\newcommand{\supWuv}{\presup11}
\newcommand{\supWuuv}{\presup21}
\newcommand{\supWuuuv}{\presup31}
\newcommand{\supWuuvv}{\presup22}

\newcommand{\subO}{i}

\newcommand{\supOu}{\presup1}
\newcommand{\supOuu}{\presup2}
\newcommand{\supOuuu}{\presup3}
\newcommand{\supOuuuu}{\presup4}

\def\showpng{1}
\def\showtikz{1}

\newcommand{\surf}[1]{{\tikz{\node[anchor=north] (char) at (0.,0.15) {\small$\mathop{#1}$}; \draw (-0.075,-0.025)--(-0.15,-0.1)--(0.05,-0.1)--(0.15,-0.)--(0.05,0);}}}
\makeatletter

\newcommand{\osetT}[2]{%
  {\mathop{{}#2}\limits^{\vbox to 5.\ex@{\kern-\tw@\ex@
   \hbox{\scriptsize #1}\vss}}}}

\newcommand{\osetL}[2]{%
  {\mathop{{}#2}\limits^{\vbox to 0.1\ex@{\kern-\tw@\ex@
   \hbox{\scriptsize $#1$}\vss}}}}
   
\newcommand{\osetLpq}[2]{%
  {\mathop{{}#2}\limits^{\vbox to 3.\ex@{\kern-\tw@\ex@
   \hbox{\scriptsize $#1$}\vss}}}}

\newcommand{\oset}[2]{%
  {\mathop{{}#2}\limits^{\vbox to 4.\ex@{\kern-\tw@\ex@
   \hbox{\scriptsize #1}\vss}}}}
   
\makeatother

\newcommand{\area}[4]{{\oset{\surf{#2}}{#1}}{\vphantom{#1}}^{\!#3}_{\!#4}}

\newcommand{\areaT}[4]{{\osetT{\surf{#2}}{#1}}{\vphantom{#1}}^{\!\!#3}_{\!\!#4}}

\newcommand{\vol}[1]{\overset{\mbox{\mancube}}{#1}}

\newcommand{\Qpt}{q}

\newcommand{\Qvol}{\vol{\Qpt}}

\newcommand{\QvolT}{\Qvol_{\subT}}

\newcommand{\QvolTe}{\QvolT}
\newcommand{\QvolTu}{\QvolTe^{\supTu}}
\newcommand{\QvolTuu}{\QvolTe^{\supTuu}}
\newcommand{\QvolTuuu}{\QvolTe^{\supTuuu}}
\newcommand{\QvolTuuuu}{\QvolTe^{\supTuuuu}}
\newcommand{\QvolTuv}{\QvolTe^{\supTuv}}
\newcommand{\QvolTuuv}{\QvolTe^{\supTuuv}}
\newcommand{\QvolTuuuv}{\QvolTe^{\supTuuuv}}
\newcommand{\QvolTuuvv}{\QvolTe^{\supTuuvv}}
\newcommand{\QvolTuvw}{\QvolTe^{\supTuvw}}
\newcommand{\QvolTuuvw}{\QvolTe^{\supTuuvw}}

\newcommand{\QptT}{\Qpt_{\subT}}

\newcommand{\QptTe}{\QptT}
\newcommand{\QptTu}{\QptTe^{\supTu}}
\newcommand{\QptTuu}{\QptTe^{\supTuu}}
\newcommand{\QptTuuu}{\QptTe^{\supTuuu}}
\newcommand{\QptTuuuu}{\QptTe^{\supTuuuu}}
\newcommand{\QptTuv}{\QptTe^{\supTuv}}
\newcommand{\QptTuuv}{\QptT^{\supTuuv}}
\newcommand{\QptTuuuv}{\QptTe^{\supTuuuv}}
\newcommand{\QptTuuvv}{\QptTe^{\supTuuvv}}
\newcommand{\QptTuvw}{\QptTe^{\supTuvw}}
\newcommand{\QptTuuvw}{\QptTe^{\supTuuvw}}

\newcommand{\QavgW}{\area{\Qpt}{z}{}{\subW}}

\newcommand{\QavgWu}{\area{\Qpt}{z}{\supWu}{\subW}}
\newcommand{\QavgWuu}{\area{\Qpt}{z}{\supWuu}{\subW}}
\newcommand{\QavgWuuu}{\area{\Qpt}{z}{\supWuuu}{\subW}}
\newcommand{\QavgWuuuu}{\area{\Qpt}{z}{\supWuuuu}{\subW}}
\newcommand{\QavgWuv}{\area{\Qpt}{z}{\supWuv}{\subW}}
\newcommand{\QavgWuuv}{\area{\Qpt}{z}{\supWuuv}{\subW}}
\newcommand{\QavgWuuuv}{\area{\Qpt}{z}{\supWuuuv}{\subW}}
\newcommand{\QavgWuuvv}{\area{\Qpt}{z}{\supWuuvv}{\subW}}

\newcommand{\QptW}{\Qpt_{\subW}}
\newcommand{\QptWe}{\QptW}
\newcommand{\QptWu}{\QptWe^{\supWu}}
\newcommand{\QptWuu}{\QptWe^{\supWuu}}
\newcommand{\QptWuuu}{\QptWe^{\supWuuu}}
\newcommand{\QptWuuuu}{\QptWe^{\supWuuuu}}
\newcommand{\QptWuv}{\QptWe^{\supWuv}}
\newcommand{\QptWuuv}{\QptWe^{\supWuuv}}
\newcommand{\QptWuuuv}{\QptWe^{\supWuuuv}}
\newcommand{\QptWuuvv}{\QptWe^{\supWuuvv}}

\newcommand{\QptO}{\Qpt_{\subO}}
\newcommand{\QptOe}{\QptO}
\newcommand{\QptOu}{\QptOe^{\supOu}}
\newcommand{\QptOuu}{\QptOe^{\supOuu}}
\newcommand{\QptOuuu}{\QptOe^{\supOuuu}}
\newcommand{\QptOuuuu}{\QptOe^{\supOuuuu}}

\newcommand{\solvar}{\vu(\vx,t)}
\newcommand{\solvarA}[1]{\vu(#1)}
\newcommand{\solvarP}{\vu(\vx,t+\Dt)}
\newcommand{\solvarM}{\vu(\vx,t-\Dt)}
\newcommand{\solvarMM}{\vu(\vx,t-2\Dt)}

\newcommand{\vE}{{\bf E}}
\newcommand{\vN}{{\bf N}}
\newcommand{\vW}{{\bf W}}
\newcommand{\vb}{{\bf b}}
\newcommand{\ve}{{\bf e}}
\newcommand{\vf}{{\bf f}}
\newcommand{\vg}{{\bf g}}

\newcommand{\vk}{{\bf k}}
\newcommand{\vm}{{\bf m}}
\newcommand{\vn}{{\bf n}}
\newcommand{\vu}{{\bf u}}
\newcommand{\vv}{{\bf v}}
\newcommand{\vw}{{\bf w}}
\newcommand{\vx}{{\bf x}}
\newcommand{\vy}{{\bf y}}
\newcommand{\vz}{{\bf z}}
\newcommand{\ex}{\ve_x}
\newcommand{\ey}{\ve_y}
\newcommand{\ez}{\ve_z}

\newcommand{\FvelL}{v}
\newcommand{\FvelFx}{\hat{\FvelL}_{\vk,x}}
\newcommand{\FvelFy}{\hat{\FvelL}_{\vk,y}}
\newcommand{\FvelFz}{\hat{\FvelL}_{\vk,z}}
\newcommand{\FvelFm}{\hat{\FvelL}_{\vk,m}}
\newcommand{\FmagFx}{\hat{\FmagL}_{\vk,x}}
\newcommand{\FmagFy}{\hat{\FmagL}_{\vk,y}}
\newcommand{\FmagFz}{\hat{\FmagL}_{\vk,z}}
\newcommand{\Fvel}{{\bf \FvelL}}
\newcommand{\Fvelx}{\FvelL_x}
\newcommand{\Fvely}{\FvelL_y}
\newcommand{\Fvelz}{\FvelL_z}
\newcommand{\fluxf}{\vf}
\newcommand{\fluxfq}{\fluxf_\Qpt}
\newcommand{\fluxfqIp}{\fluxf^+_{\Qpt,i+1/2,j,k}}
\newcommand{\fluxfqIm}{\fluxf^-_{\Qpt,i+1/2,j,k}}
\newcommand{\fluxfqIpm}{\fluxf^\pm_{\Qpt,i+1/2,j,k}}
\newcommand{\fluxLLFqijk}{\fluxf^{LLF}_{\Qpt,i+1/2,j,k}}
\newcommand{\fluxfqt}[1]{\fluxf_{#1}}
\newcommand{\FmagL}{b}
\newcommand{\Fmag}{{\bf \FmagL}}
\newcommand{\Fmagx}{\FmagL_x}
\newcommand{\Fmagy}{\FmagL_y}
\newcommand{\Fmagz}{\FmagL_z}
\newcommand{\FmagJL}{j}
\newcommand{\FmagJ}{{\bf \FmagJL}}
\newcommand{\press}{p}
\newcommand{\edens}{e}
\newcommand{\Imat}{\underline{\mathbf{I}}}
\newcommand{\vnabla}{\nabla}

\newcommand{\pal}{\partial_l}

\newcommand{\pax}{\partial_x}
\newcommand{\pay}{\partial_y}
\newcommand{\paz}{\partial_z}
\newcommand{\pat}{\partial_t}

\newcommand{\resis}{\eta}
\newcommand{\visc}{\mu}
\newcommand{\viscstress}{\underline{\boldsymbol \sigma}}

\newcommand{\beq}{\begin{equation}}
\newcommand{\eeq}{\end{equation}}

\newcommand{\eq}[1]{eq. \eqref{#1}}
\newcommand{\eqs}[2]{eqs. \eqref{#1}-\eqref{#2}}
\newcommand{\eqsN}[2]{\eqref{#1}-\eqref{#2}}
\newcommand{\eqsa}[2]{eqs. \eqref{#1} and \eqref{#2}}
\newcommand{\eqsaa}[3]{eqs. \eqref{#1}, \eqref{#2} and \eqref{#3}}

\newcommand{\autocorrelation}{autocorrelation\xspace}
\newcommand{\formulas}{formulas\xspace}
\newcommand{\timestep}{time-step\xspace}
\newcommand{\nonlinear}{nonlinear\xspace}
\newcommand{\nonideal}{non-ideal\xspace}
\newcommand{\multistep}{multi-step\xspace}

\newcommand{\tab}[1]{table \ref{#1}}
\newcommand{\Tab}[1]{Table \ref{#1}}
\newcommand{\fig}[1]{fig. \ref{#1}}
\newcommand{\Fig}[1]{Fig. \ref{#1}}
\newcommand{\SFig}[2]{Fig. \ref{#1}.$#2$}
\newcommand{\Sfig}[2]{fig. \ref{#1}.$#2$}
\newcommand{\SFigs}[3]{Fig. \ref{#1}.$#2$ and $#3$}
\newcommand{\Sfigs}[3]{fig. \ref{#1}.$#2$ and $#3$}

\newcommand{\rhs}{right-hand side\xspace}

\newcommand{\beqa}{\begin{eqnarray}}
\newcommand{\eeqa}{\end{eqnarray}}

\newcommand{\cs}{c_s}
\newcommand{\gadia}{\gamma}

\newcommand{\Dx}{\Delta x}
\newcommand{\Dy}{\Delta y}
\newcommand{\Dz}{\Delta z}
\newcommand{\dvol}{dxdydz}

\newcommand{\dd}{d}
\newcommand{\dx}{\dd x}
\newcommand{\dy}{\dd y}
\newcommand{\dz}{\dd z}
\newcommand{\dydz}{\dy\dz}
\newcommand{\dzdx}{\dz\dx}
\newcommand{\dxdy}{\dx\dy}

\newcommand{\Dt}{\Delta t}
\newcommand{\DtCFL}{\Delta t_{CFL}}
\newcommand{\DtVISCO}{\Delta t_{VISCO}}
\newcommand{\DtRESI}{\Delta t_{RESI}}
\newcommand{\DtEsink}{\Delta t_{Esink}}

\newcommand{\Surf}{S}
\newcommand{\Snx}{\Surf_{x,j,k}}
\newcommand{\Sny}{\Surf_{y,k,i}}
\newcommand{\Snz}{\Surf_{z,i,j}}

\newcommand{\Cell}{\Omega}
\newcommand{\Cellijk}{\Cell_{ijk}}
\newcommand{\Cellipjk}{\Cell_{i+1,j,k}}
\newcommand{\Cellijkdef}{[x_i-\Dx/2,x_i+\Dx/2] \times [y_j-\Dy/2,y_j+\Dy/2] \times [z_k-\Dz/2,z_k+\Dz/2]}

\newcommand{\ipjk}{i+1/2,j,k}
\newcommand{\imjk}{i-1/2,j,k}
\newcommand{\ijpk}{i,j+1/2,k}
\newcommand{\ijmk}{i,j-1/2,k}
\newcommand{\ijkp}{i,j,k+1/2}
\newcommand{\ijkm}{i,j,k-1/2}
\newcommand{\ipmjk}{i\pm1/2,j,k}

\newcommand{\lra}{\leftrightarrow}

\newcommand{\AaSx}[2]{\area{#1}{x}{}{#2}} 
\newcommand{\AaSy}[2]{\area{#1}{y}{}{#2}} 
\newcommand{\AaSz}[2]{\area{#1}{z}{}{#2}}

\newcommand{\AaFn}[3]{\area{#1}{n}{#2}{#3}} 
\newcommand{\AaFx}[3]{\area{#1}{x}{#2}{#3}} 
\newcommand{\AaFy}[3]{\area{#1}{y}{#2}{#3}} 
\newcommand{\AaFz}[3]{\area{#1}{z}{#2}{#3}}

\newcommand{\fluxFaceqn}{\AaFn{F}{\Qpt}{\phantom{i}}}

\newcommand{\fluxFaceqxipjk}{\AaFx{F}{\Qpt}{\ipjk}}
\newcommand{\fluxFaceqximjk}{\AaFx{F}{\Qpt}{\imjk}}
\newcommand{\fluxFaceqyijpk}{\AaFy{F}{\Qpt}{\ijpk}}
\newcommand{\fluxFaceqyijmk}{\AaFy{F}{\Qpt}{\ijmk}}
\newcommand{\fluxFaceqzijkp}{\AaFz{F}{\Qpt}{\ijkp}}
\newcommand{\fluxFaceqzijkm}{\AaFz{F}{\Qpt}{\ijkm}}
\newcommand{\fluxFaceqxipmjk}{\AaFx{F}{\Qpt}{\ipmjk}}

\newcommand{\clrJ}{yellow}
\newcommand{\clrB}{orange}
\newcommand{\clrE}{red}

\newcommand{\clrH}{blue}
\newcommand{\clrV}{green}
\newcommand{\clrDV}{cyan}

\newcommand{\clrU}{purple}

\newcommand{\clrBgen}{green}
\newcommand{\clrBNID}{purple}
\newcommand{\clrBreco}{orange}
\newcommand{\clrBcolella}{blue}
\newcommand{\clrBriemann}{gray}
\newcommand{\clrBDVDB}{yellow}
\newcommand{\clrBcalc}{yellow}
\newcommand{\clrBcond}{blue}
\newcommand{\clrBdiag}{orange}
\newcommand{\clrBforce}{cyan}
\newcommand{\clrForce}{cyan}

\usepackage[left=1.5cm, right=1.5cm, top=2.cm, bottom=2.cm]{geometry}

\usepackage{titling}
\thanksmarkseries{arabic}
\begin{document}

\author{J.-M. Teissier\thanks{Technische Universit\"at Berlin, ER 3-2, Hardenbergstr. 36a, D-10623 Berlin, Germany} and W.-C. M\"{u}ller$^{1,}$\thanks{Max-Planck/Princeton Center for Plasma Physics} \\ \footnotesize jm.teissier@astro.physik.tu-berlin.de}
  \date{June 28, 2024}

\title{High-order finite-volume integration schemes for subsonic magnetohydrodynamics}%

\maketitle

\begin{abstract}
We present an efficient dimension-by-dimension finite-volume method which solves the adiabatic magnetohydrodynamics equations at high discretization order, using the constrained-transport approach on Cartesian grids. Results are presented up to tenth order of accuracy. This method requires only one reconstructed value per face for each computational cell. A passage through high-order point values leads to a modest growth of computational cost with increasing discretization order. At a given resolution, these high-order schemes present significantly less numerical dissipation than commonly employed lower-order approaches. Thus, results of comparable accuracy are achievable at a substantially coarser resolution, yielding overall performance gains. We also present a way to include physical dissipative terms: viscosity, magnetic diffusivity and cooling functions, respecting the finite-volume and constrained-transport frameworks.
\end{abstract}


\section{Introduction}

Turbulence is an ubiquitous phenomenon in nature which still lacks a comprehensive theory. Progress in its understanding can be achieved by direct numerical simulation (DNS) of the Navier-Stokes equations (NSE) for neutral fluids. For electrically conductive media, such as liquid planetary cores or plasmas, ionized gases, found in many astrophysical environments, the equations of magnetohydrodynamics (MHD) are a convenient single-fluid description. This approximation is valid on large spatial and long temporal scales as compared to the microscopic characteristics of electric and magnetic interactions between the charged constituents. Therefore, MHD is widely used when dealing with geo- and astrophysical problems.

In natural systems, the Reynolds number $\Rey$, which characterizes the spectral bandwidth of turbulent fluctuations, can be orders of magnitude above the maximum values achievable by DNS on today's high-performance computing facilities. A main culprit limiting the achievable Reynolds number in DNS is insufficient numerical resolution of space and time and the consequential dissipation introduced via the numerical integration scheme and resulting in amplitude and phase errors of the numerical solution. Numerical dissipation is generated by the loss of information due to discretization errors, which is dominant at small scales and represents a lower bound on the physical dissipation one can apply in DNS via corresponding terms in the underlying differential equations.

A reduction of algorithmic numerical dissipation directly increases the spectral bandwidth of a DNS at the same numerical resolution, and ideally leads to more accurate numerical experiments at affordable computational cost. For this aim, more accurate Riemann solvers, which numerically approximate the time evolution of two neighbouring fluid states, have been developed \cite{ROE81,TSS94,POW97,TOR09,MRE15,MIM21}. Another approach is to raise the discretization order, i.e. to go beyond the standard second-order approach by including terms of higher order in the Taylor expansion of the solution (assumed to be analytic within each computational cell). Even though at a given resolution higher-order discretization is computationally more demanding, such schemes can lead to a large gain in efficiency overall. Indeed, higher-order solvers typically require a lower resolution to reach a similar accuracy as compared to second-order schemes \cite{VTH19}.

However, multi-dimensional polynomial reconstruction becomes increasingly expensive with a growing order of accuracy. A specific strategy to alleviate this problem is realized in ``dimension-by-dimension'' solvers. Instead of solving an expensive 3D reconstruction problem at once, they solve three one-dimensional cost-effective problems separately. A further gain in performance is achieved by a dimension-by-dimension approach which requires only one high-order reconstruction along each dimension \cite{COC11,BUH14}. It relies on a transformation of reconstructed area-averages to high-order point values in the middle of the faces, from which high-order point fluxes are obtained that lead to high-order area averages.

The present work employs this numerical technique. It extends a dimension-by-dimension fourth-order finite-volume solver for the compressible MHD equations using the con\-strained-transport approach \cite{VTH19} up to order ten. The schemes are based on Weighted Essentially Non-Oscillatory (WENO) 1D reconstructions, associated with a passage through point values. Volume/area/line-averages$\lra$point value transformations are derived using the method described in \cite{COC11,BUH14}.

Only subsonic flows are considered in this work. For supersonic flows, where shocks and discontinuities are prevalent, avoiding the generation of spurious oscillations near sharp gradients is a very challenging issue \cite{ROM16,VTH19,BAG16,WSH19} which we do not consider here. Our aim in this work is to provide a relatively simple framework to build computationally efficient higher-order 3D compressible MHD solvers.

The rest of this paper is organized as follows: section~\ref{sec:equations} describes the MHD equations to be solved as well as the chosen discretization (finite-volume in a Cartesian domain and constrained-transport for the magnetic field). Section~\ref{sec:numsolver} presents how the higher-order numerical schemes are built. Their numerical accuracy is verified in section~\ref{sec:numtests} through convergence tests. That section presents as well an application example: driven turbulence in a statistically stationary state. Finally, concluding remarks are given in section~\ref{sec:conclusion}.

\section{Governing equations and discretization}
\label{sec:equations}
\subsection{Compressible MHD equations}

The single fluid compressible adiabatic MHD equations describe the time evolution of the mass density $\rho$, the momentum $\rho \Fvel$ (with velocity $\Fvel$), the total energy density $\edens$ and the magnetic field $\Fmag$ through:

\newcommand{\Kforce}{\vf_K}
\newcommand{\Kforceijk}{\vf_{K,i,j,k}}
\newcommand{\Mforce}{\vf_M}
\newcommand{\Mforceijk}{\vf_{M,i,j,k}}
\newcommand{\KMforceE}{\vf_{eKM}}
\newcommand{\KMforceEijk}{\vf_{eKM,i,j,k}}
\newcommand{\Usink}{S_e}

\newcommand{\intE}{U}
\newcommand{\UsinkF}{\lambda}

\beqa
	\label{eq:beginMHD}
	\label{eq:dtrho}
    \pat \rho  &=& -\vnabla \cdot (\rho \Fvel),\\
	\label{eq:dtrhov}
    \pat (\rho \Fvel) &=& - \vnabla \cdot \left( \rho \Fvel \Fvel^T + (\press + \frac{1}{2}|\Fmag|^2)\Imat - \Fmag \Fmag^T -\viscstress \right)+\Kforce,\\
	\label{eq:dtedens}
    \pat  \edens  &=& -\vnabla \cdot \left((\edens + \press + \frac{1}{2}|\Fmag|^2)\Fvel  -  (\Fvel\cdot\Fmag)\Fmag \right)-\Usink+\KMforceE,\\
\label{eq:dtb}
    \pat \Fmag &=& - \vnabla \times (- \Fvel \times \Fmag + \resis \vnabla \times \Fmag)+\Mforce,\\
\label{eq:endMHD}
\label{eq:solB}
\vnabla \cdot \Fmag&=&0.
\eeqa

The $3 \times 3$ identity matrix is denoted by $\Imat$. The total energy density $\edens$ is the sum of the kinetic, $\rho |\Fvel|^2/2$, magnetic, $|\Fmag|^2/2$ and internal, $\press/(\gadia - 1)$, energy densities. Thus, the thermal pressure, $\press$, is given by:

\beq
	\label{eq:pressure}
	\press=(\gadia-1)(\edens-\frac12 \rho|\Fvel|^2-\frac12 |\Fmag|^2).
\eeq

In driven turbulence, energy is injected by a kinetic $\Kforce$ and a magnetic $\Mforce$ forcing. These two forcing mechanisms also affect the total energy density through the $\KMforceE$ term; these are described in section \ref{sec:forcing}. A statistically stationary state is reached through the dynamical balance between the forcing terms and dissipation (both of numerical and physical nature). The physical dissipative terms are represented by the magnetic diffusivity $\resis$ (caused by resistive effects), the internal energy sink modelled as Stefan-Boltzmann-like radiative losses:
\beq
\label{eq:Usink}
\Usink=\UsinkF \intE^4,
\eeq
with $\intE=\press/(\gadia-1)$ and $\UsinkF$ a constant, as well as the divergence of the linear viscous stress tensor:

\beq
	\label{eq:divviscstress}
	\vnabla \cdot \viscstress = \visc \vnabla^2 \Fvel+\frac{1}{3} \visc\vnabla(\vnabla \cdot \Fvel),
\eeq

with $\visc$ the dynamic viscosity (the bulk viscosity is set to zero here). 

In the isothermal case $\gadia \to 1$, \eq{eq:dtedens} does not need to be solved and the pressure is defined by $\press=\rho c_s^2$ with $c_s$ the constant sound-speed.

\subsection{Finite-volume discretization}

The hydrodynamic quantities $(\rho,\rho \Fvel)$ and the total energy density $\edens$ are discretized as volume averages in a Cartesian coordinate system. The cubic domain $[0,\lx]\times[0,\ly]\times[0,\lz]$ consists of $\resx \times \resy \times \resz$ cells in the $\vx$-, $\vy$- and $\vz$- directions with constant grid-sizes $\Dx=\lx/\resx, \Dy=\ly/\resy$ and $\Dz=\lz/\resz$. The cells are centered at $x_i=(i+1/2)\Dx, y_j=(j+1/2)\Dy, z_k=(k+1/2)\Dz$. Their volume is $\Cellijk=\Cellijkdef$.

As appropriate for finite-volume solvers, \eqs{eq:dtrho}{eq:dtedens} are written in conservative form: $\pat \Qpt=- \vnabla \cdot \fluxfq$. By the divergence theorem, the volume average of $\Qpt$ over the cell $\Cellijk$, denoted by $\QvolT$, evolves in time as:

\beqa
\nonumber \pat \QvolT&=&\frac{1}{\Dx\Dy\Dz}\iiint_{\Cellijk} \pat \Qpt \dvol,\\
\label{eq:patQavgT}	&=&-\frac{\fluxFaceqxipjk-\fluxFaceqximjk}{\Dx}-\frac{\fluxFaceqyijpk-\fluxFaceqyijmk}{\Dy}-\frac{\fluxFaceqzijkp-\fluxFaceqzijkm}{\Dz},
\eeqa

with $\fluxFaceqn, n\in\{x,y,z\}$ representing the flux' area average on the cell face normal to the $\vn$-direction. For example:

\beq
	\fluxFaceqxipmjk=\frac{1}{\Dy\Dz}\iint_{\Snx} \fluxfq(x_i\pm\Dx/2,y,z) \cdot \ex \dydz,
\eeq

	with $\Snx=[y-\Dy/2,y+\Dy/2]\times[z-\Dz/2,z+\Dz/2]$ and $\ex$ the unit vector along the $\vx$-direction. The averaged fluxes are defined in an analogous manner in the other directions. The explicit expressions for the fluxes, projected along $\ex$, are derived from \eqs{eq:dtrho}{eq:dtedens}:

\beq
	\label{eq:fluxexprx}
\begin{pmatrix} \fluxfqt{\rho} \cdot \ex \\ \fluxfqt{\rho\FvelL_x} \cdot \ex\\ \fluxfqt{\rho\FvelL_y} \cdot \ex\\ \fluxfqt{\rho\FvelL_z} \cdot \ex\\ \fluxfqt{\edens} \cdot \ex \end{pmatrix} = \begin{pmatrix} \rho\Fvelx \\ \rho\Fvelx^2+\press+|\Fmag|^2/2-\Fmagx^2 - \visc(\pax \Fvelx + \frac{1}{3}\vnabla \cdot \Fvel)\\ \rho\Fvelx\Fvely-\Fmagx\Fmagy -\visc\pax\Fvely \\ \rho\Fvelx\Fvelz-\Fmagx\Fmagz -\visc\pax\Fvelz\\ (\edens+\press+|\Fmag|^2/2)\Fvelx-\Fmagx(\Fvel\cdot\Fmag) \end{pmatrix}.
\eeq

Appropriate circular permutations of the spatial dimensions $(x\to y \to z \to x)$ give the projections along $\ey$ and $\ez$.

The finite-volume approach guarantees that the physically conserved quantities (total mass, total impulse, total energy) are numerically conserved as well, up to machine precision. Indeed: everything that exits a cell through an interfacial flux enters a neighbouring one (e.g. $\fluxFaceqxipjk$ is a loss for cell $\Cellijk$ but a gain for cell $\Cellipjk$).

\subsection{Constrained-transport discretization}
The constrained-transport approach discretizes the magnetic field such that its divergence (\eq{eq:solB}) is conserved up to machine precision \cite{EVH88}. Each magnetic field component is discretized as an area average on the faces normal to its respective direction. For cell $\Cellijk$:

\label{sec:CTd}

\beqa
	\label{eq:defBx}
	\AaSx{\FmagL}{x,i-1/2,j,k}&=&\iint_{\Snx} \FmagL_x(x_i-\Dx/2,y,z)\dydz,\\
	\AaSy{\FmagL}{y,i,j-1/2,k}&=&\iint_{\Sny} \FmagL_y(x,y_j-\Dy/2,z)\dzdx,\\
	\AaSz{\FmagL}{z,i,j,k-1/2}&=&\iint_{\Snz} \FmagL_z(x,y,z_k-\Dz/2)\dxdy,
\eeqa

Applying Stoke's theorem on \eq{eq:dtb} yields:

\newcommand{\linea}[3]{\osetL{#2}{\overline{#1}}_{#3}}
\newcommand{\lineapq}[3]{\osetLpq{#2}{\overline{#1}}_{#3}}
\newcommand{\lineasup}[4]{\osetL{#2}{\overline{#1}}^{#4}_{#3}}

\newcommand{\EfL}{E}
\newcommand{\EfLA}[2]{\linea{\EfL}{#1}{#2}}
\newcommand{\bfLA}[2]{\linea{\FmagL}{#1}{#2}}
\newcommand{\EfLApq}[2]{\lineapq{\EfL_#2}{#1}{}}
\newcommand{\bfLApq}[2]{\lineapq{\FmagL_#2}{#1}{}}
\beq
	\label{eq:dtmagCT} \pat \AaSx{\FmagL}{x,i-1/2,j,k}=-\frac{\EfLA{z}{z,i-1/2,j+1/2,k}-\EfLA{z}{z,i-1/2,j-1/2,k}}{\Dy}+\frac{\EfLA{y}{y,i-1/2,j,k+1/2}-\EfLA{y}{y,i-1/2,j,k-1/2}}{\Dz},
\eeq

with $\vE=-\Fvel \times \Fmag + \resis \nabla \times \Fmag$ the electric field, the overbar meaning a line-average along its respective direction, for example:

\beq
	\EfLA{z}{z,i-1/2,j+1/2,k}=\frac{1}{\Dz} \int_{z_k-\Dz/2}^{z_k+\Dz/2} \EfL_{z}(x_{i}-\Dx/2,y_j+\Dy/2,z) \dz.
\eeq

Similar relations apply for $\pat \AaSy{\FmagL}{y}$ and $\pat \AaSz{\FmagL}{z}$ with appropriate variable permutations. This staggered area-average definition of the magnetic field components leads to the conservation of a discretized formulation of $\nabla \cdot \Fmag$ up to second-order approximation:

\beq
	(\nabla \cdot \Fmag)_{i,j,k} \approx \frac{\AaSx{\FmagL}{x,i+1/2,j,k}-\AaSx{\FmagL}{x,i-1/2,j,k}}{\Dx}+\frac{\AaSy{\FmagL}{y,i,j+1/2,k}-\AaSy{\FmagL}{y,i,j-1/2,k}}{\Dy}+\frac{\AaSz{\FmagL}{z,i,j,k+1/2}-\AaSz{\FmagL}{z,i,j,k-1/2}}{\Dz},
\eeq

since the terms of this approximation cancel pairwise \cite{EVH88,ZIE04} (e.g. the circulation of $\vE$ on the edge $x=(i-1/2)\Dx, y=(j-1/2)\Dy, z\in [(k-1/2)\Dz,(k+1/2)\Dz]$ enters $\pat \AaSx{\FmagL}{x,i-1/2,j,k}$ and $\pat \AaSy{\FmagL}{y,i,j-1/2,k}$ with opposite signs). In this way, the con\-strained-transport approach is a counterpart to the finite-volume approach, where the flux exiting a cell cancels out with the one entering a neighbouring one. This implies that, if the divergence of $\Fmag$ is zero initially, it remains zero up to machine precision at all times.

\newcommand{\module}[1]{{\tt #1}\xspace}
\newcommand{\modBinterp}{\module{Binterp}}
\newcommand{\modiBinterp}{\module{(Binterp)}}
\newcommand{\modHBreco}{\module{HBreco}}
\newcommand{\modHBflux}{\module{HBflux}}
\newcommand{\modiHBflux}{\module{(HBflux)}}
\newcommand{\modFor}{\module{Forcing}}
\newcommand{\modiFor}{\module{(Forcing)}}
\newcommand{\modtoHBrhs}{\module{toHBrhs}}
\newcommand{\moditoHBrhs}{\module{(toHBrhs)}}
\newcommand{\modCT}{\module{CT}}
\newcommand{\modiCT}{\module{(CT)}}
\newcommand{\modWENO}{\module{WENO}}
\newcommand{\modiWENO}{\module{(WENO)}}
\newcommand{\modSumF}{\module{SumFlux}}
\newcommand{\modUpdateF}{\module{UpdateF}}
\newcommand{\modUpdateR}{\module{UpdateR}}
\newcommand{\modFEMK}{\module{TotE}}
\newcommand{\modVforce}{\module{PtoV}}

\newcommand{\modRESI}{\module{RESI}}
\newcommand{\modVISCO}{\module{VISCO}}
\newcommand{\modEsink}{\module{Esink}}

\newcommand{\modiRESI}{\module{(RESI)}}
\newcommand{\modiVISCO}{\module{(VISCO)}}
\newcommand{\modiEsink}{\module{(Esink)}}

\newcommand{\modAtoP}{\module{AtoP}}
\newcommand{\modiAtoPPtoA}{\module{(AtoP} and \module{PtoA)}}
\newcommand{\modLLF}{\module{1DRieS}}
\newcommand{\modMLLF}{\module{2DRieS}}
\newcommand{\modiLLF}{\module{(1DRieS)}}
\newcommand{\modPtoA}{\module{PtoA}}

\newcommand{\modPvbtoAE}{\module{PvbtoAE}}

\newcommand{\modPtoAv}{\module{PtoAv}}
\newcommand{\modAvtoVdv}{\module{AvtoVdv}}
\newcommand{\modAbtoVj}{\module{AbtoVj}}
\newcommand{\modVtoP}{\module{VtoP}}
\newcommand{\modPtoVSe}{\module{PtoVSe}}

\newcommand{\modRHS}{\module{RHS}}
\newcommand{\modSSPRK}{\module{SSPRK}}
\newcommand{\modiSSPRK}{\module{(SSPRK)}}

\input{tikzsets_arxiv.tex}

\ifx\showtikz\undefined
\else
\def\bsize{1.3}
\def\xminBBlockRHS{1.5+\lreco*0.25}
\def\xmaxBBlockRHS{\baseDUDT+\wDUDT+\resarrXoffs*0.5}
\def\yminBBlockRHS{\yarrLO-1.5}
\def\ymaxBBlockRHS{\tblockSe+0.2}

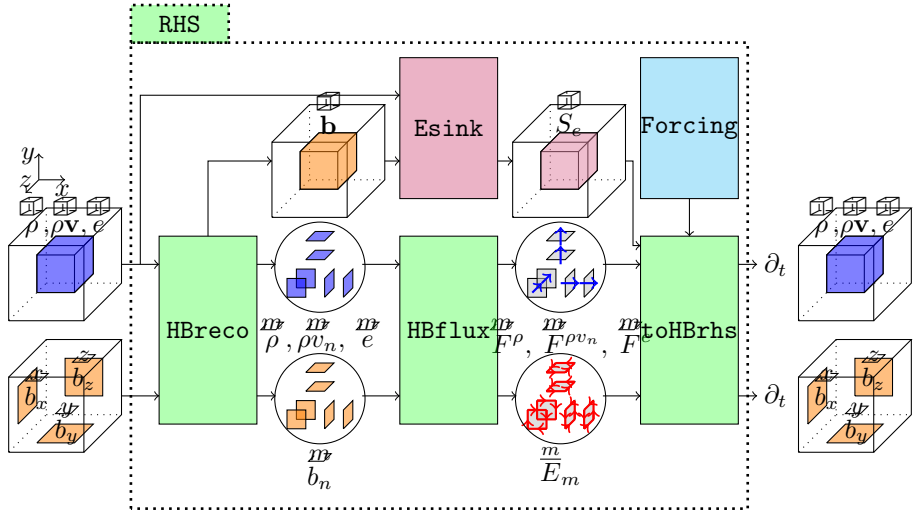
\begin{figure}[h]
\centering
\begin{tikzpicture}
\draw[dotted, line width=1pt] (\xminBBlockRHS,\yminBBlockRHS) rectangle (\xmaxBBlockRHS,\ymaxBBlockRHS);
\draw[dotted, line width=1pt,fill=\clrBgen,fill opacity=0.3] (\xminBBlockRHS,\ymaxBBlockRHS) rectangle (\xminBBlockRHS+\bsize,\ymaxBBlockRHS+0.5);
\node at (\xminBBlockRHS+\bsize*0.5,\ymaxBBlockRHS+0.25) {\modRHS};
\path (0,\yarrHI-\yarrcROffs) pic {volcube=\clrH};

\node at (0+\volnodeX,\yarrHI-\yarrcROffs+\volnodeY) {$\volH$};

\path (0,\yarrLO-\yarrcROffs) pic {areacube=\clrB};

\node at (0+\aXnodeX,\yarrLO-\yarrcROffs+\aXnodeY) {$\areaTbx$};
\node at (0+\aYnodeX,\yarrLO-\yarrcROffs+\aYnodeY) {$\areaTby$};
\node at (0+\aZnodeX,\yarrLO-\yarrcROffs+\aZnodeY) {$\areaTbz$};

\path (0.4,3.75) pic[scale=0.35] {axes};

\draw[->] (1.5,\yarrHI)--(1.5+\lreco,\yarrHI);
\draw[->] (1.5,\yarrLO)--(1.5+\lreco,\yarrLO);

\draw[->] (1.5+\lreco*0.5,\yarrHI)--(1.5+\lreco*0.5,\yarrHI+\harrSe)--(\baseFCT,\yarrHI+\harrSe);

\draw[fill=\clrBgen,fill opacity=0.3] (1.5+\lreco,\bblock) rectangle (1.5+\lreco+\wreco,\tblock);
\node at (1.5+\lreco+\wreco*0.5,1.75) {\modHBreco};

\draw[->] (1.5+\lreco+\wreco*0.5,\tblock) -- (1.5+\lreco+\wreco*0.5,\tblock+\harrbSe) -- (1.5+\lreco+\wreco+\rescXoffs,\tblock+\harrbSe);
\path (1.5+\lreco+\wreco+\rescXoffs,\mblockSe-\yarrcROffs) pic {volcube=\clrB};
\node at (1.5+\lreco+\wreco+\rescXoffs+\volnodeX,\mblockSe-\yarrcROffs+\volnodeY) {$\vol{\Fmag}$};
\draw[->] (1.5+\lreco+\wreco+\rescXoffs+1.5,\tblock+\harrbSe) -- (\baseFCT,\tblock+\harrbSe);

\draw[->] (1.5+\lreco+\wreco,\yarrHI) -- (1.5+\lreco+\wreco+\resarrXoffs,\yarrHI);
\draw[->] (1.5+\lreco+\wreco,\yarrLO) -- (1.5+\lreco+\wreco+\resarrXoffs,\yarrLO);

\path (1.5+\lreco+\wreco+\resXoffs,\yarrHI-\yarrROffs) pic {reco=\clrH};
\path (1.5+\lreco+\wreco+\resXoffs,\yarrLO-\yarrROffs) pic {reco=\clrB};

\node[anchor=north] at (1.5+\lreco+\wreco+\resXoffs+\underCircleX,\yarrHI-\yarrROffs-\underCircleY) {$\areaTH$};
\node[anchor=north] at (1.5+\lreco+\wreco+\resXoffs+\underCircleX,\yarrLO-\yarrROffs-\underCircleY) {$\areaTB$};

\draw[->] (1.5+\lreco+\wreco+\resXoffs+1.,\yarrHI)--(\baseFCT,\yarrHI);
\draw[->] (1.5+\lreco+\wreco+\resXoffs+1.,\yarrLO)--(\baseFCT,\yarrLO);

\draw[fill=\clrBNID,fill opacity=0.3] (\baseFCT,\bblockSe) rectangle (\baseFCT+\wFCT,\tblockSe);
\node at (\baseFCT+\wFCT*0.5,\tblockSe*0.5+\bblockSe*0.5) {\modEsink};

\draw[->] (\baseFCT+\wFCT,\mblockSe) -- (\baseFCT+\wFCT+\resarrcXoffs,\mblockSe);
\path (\baseFCT+\wFCT+\rescXoffs,\mblockSe-\yarrcROffs) pic {volcube=\clrU!50};
\node at (\baseFCT+\wFCT+\rescXoffs+\volnodeX,\mblockSe-\yarrcROffs+\volnodeY) {$\vol{\Usink}$};

\draw[->] (\baseFCT+\wFCT+\rescXoffs+1.5,\mblockSe) -- (\baseFCT+\wFCT+\rescXoffs+1.5+0.1,\mblockSe)-- (\baseFCT+\wFCT+\rescXoffs+1.5+0.1,\yarrHI+0.25)--(\baseDUDT,\yarrHI+0.25);

\draw[->] (\baseFCT+\wFCT,\yarrHI) -- (\baseFCT+\wFCT+\resarrXoffs,\yarrHI);
\draw[->] (\baseFCT+\wFCT,\yarrLO) -- (\baseFCT+\wFCT+\resarrXoffs,\yarrLO);
\path (\baseFCT+\wFCT+\resXoffs,\yarrHI-\yarrROffs) pic {hflux=gray!50/\clrH};
\node[anchor=north] at (\baseFCT+\wFCT+\resXoffs+\underCircleX+0.1,\yarrHI-\yarrROffs-\underCircleY) {$\areaTFH$};
\path (\baseFCT+\wFCT+\resXoffs,\yarrLO-\yarrROffs) pic {bflux=gray!50/\clrE};
\node[anchor=north] at (\baseFCT+\wFCT+\resXoffs+\underCircleX,\yarrLO-\yarrROffs-\underCircleYlin) {$\lineTEm$};

\draw[fill=\clrBgen,fill opacity=0.3] (\baseFCT,\bblock) rectangle (\baseFCT+\wFCT,\tblock);
\node at (\baseFCT+\wFCT*0.5,1.75) {\modHBflux};

\draw[fill=\clrBforce,fill opacity=0.3] (\baseDUDT,\bblockSe) rectangle (\baseDUDT+\wDUDT,\tblockSe);
\node at (\baseDUDT+\wDUDT*0.5,\tblockSe*0.5+\bblockSe*0.5) {\modFor};

\draw[->] (\baseDUDT+\wDUDT*0.5,\bblockSe) -- (\baseDUDT+\wDUDT*0.5,\tblock);

\draw[->] (\baseFCT+\wFCT+\resXoffs+1.,\yarrHI)--(\baseDUDT,\yarrHI);
\draw[->] (\baseFCT+\wFCT+\resXoffs+1.,\yarrLO)--(\baseDUDT,\yarrLO);

\draw[fill=\clrBgen,fill opacity=0.3] (\baseDUDT,\bblock) rectangle (\baseDUDT+\wDUDT,\tblock);
\node at (\baseDUDT+\wDUDT*0.5,1.75) {\modtoHBrhs};

\draw[->] (\baseDUDT+\wDUDT,\yarrHI) -- (\baseDUDT+\wDUDT+\resarrXoffs,\yarrHI);
\draw[->] (\baseDUDT+\wDUDT,\yarrLO) -- (\baseDUDT+\wDUDT+\resarrXoffs,\yarrLO);

\node at (\baseDUDT+\wDUDT+\rescXoffs+0.3,\yarrHI) {$\pat$};
\path (\baseDUDT+\wDUDT+\rescXoffs+0.6,\yarrHI-\yarrcROffs) pic {volcube=\clrH};
\node at (\baseDUDT+\wDUDT+\rescXoffs+0.6+\volnodeX,\yarrHI-\yarrcROffs+\volnodeY) {$\volH$};

\node at (\baseDUDT+\wDUDT+\rescXoffs+0.3,\yarrLO) {$\pat$};
\path (\baseDUDT+\wDUDT+\rescXoffs+0.6,\yarrLO-\yarrcROffs) pic {areacube=\clrB};
\node at (\baseDUDT+\wDUDT+\rescXoffs+0.6+\aXnodeX,\yarrLO-\yarrcROffs+\aXnodeY) {$\areaTbx$};
\node at (\baseDUDT+\wDUDT+\rescXoffs+0.6+\aYnodeX,\yarrLO-\yarrcROffs+\aYnodeY) {$\areaTby$};
\node at (\baseDUDT+\wDUDT+\rescXoffs+0.6+\aZnodeX,\yarrLO-\yarrcROffs+\aZnodeY) {$\areaTbz$};

\end{tikzpicture}
\caption{\label{fig:code_workflow}Right-hand-side computation's workflow.}
\end{figure}
\fi

\section{Numerical solver}
\label{sec:numsolver}
	The fourth-order solver presented in \cite{VTH19} is extended up to tenth order of accuracy. The main structure of the \rhs computation is shown in \fig{fig:code_workflow}: the solver reconstructs the fields at the interfaces between cells (\modHBreco, section \ref{sec:HBreco}), deduces the fluxes (\modHBflux, section \ref{sec:HBflux}) and from them the temporal derivatives (\modtoHBrhs, section \ref{sec:toHBrhs}). The computation of the non-ideal terms: explicit viscous and resistive terms, alongside the internal energy sink, is described in section \ref{sec:NID}. The driving of the system is discussed in section \ref{sec:forcing}. Finally, section \ref{sec:timeINT} presents the time integrator and section \ref{sec:numSummary} summarizes this method's keypoints.

\subsection{Reconstruction module (\modHBreco)}
\label{sec:HBreco}
This subsection describes the reconstruction module (\modHBreco in \fig{fig:code_workflow}), which is composed of two parts (\fig{fig:code_HBreco}): 
\begin{enumerate}
\item \modBinterp, which reconstructs volume averages of the magnetic field's components from the staggered constrained-transport discretization. This occurs through a polynomial interpolation. The magnetic field's volume averages are also needed to compute the internal energy sink (section \ref{sec:Esink}).
\item \modWENO, which stands for Weighted Essentially Non-Oscillatory. It is a reconstruction procedure computing area averages of all quantities at the cell's interfaces from their volume average.
\end{enumerate}
These two parts are described in the following.

\ifx\showtikz\undefined
\else
\def\xPc{1.0375}
\def\xMc{0.1625}
\def\lin{0.5}
\def\lout{0.2}
\def\lrout{0.25}
\def\baseBinterp{1.5+\lin}
\def\baseWENO{\baseBinterp+\bsize+\lout+1.5+\lin}
\def\xminBBlockHBreco{1.5+\lin*0.5}
\def\xmaxBBlockHBreco{\baseWENO+\bsize+\lout*0.5}
\def\yminBBlockHBreco{-0.2}
\def\ymaxBBlockHBreco{\yarrHI-\yarrcROffs+1.5+0.2}

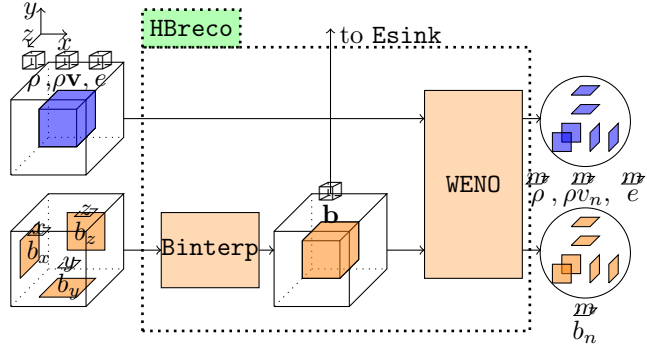
\begin{figure}[h]
\centering
\begin{tikzpicture}

\draw[dotted, line width=1pt] (\xminBBlockHBreco,\yminBBlockHBreco) rectangle (\xmaxBBlockHBreco,\ymaxBBlockHBreco);
\draw[dotted, line width=1pt,fill=\clrBgen,fill opacity=0.3] (\xminBBlockHBreco,\ymaxBBlockHBreco) rectangle (\xminBBlockHBreco+\bsize,\ymaxBBlockHBreco+0.5);
\node at (\xminBBlockHBreco+\bsize*0.5,\ymaxBBlockHBreco+0.25) {\modHBreco};

\path (0,\yarrHI-\yarrcROffs) pic {volcube=\clrH};

\node at (0+\volnodeX,\yarrHI-\yarrcROffs+\volnodeY) {$\volH$};

\path (0,\yarrLO-\yarrcROffs) pic {areacube=\clrB};

\node at (0+\aXnodeX,\yarrLO-\yarrcROffs+\aXnodeY) {$\areaTbx$};
\node at (0+\aYnodeX,\yarrLO-\yarrcROffs+\aYnodeY) {$\areaTby$};
\node at (0+\aZnodeX,\yarrLO-\yarrcROffs+\aZnodeY) {$\areaTbz$};

\path (0.4,3.75) pic[scale=0.35] {axes};

\draw[->] (1.5,\yarrHI)--(\baseWENO,\yarrHI);

\draw[->] (1.5,\yarrLO)--(\baseBinterp,\yarrLO);
\draw[fill=\clrBreco,fill opacity=0.3] (\baseBinterp,\yarrLO-0.5) rectangle (\baseBinterp+\bsize,\yarrLO+0.5);
\node at (\baseBinterp+\bsize*0.5,\yarrLO) {\modBinterp};
\draw[->] (\baseBinterp+\bsize,\yarrLO)--(\baseBinterp+\bsize+\lout,\yarrLO);
\path (\baseBinterp+\bsize+\lout,\yarrLO-0.75) pic {volcube=\clrB};
\node at (\baseBinterp+\bsize+\lout+\volnodeX,\yarrLO-0.75+\volnodeY) {$\vol{\Fmag}$};

\draw[->] (\baseBinterp+\bsize+\lout+1.5,\yarrLO)--(\baseWENO,\yarrLO);
\draw[fill=\clrBreco,fill opacity=0.3] (\baseWENO,\bblock) rectangle (\baseWENO+\bsize,\tblock);
\node at (\baseWENO+\bsize*0.5,1.75) {\modWENO};

\draw[->] (\baseBinterp+\bsize+\lout+0.75,\yarrLO-0.75+1.5)--(\baseBinterp+\bsize+\lout+0.75,\ymaxBBlockHBreco+0.25);
\node[anchor=north west] at (\baseBinterp+\bsize+\lout+0.75,\ymaxBBlockHBreco+0.25+0.15) {to \modEsink};

\draw[->] (\baseWENO+\bsize,\yarrHI) -- (\baseWENO+\bsize+\lrout,\yarrHI);
\draw[->] (\baseWENO+\bsize,\yarrLO) -- (\baseWENO+\bsize+\lrout,\yarrLO);

\path (\baseWENO+\bsize+\resXoffs,\yarrHI-\yarrROffs) pic {reco=\clrH};
\path (\baseWENO+\bsize+\resXoffs,\yarrLO-\yarrROffs) pic {reco=\clrB};

\node[anchor=north] at (\baseWENO+\bsize+\resXoffs+\underCircleX,\yarrHI-\yarrROffs-\underCircleY) {$\areaTH$};
\node[anchor=north] at (\baseWENO+\bsize+\resXoffs+\underCircleX,\yarrLO-\yarrROffs-\underCircleY) {$\areaTB$};

\end{tikzpicture}
\caption{\label{fig:code_HBreco}Reconstruction module.}
\end{figure}
\fi

\subsubsection{Magnetic field interpolation \modiBinterp}
\label{sec:Binterp}

The constrained-transport approach evolves each magnetic field component $\FmagL_n$ on the faces normal to its respective direction $\vn$. However, $\Fmag$ needs to be known on all faces in order to compute, e.g., the fluxes (\eq{eq:fluxexprx}). A possibility is to first deduce volume averages, which are subsequently treated in the same way as the hydrodynamic quantities $(\rho,\rho\Fvel)$ and the total energy density $\edens$ (see section \ref{sec:WENO}).

\newcommand{\orderindex}{m}
\newcommand{\porderindex}{q}
\newcommand{\itindex}{l}
\newcommand{\Nset}{\mathbb{N}}

\newcommand{\VFmagL}{B}
\newcommand{\BAfx}[1]{\area{\FmagL}{x}{\widetilde{(#1)}}{x,j,k}}
\newcommand{\BAfn}[1]{\area{\FmagL}{n}{\widetilde{(#1)}}{x,j,k}}

The volume averages are computed as follows. For a given discretization order $\orderindex=2\porderindex, \porderindex\in\Nset$, there is a unique polynomial $P$ of degree at most $\orderindex-1$ verifying:

\beq
P(x_\itindex-\Dx/2)=\AaSx{\FmagL}{x,\itindex-1/2,j,k}, \forall \itindex \in \{i-\porderindex+1,\ldots,i+\porderindex\}.
\eeq

Its line-average along $\vx$ gives a $\orderindex^{th}$ order approximation of the volume-averaged $\vol{\FmagL_x}$:

\beq
\frac{1}{\Dx} \int_{x_i-\Dx/2}^{x_i+\Dx/2} P(x)dx=\frac{1}{\Dx}\int_{x_i-\Dx/2}^{x_i+\Dx/2}\Big(\AaSx{\FmagL}{x}(x,y_j,z_k)+\bigO(\Dx^\orderindex)\Big)dx=\vol{\FmagL}_{x,i,j,k}+\bigO(\Dx^\orderindex).
\eeq

All algebra done, this method yields:

\beqa
\label{eq:magint2} \vol{\FmagL_x}_{i,j,k}&=&\frac{1}{2}\Big( \BAfx{\frac{1}{2}} \Big)+O(\Dx^2),\\
\label{eq:magint4}  &=&\frac{1}{24}\Big(13\BAfx{\frac{1}{2}}-\BAfx{\frac{3}{2}}\Big)+O(\Dx^4),\\
\label{eq:magint6} &=&\frac{1}{1440}\Big(802\BAfx{\frac{1}{2}}-93\BAfx{\frac{3}{2}}+11\BAfx{\frac{5}{2}}\Big)+O(\Dx^6),\\
\label{eq:magint8} &=&\frac{1}{120960}\Big(68323\BAfx{\frac{1}{2}}-9531\BAfx{\frac{3}{2}}+ 1879\BAfx{\frac{5}{2}}-191\BAfx{\frac{7}{2}} \Big)+O(\Dx^8),\\
\label{eq:magint10} &=&\frac{1}{7257600}\Big(4134338\BAfx{\frac{1}{2}}-641776\BAfx{\frac{3}{2}}  + 162680\BAfx{\frac{5}{2}}-28939\BAfx{\frac{7}{2}}+2497\BAfx{\frac{9}{2}}\Big)\!+\!O(\Dx^{10}),
\eeqa

with $\BAfx{p}=\AaSx{\FmagL}{x,i+p,j,k}+\AaSx{\FmagL}{x,i-p,j,k}$. The extension to the $\vy-$ and $\vz$-directions is straightforward.

Even though a non-oscillatory reconstruction would enhance the stability (cf. \cite{BMD13}), a polynomial interpolation is sufficient because each magnetic field component is continuous along its respective direction \cite{LOZ00}.

\subsubsection{Weighted Essentially Non-Oscillatory reconstruction \modiWENO}
\label{sec:WENO}

Weighted Essentially Non-Oscillatory (WENO) schemes allow high-order reconstruction in smooth regions while limiting oscillations in non-smooth ones. They are improvements of the Essentially Non-Oscillatory (ENO) schemes \cite{HEO87}. The first WENO scheme \cite{LOC94}, third-order accurate, appeared in the 90s and has rapidly been extended to fifth-order \cite{JIS96}. A method to derive WENO schemes of arbitrarily high order of accuracy is provided in \cite{SHU97}, and very high-order WENO schemes, with a discretization order up to 17 are given in \cite{BAS00,GSV09}.

\newcommand{\optP}{P^{opt}}
\newcommand{\vfa}{m}
\newcommand{\stenn}{n}
\newcommand{\stenP}{P}
\newcommand{\wdweno}{d}
\newcommand{\wwweno}{w}

The main idea of a WENO reconstruction is to split the reconstruction stencil in substencils associated with different weights \cite{SHU97}. The weights are chosen so that, on the one hand, maximum accuracy is obtained in smooth regions, and, on the other hand, weights associated with non-smooth regions are vanishingly small. The method is presented here for a reconstruction in the $\vx$-direction and can be straightforwardly transposed to the $\vy$-and $\vz$-directions.  Given a stencil $\{x_{i-\stenn},x_{i-\stenn+1},\ldots,x_i,\ldots,x_{i+\stenn}\}$ and a quantity $q$, there is a unique polynomial $\optP$ of degree at most $2\stenn$ such that:

\beq
\label{eq:popt}
\frac{1}{\Dx}\int_{x_{i+\vfa-1/2}}^{x_{i+\vfa+1/2}} \optP(x) \dx=\vol{q}_{i+\vfa,j,k}, \forall \vfa \in \{-\stenn,\ldots,\stenn\}.
\eeq

In smooth regions, this often-called ``optimal polynomial'' verifies:
\beq
\optP(x)=\AaSx{q}{}(x,y_j,z_k)+\bigO(\Dx^{2\stenn+1}).
\eeq

However, in regions containing strong gradients or even discontinuities, such as shock fronts, it would give an oscillatory reconstruction exhibiting over-and undershoots. To solve this issue, the $(2\stenn+1)$-long stencil is subdivided in the $\stenn+1$ convex sub-stencils containing $x_i$: $\{x_{i-n},\ldots,x_i\}$, $\{x_{i-n+1},\ldots,x_i,x_{i+1}\}$, $\ldots$, $\{x_i,\ldots,x_{i+n}\}$. On each substencil, there is a unique optimal polynomial $\stenP^{\vfa}$ of degree at most $\stenn$ which verifies \eq{eq:popt} for each cell in its substencil. At the cell boundaries, the aim is to write the ``global'' optimal polynomial $\optP$ as a linear combination of these ``local'' optimal polynomials $\stenP^{\vfa}$:

\beq
\label{eq:optPsumwdweno}
\optP(x_i\pm\frac{\Dx}{2})=\sum_\vfa \wdweno^\pm_{\vfa} \stenP^{\vfa}(x_i\pm\frac{\Dx}{2}),
\eeq

with two sets of positive weights $(\wdweno^\pm_{\vfa})$, $\sum \wdweno^\pm_{\vfa}=1$. 

Then, one designs other positive weights $(\wwweno^\pm_{\vfa})$ with $\sum_m \wwweno^\pm_{\vfa}=1$ such that, on the one hand, $(\wwweno^\pm_{\vfa}) \approx (\wdweno^\pm_{\vfa})$ in smooth regions, giving high-order accuracy, and, on the other hand, $\wwweno^\pm_{\vfa}\to 0$ in substencils containing discontinuities. A design possibility is:

\beq
\label{eq:wwweno}
\wwweno^\pm_{\vfa}=\frac{\alpha^\pm_\vfa}{\sum_q \alpha^\pm_q},\quad\quad\quad\alpha^\pm_\vfa=\frac{\wdweno^\pm_\vfa}{(\epsilon+IS_\vfa)^p},
\eeq

with $\epsilon=10^{-6}, p=2$ as in \cite{LPR99,VTH19} and $IS_\vfa$ is a so-called ``smoothness indicator'', defined as \cite{JIS96}:

\beq
	\label{eq:IS}
	IS^\vfa_i=\sum_{l=1}^{\stenn} \int_{x_i-\Dx/2}^{x_i+\Dx/2} (\Dx)^{2l-1} (\stenP_i^{\vfa,(l)}(x))^2 \dx,
\eeq

with $\stenP^{\vfa,(l)}$ the $l^{th}$ derivative of the $\stenP^{\vfa}$ polynomial. The reconstructed values are finally:

\beq
\label{eq:recoval}
\AaSx{q}{}(x_i\pm\Dx/2,y_j,z_k)=\sum_\vfa \wwweno^\pm_{\vfa} \stenP^{\vfa}(x_i\pm\frac{\Dx}{2})+\epsilon_{reco},
\eeq

with $\epsilon_{reco}$ the discretization error, which goes as $\bigO(\Dx^{2\stenn+1})$ in smooth regions.

For the sake of brevity, we do not repeat here the numerical expressions for the smoothness indicators $(IS^{\vfa})$, the weights $(\wdweno^\pm_{\vfa})$ and the reconstructed values $\AaSx{q}{i\pm 1/2,j,k}$ as a function of the volume averages $\vol{q}_{i,j,k}$, but refer to the literature \cite{JIS96,BAS00}.

\textbf{Remarks:}
\begin{itemize}
\item A reconstruction of the magnetic field $\FmagL_n$ component is not needed in the $\vn$-direction: this area-average is already known in the constrained-transport framework.
\item The procedure described here for volume-average$\to$area-average reconstruction is also used for area-average$\to$line-average reconstructions (section \ref{sec:CTm}). The principle is the same, replacing the superscripts $\vol{\cdot}$ by $\area{\cdot}{n}{}{}$ and $\area{\cdot}{m}{}{}$ by $\linea{\cdot}{m}{}$.
\end{itemize}

\subsection{Flux module \modiHBflux}
\label{sec:HBflux}

\ifx\showtikz\undefined
\else
\def\baseAtoP{\xPc+\lin}
\def\baseRiem{\baseAtoP+\bsize+\lout+\lin*2+\xPc}
\def\basePtoA{\baseRiem+\bsize+\lout+\lin+\xPc}

\def\yarrLLO{\yarrLO-1.75}
\def\yarrHHI{\yarrHI+1.75}

\def\hSblock{0.75}

\def\xminBBlockHBflux{\xPc+\lin*0.25}
\def\xmaxBBlockHBflux{\basePtoA+\bsize+\lout+0.5+\lout*0.5}
\def\yminBBlockHBflux{\yarrLLO-1.25}
\def\ymaxBBlockHBflux{\yarrHHI+0.75}

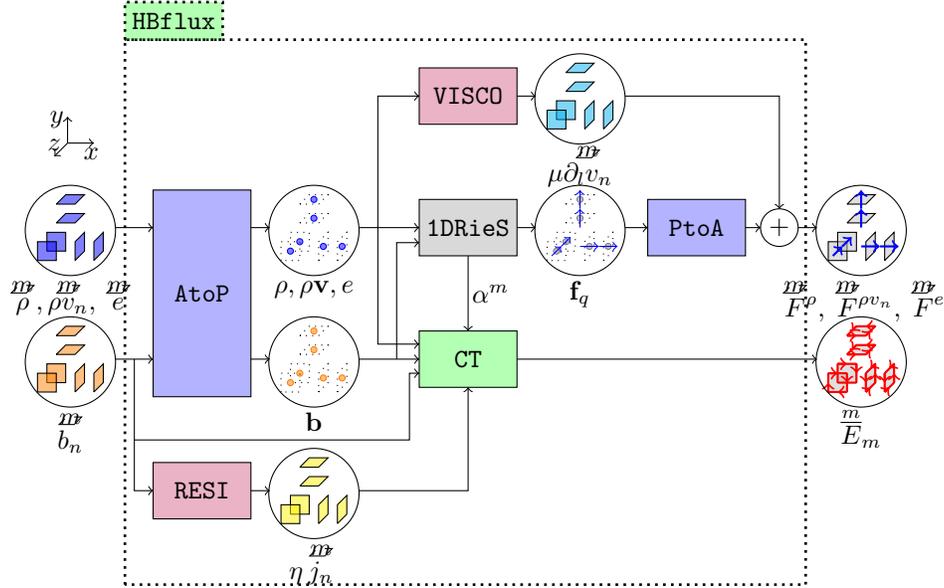
\begin{figure}[h]
\centering
\begin{tikzpicture}

\draw[dotted, line width=1pt] (\xminBBlockHBflux,\yminBBlockHBflux) rectangle (\xmaxBBlockHBflux,\ymaxBBlockHBflux);
\draw[dotted, line width=1pt,fill=\clrBgen,fill opacity=0.3] (\xminBBlockHBflux,\ymaxBBlockHBflux) rectangle (\xminBBlockHBflux+\bsize,\ymaxBBlockHBflux+0.5);
\node at (\xminBBlockHBflux+\bsize*0.5,\ymaxBBlockHBflux+0.25) {\modHBflux};

\path (0,\yarrHI-\yarrROffs) pic {reco=\clrH};

\path (0,\yarrLO-\yarrROffs) pic {reco=\clrB};

\node[anchor=north] at (0+\underCircleX,\yarrHI-\yarrROffs-\underCircleY) {$\areaTH$};
\node[anchor=north] at (0+\underCircleX,\yarrLO-\yarrROffs-\underCircleY) {$\areaTB$};

\path (0.4,3.75) pic[scale=0.35] {axes};

\draw[->] (1.0375,\yarrHI)--(\baseAtoP,\yarrHI);
\draw[->] (1.0375,\yarrLO)--(\baseAtoP,\yarrLO);

\draw[->] (1.0375+\lin*0.5,\yarrLO)--(1.0375+\lin*0.5,\yarrLLO) -- (\baseAtoP,\yarrLLO);

\draw[->] (1.0375+\lin*0.5,\yarrLLO*0.5-0.2)--(\baseRiem-\bsize*0.1,\yarrLLO*0.5-0.2) -- (\baseRiem-\bsize*0.1,\yarrLO-\hSblock*0.25) -- (\baseRiem,\yarrLO-\hSblock*0.25);

\draw[fill=\clrBcolella,fill opacity=0.3] (\baseAtoP,\yarrLO-0.5) rectangle (\baseAtoP+\bsize,\yarrHI+0.5);
\node at (\baseAtoP+\bsize*0.5,\yarrLO*0.5+\yarrHI*0.5) {\modAtoP};
\draw[->] (\baseAtoP+\bsize,\yarrHI)--(\baseAtoP+\bsize+\lrout,\yarrHI);
\path (\baseAtoP+\bsize+\resXoffs,\yarrHI-\yarrROffs) pic {pval=dotted/\clrH};
\draw[->] (\baseAtoP+\bsize,\yarrLO)--(\baseAtoP+\bsize+\lrout,\yarrLO);
\path (\baseAtoP+\bsize+\resXoffs,\yarrLO-\yarrROffs) pic {pval=dotted/\clrB};

\node[anchor=north] at (\baseAtoP+\bsize+\resXoffs+\underCircleX,\yarrHI-\yarrROffs-\underCircleYp) {$\rho,\rho\Fvel,\edens$};
\node[anchor=north] at (\baseAtoP+\bsize+\resXoffs+\underCircleX,\yarrLO-\yarrROffs-\underCircleYp) {$\Fmag$};

\draw[fill=\clrBNID,fill opacity=0.3] (\baseAtoP,\yarrLLO-\hSblock*0.5) rectangle (\baseAtoP+\bsize,\yarrLLO+\hSblock*0.5);
\node at (\baseAtoP+\bsize*0.5,\yarrLLO) {\modRESI};
\draw[->] (\baseAtoP+\bsize,\yarrLLO)--(\baseAtoP+\bsize+\lrout,\yarrLLO);
\path (\baseAtoP+\bsize+\resXoffs,\yarrLLO-\yarrROffs) pic {reco=\clrJ};

\node[anchor=north] at (\baseAtoP+\bsize+\resXoffs+\underCircleX,\yarrLLO-\yarrROffs-\underCircleY) {$\eta\areaTJ$};

\draw[->] (\baseAtoP+\bsize+\resXoffs+\xPc+\lin,\yarrLO)--(\baseAtoP+\bsize+\resXoffs+\xPc+\lin,\yarrHI-0.2)--(\baseRiem,\yarrHI-0.2);
\draw[->] (\baseAtoP+\bsize+\resXoffs+\xPc,\yarrHI)--(\baseRiem,\yarrHI);
\draw[fill=\clrBriemann,fill opacity=0.3] (\baseRiem,\yarrHI-\hSblock*0.5) rectangle (\baseRiem+\bsize,\yarrHI+\hSblock*0.5);
\node at (\baseRiem+\bsize*0.5,\yarrHI) {\modLLF};

\draw[->] (\baseRiem+\bsize,\yarrHI)--(\baseRiem+\bsize+\lrout,\yarrHI);
\path (\baseRiem+\bsize+\resXoffs,\yarrHI-\yarrROffs) pic {pflux=dotted/gray/\clrH};
\node[anchor=north] at (\baseRiem+\bsize+\resXoffs+\underCircleX,\yarrHI-\yarrROffs-\underCircleYp) {$\fluxfq$};

\draw[->] (\baseAtoP+\bsize+\resXoffs+\xPc+\lin*0.5,\yarrHI)--(\baseAtoP+\bsize+\resXoffs+\xPc+\lin*0.5,\yarrLO+0.2)--(\baseRiem,\yarrLO+0.2);
\draw[->] (\baseAtoP+\bsize+\resXoffs+\xPc,\yarrLLO)--(\baseRiem+\bsize*0.5,\yarrLLO)--(\baseRiem+\bsize*0.5,\yarrLO-\hSblock*0.5);
\draw[->] (\baseAtoP+\bsize+\resXoffs+\xPc,\yarrLO)--(\baseRiem,\yarrLO);
\draw[fill=\clrBgen,fill opacity=0.3] (\baseRiem,\yarrLO-\hSblock*0.5) rectangle (\baseRiem+\bsize,\yarrLO+\hSblock*0.5);
\node at (\baseRiem+\bsize*0.5,\yarrLO) {\modCT};

	\draw[->] (\baseRiem+\bsize*0.5,\yarrHI-\hSblock*0.5)--(\baseRiem+\bsize*0.5,\yarrLO+\hSblock*0.5);
	\node at (\baseRiem+\bsize*0.5+0.3,\yarrLO*0.5+\yarrHI*0.5) {$\alpha^{m}$};

\draw[->] (\baseAtoP+\bsize+\resXoffs+\xPc+\lin*0.5,\yarrHI)--(\baseAtoP+\bsize+\resXoffs+\xPc+\lin*0.5,\yarrHHI)--(\baseRiem,\yarrHHI);
\draw[fill=\clrBNID,fill opacity=0.3] (\baseRiem,\yarrHHI-\hSblock*0.5) rectangle (\baseRiem+\bsize,\yarrHHI+\hSblock*0.5);
\node at (\baseRiem+\bsize*0.5,\yarrHHI) {\modVISCO};

\draw[->] (\baseRiem+\bsize,\yarrHHI)--(\baseRiem+\bsize+\lrout,\yarrHHI);
\path (\baseRiem+\bsize+\resXoffs,\yarrHHI-\yarrROffs) pic {reco=\clrDV};
\node[anchor=north] at (\baseRiem+\bsize+\resXoffs+\underCircleX,\yarrHHI-\yarrROffs-\underCircleY+0.05) {$\mu\areaTDV$};

\draw[->] (\baseRiem+\bsize+\resXoffs+\xPc,\yarrHI)--(\basePtoA,\yarrHI);
\draw[fill=\clrBcolella,fill opacity=0.3] (\basePtoA,\yarrHI-\hSblock*0.5) rectangle (\basePtoA+\bsize,\yarrHI+\hSblock*0.5);
\node at (\basePtoA+\bsize*0.5,\yarrHI) {\modPtoA};


\draw[->] (\baseRiem+\bsize+\resXoffs+\xPc,\yarrHHI)--(\basePtoA+\bsize+\lout+0.25,\yarrHHI)--(\basePtoA+\bsize+\lout+0.25,\yarrHI+0.25);

\draw[->] (\basePtoA+\bsize,\yarrHI)--(\basePtoA+\bsize+\lout,\yarrHI);
\draw (\basePtoA+\bsize+\lout+0.25,\yarrHI) circle (0.25);
\node at (\basePtoA+\bsize+\lout+0.25,\yarrHI) {+};

\draw[->] (\basePtoA+\bsize+\lout+0.25+0.25,\yarrHI)--(\basePtoA+\bsize+\lout+0.25+0.25+\lrout,\yarrHI);
\path (\basePtoA+\bsize+\lout+0.25+0.25+\resXoffs,\yarrHI-\yarrROffs) pic {hflux=gray!50/\clrH};
\node[anchor=north] at (\basePtoA+\bsize+\lout+0.25+0.25+\resXoffs+\underCircleX,\yarrHI-\yarrROffs-\underCircleY) {$\areaTFH$};

\draw[->] (\baseRiem+\bsize,\yarrLO)--(\basePtoA+\bsize+\lout+0.25+0.25+\lrout,\yarrLO);
\path (\basePtoA+\bsize+\lout+0.25+0.25+\resXoffs,\yarrLO-\yarrROffs) pic {bflux=gray!50/\clrE};
\node[anchor=north] at (\basePtoA+\bsize+\lout+0.25+0.25+\resXoffs+\underCircleX,\yarrLO-\yarrROffs-\underCircleYlin) {$\lineTEm$};
\end{tikzpicture}
\caption{\label{fig:code_HBflux}Flux computation module.}
\end{figure}
\fi

The core of the flux computation module (\fig{fig:code_HBflux}) determines:

\begin{itemize}
\item the flux of the hydrodynamic quantities and the total energy density by the interfacial one-dimensional Riemann solver (\modLLF),
\item the magnetic field's flux by the constrained-transport module (\modCT).
\end{itemize}

These are described in sections \ref{sec:Rusanov} and \ref{sec:CTm} respectively.

In order to obtain an order of accuracy strictly greater than two, our approach requires a passage through point values (\modAtoP and \modPtoA), described in section \ref{sec:colella}.

Finally, the non-ideal contributions through viscosity (\modVISCO) and magnetic diffusivity (\modRESI) are discussed in sections \ref{sec:VISCO} and \ref{sec:RESI}.

\subsubsection{One-dimensional Riemann solver \modiLLF}
\label{sec:Rusanov}
\label{sec:LLF}

We employ the Rusanov approximation, also called ``local Lax-Friedrichs'' (LLF) \cite{RUS61,SHO89}, which explicitly considers only the fastest propagating wave mode on either side of a grid cell boundary (the fast magneto-sonic speed). This rough approximation results in additional numerical dissipation, but is chosen here both for the sake of simplicity and to underline the benefits of higher-order discretizations. At the boundary $(i+1/2,j,k)$, the LLF flux of field $\Qpt$ is given by:

\newcommand{\maxspeedS}{a}
\newcommand{\maxspeedxijk}{\maxspeedS^x_{i+1/2,j,k}}
\newcommand{\stateL}{\Qpt^-_{i+1/2,j,k}}
\newcommand{\stateR}{\Qpt^+_{i+1/2,j,k}}
\newcommand{\cfx}{c^x_f}
\newcommand{\cfi}{c^i_f}
\newcommand{\cA}{c_A}
\beq
	\label{eq:fluxLLF}
	\fluxLLFqijk=\frac{1}{2}( (\fluxfqIp+\fluxfqIm) \cdot \ex) - \frac{\maxspeedxijk}{2}(\stateR-\stateL),
\eeq

with the superscripts $\pm$ denoting the states ``right'' and ``left'' of the interface: $x\to (x_i+\Dx/2)^\pm$, reconstructed through the stencils centered at $(i+1,j,k)$ and $(i,j,k)$ respectively. The fluxes $\fluxfqIpm$ (see \eq{eq:fluxexprx}) are obtained by plugging the quantities at the corresponding side of the interface. The treatment of the viscous term is described in section \ref{sec:VISCO}. The maximum speed of propagation of information, $\maxspeedxijk$, is estimated as:

\beq
	\label{eq:maxspeed}
	\maxspeedxijk=\mathrm{max}\big( (|\Fvel_x|+\cfx)^+, (|\Fvel_x|+\cfx)^- \big),
\eeq

with $\cfx$ the fast magneto-sonic speed:

\beq
	\cfx=\sqrt{\frac{1}{2}\Big( \cs^2+\cA^2+ \sqrt{(\cs^2+\cA^2)^2-4\cs^2\frac{\FmagL^2_x}{\rho}} \Big)},
\eeq

where $\cs=(\gamma \press/\rho)^{1/2}$ and $\cA=(|\Fmag|^2/\rho)^{1/2}$ are the sound speed and the Alfv\'{e}n speed, respectively. The \formulas are analogous in the $\vy$-and $\vz$-directions.

The expressions above are strictly valid only for ideal MHD, i.e. for $\visc=\resis=0$. In practice one can still use them in non-ideal MHD when the physical dissipation is small enough: $\max(\resis,\visc/\rho)/\min(\Dx,\Dy,\Dz) \ll \min_{i\in\{x,y,z\}}(\cfi)$ (see section \ref{sec:ssturb}). Nevertheless, when solving for the dissipative terms alone, the Riemann solver has to be changed in order to have a stable scheme, see section \ref{sec:testNID}.

\subsubsection{Passage through point values \modiAtoPPtoA}
\label{sec:colella}

The reconstruction procedure (section \ref{sec:WENO}) provides area-averages of arbitrarily high order of accuracy. However, one cannot plug them directly in the flux \formulas (\eq{eq:fluxexprx}) without loss of numerical accuracy. Indeed, the fluxes consist of products and quotients of the fields $(\rho,\rho\Fvel,\edens,\Fmag)$, e.g.: 

\beq
\label{eq:ptflux}
\rho \FvelL_x\FvelL_y=(\rho \FvelL_x)(\rho \FvelL_y)/\rho.
\eeq

This equality is not valid for area-averages:

\beq
\label{eq:areaflux}
\area{(\rho \FvelL_x\FvelL_y)}{x}{}{} \neq \area{(\rho \FvelL_x)}{x}{}{}\area{(\rho \FvelL_y)}{x}{}{}/\area{\rho}{x}{}{}.
\eeq

Identifying the area-averaged quantity with its point-value in the middle of the considered face gives a second-order error term \cite{COC11,BUH14}. In the one-dimensional case, a Taylor expansion of $\linea{\Qpt}{x}{i}$, line-average of $\Qpt$ over cell $i$, illustrates this:

\beq
\label{eq:simple2O} \linea{\Qpt}{x}{i}=\frac{1}{\Dx} \int_{-\Dx/2}^{\Dx/2} q(x_i+\epsilon) \deps=\frac{1}{\Dx} \int_{-\Dx/2}^{\Dx/2} (q(x_i)+\epsilon q'(x_i) + \bigO(\epsilon^2))\deps=q_i+\bigO(\Dx^2).
\eeq

A work-around is proposed in \cite{COC11,BUH14}: 
\begin{enumerate}
\item the reconstructed area-averages are transformed into point values at the desired order of accuracy (\modAtoP in \fig{fig:code_HBflux}),
\item from them, the Riemann solver computes point-valued fluxes (\modLLF, section \ref{sec:Rusanov}),
\item the point-valued fluxes are finally transformed into area-averaged fluxes (\modPtoA). They have the same order of accuracy as the point values obtained in step 1.
\end{enumerate}

The method described in appendix \ref{app:transf_form_deriv} produces the area-average$\leftrightarrow$point value transformations written below. The area-average$\to$point value transformation \formulas are, up to tenth-order of accuracy:

\begin{eqnarray}
	\label{eq:ftop2}
	\QptW&=&\QavgW+\bigO(h^{2}),\\
	\label{eq:ftop4}
	\QptW&=&\frac{7}{6}\QavgW-\frac{1}{24}\QavgWu+\bigO(h^{4}),\\
	\label{eq:ftop6}
	\QptW&=&\frac{1771}{1440}\QavgW-\frac{23}{360}\QavgWu+\frac{3}{640}\QavgWuu+\frac{1}{576}\QavgWuv+\bigO(h^{6}),\\
\label{eq:ftop8} 	\QptW&=&\frac{25451}{20160}\QavgW-\frac{4973}{64512}\QavgWu+\frac{83}{8960}\QavgWuu-\frac{5}{7168}\QavgWuuu+\frac{19}{5760}\QavgWuv-\frac{1}{5120}\QavgWuuv+\bigO(h^{8}),\\
\nonumber 	\QptW&=&\frac{5514407}{4300800}\QavgW-\frac{69119}{806400}\QavgWu+\frac{42149}{3225600}\QavgWuu-\frac{55}{32256}\QavgWuuu\\
\label{eq:ftop10} &+&\frac{35}{294912}\QavgWuuuu+\frac{29173}{6451200}\QavgWuv-\frac{41}{89600}\QavgWuuv+\frac{5}{172032}\QavgWuuuv+\frac{9}{409600}\QavgWuuvv+\bigO(h^{10}),
\end{eqnarray}

\label{sec:explsigmamn}
where $h=\Dx+\Dy$ and $\area{q}{z}{\Sigma mn}{ij}, m\geq n\geq 0$ is the sum of all area-averages where one offset with respect to $(i,j)$ is $\pm m$ and the other $\pm n$. Each possible combination of offsets $(\pm m,\pm n),(\pm n,\pm m)$ is considered only once, thus, this sum contains eight terms, unless $n=0$ or $m=n$, in which case it contains four terms (see \fig{ref:fig_colstencils}). Thus, for $m>0$:

\newcommand{\QavgAZ}[1]{\area{q}{z}{}{#1}}

\beqa
\area{q}{z}{\Sigma mn}{ij}&\underset{m>n>0}{=}&\QavgAZ{i+m,j+n}+\QavgAZ{i+m,j-n}+\QavgAZ{i-m,j+n}+\QavgAZ{i-m,j-n}+\QavgAZ{i+n,j+m}+\QavgAZ{i+n,j-m}+\QavgAZ{i-n,j+m}+\QavgAZ{i-n,j-m},\\
\area{q}{z}{\Sigma mm}{ij}&=&\QavgAZ{i+m,j+m}+\QavgAZ{i+m,j-m}+\QavgAZ{i-m,j+m}+\QavgAZ{i-m,j-m},\\
\area{q}{z}{\Sigma m0}{ij}&=&\QavgAZ{i+m,j}+\QavgAZ{i-m,j}+\QavgAZ{i,j+m}+\QavgAZ{i,j-m}.
\eeqa

\ifx\showtikz\undefined
\else

\definecolor{clrWu}{rgb}{0.2,0.3,0.7}

\definecolor{clrWuu}{rgb}{0.99,0.7,0.5}
\definecolor{clrWuv}{rgb}{0.99,0.5,0.}

\definecolor{clrWuuu}{rgb}{0.9,0.5,0.9}
\definecolor{clrWuuv}{rgb}{0.6,0.,0.8}

\definecolor{clrWuuuu}{rgb}{0.85,0.99,0.99}
\definecolor{clrWuuvv}{rgb}{0.4,0.9,0.99}
\definecolor{clrWuuuv}{rgb}{0.,0.6,0.99}

\definecolor{grayA}{rgb}{0.2,0.2,0.2}
\definecolor{grayB}{rgb}{0.5,0.5,0.5}
\definecolor{grayC}{rgb}{0.8,0.8,0.8}

\definecolor{blueA}{rgb}{0.,0.,0.7}

\newcommand{\sAarg}{fill=black}

\def\gsize{0.75}
\def\gcent{4}

\def\gs{0.75}
\def\gc{4}

\newcommand{\squareS}[3]{(#1*\gsize,#2*\gsize)--(#1*\gsize+#3*\gsize,#2*\gsize)--(#1*\gsize+#3*\gsize,#2*\gsize+#3*\gsize)--(#1*\gsize,#2*\gsize+#3*\gsize)--cycle}
\newcommand{\squareT}[2]{(#1*\gsize+\gcent*\gsize,#2*\gsize+\gcent*\gsize)--(#1*\gsize+\gsize+\gcent*\gsize,#2*\gsize+\gcent*\gsize)--(#1*\gsize+\gsize+\gcent*\gsize,#2*\gsize+\gsize+\gcent*\gsize)--(#1*\gsize+\gcent*\gsize,#2*\gsize+\gsize+\gcent*\gsize)--cycle}
\newcommand{\squareC}[3]{\filldraw[#1] \squareT{#2}{#3}; \filldraw[#1] \squareT{-#2}{#3}; \filldraw[#1] \squareT{#2}{-#3}; \filldraw[#1] \squareT{-#2}{-#3}; \filldraw[#1] \squareT{#3}{#2}; \filldraw[#1] \squareT{-#3}{#2}; \filldraw[#1] \squareT{#3}{-#2}; \filldraw[#1] \squareT{-#3}{-#2};}

\newcommand{\atPos}[2]{at (#1*\gsize+0.5*\gsize,#2*\gsize+0.5*\gsize)}

\newcommand{\labelNode}[3]{\node at (#2*\gsize+0.5*\gsize+\gsize*\gcent,#3*\gsize+0.5*\gsize+\gsize*\gcent) {#1}; }
\newcommand{\squareD}[4]{\squareC{fill=#4}{#2}{#3}\labelNode{#1}{#2}{#3}\labelNode{#1}{-#2}{#3}\labelNode{#1}{#2}{-#3}\labelNode{#1}{-#2}{-#3}\labelNode{#1}{#3}{#2}\labelNode{#1}{-#3}{#2}\labelNode{#1}{#3}{-#2}\labelNode{#1}{-#3}{-#2}}

\newcommand{\argOlimIV}{ultra thick}
\newcommand{\argOlimVI}{ultra thick}

\newcommand{\legend}[3]{\filldraw[#1] \squareS{9.5}{#3}{0.5};}
\begin{figure}
\centering
\begin{tikzpicture}
\draw[step=\gsize,black,thin] (0.,0.) grid (9.*\gsize,9.*\gsize);

\node at (\gcent*\gsize+0.5*\gsize,\gcent*\gsize+0.5*\gsize) {$(i,j)$};
\squareD{$\supWu$}{1}{0}{clrWu}
\squareD{$\supWuu$}{2}{0}{clrWuu}
\squareD{$\supWuuu$}{3}{0}{clrWuuu}
\squareD{$\supWuuuu$}{4}{0}{clrWuuuu}
\squareD{$\supWuv$}{1}{1}{clrWuv}
\squareD{$\supWuuvv$}{2}{2}{clrWuuvv}
\squareD{$\supWuuv$}{2}{1}{clrWuuv}
\squareD{$\supWuuuv$}{3}{1}{clrWuuuv}

\def\facyaxis{0.25}
\node at (0*\gsize+0.5*\gsize,-\facyaxis*\gsize) {$i-4$};
\node at (1*\gsize+0.5*\gsize,-\facyaxis*\gsize) {$i-3$};
\node at (2*\gsize+0.5*\gsize,-\facyaxis*\gsize) {$i-2$};
\node at (3*\gsize+0.5*\gsize,-\facyaxis*\gsize) {$i-1$};
\node at (4*\gsize+0.5*\gsize,-\facyaxis*\gsize) {$i$};
\node at (5*\gsize+0.5*\gsize,-\facyaxis*\gsize) {$i+1$};
\node at (6*\gsize+0.5*\gsize,-\facyaxis*\gsize) {$i+2$};
\node at (7*\gsize+0.5*\gsize,-\facyaxis*\gsize) {$i+3$};
\node at (8*\gsize+0.5*\gsize,-\facyaxis*\gsize) {$i+4$};
\def\facxaxis{0.75}
\node[text width=1cm, align=right] at (-\facxaxis*\gsize,0*\gsize+0.5*\gsize) {$j-4$};
\node[text width=1cm, align=right] at (-\facxaxis*\gsize,1*\gsize+0.5*\gsize) {$j-3$};
\node[text width=1cm, align=right] at (-\facxaxis*\gsize,2*\gsize+0.5*\gsize) {$j-2$};
\node[text width=1cm, align=right] at (-\facxaxis*\gsize,3*\gsize+0.5*\gsize) {$j-1$};
\node[text width=1cm, align=right] at (-\facxaxis*\gsize,4*\gsize+0.5*\gsize) {$j$};
\node[text width=1cm, align=right] at (-\facxaxis*\gsize,5*\gsize+0.5*\gsize) {$j+1$};
\node[text width=1cm, align=right] at (-\facxaxis*\gsize,6*\gsize+0.5*\gsize) {$j+2$};
\node[text width=1cm, align=right] at (-\facxaxis*\gsize,7*\gsize+0.5*\gsize) {$j+3$};
\node[text width=1cm, align=right] at (-\facxaxis*\gsize,8*\gsize+0.5*\gsize) {$j+4$};

\end{tikzpicture}
\caption{\label{ref:fig_colstencils}Illustration of the different terms appearing in the area-averages$\lra$point value transformations. A term with a superscript of the form $\Sigma mn$ with $m,n\in[0,4]$ corresponds to the sum of the values in the cells where this superscript appears: it contains either four or eight terms.}
\end{figure}

Using the same notation for the point values: $\Qpt^{\Sigma mn}$, the point-to-area transformations needed for the transformation of the point-valued fluxes (\modPtoA block) are given by:
\fi

\begin{eqnarray}
	\label{eq:ptof2}
	\QavgW&=&\QptW+\bigO(h^{2}),\\
	\label{eq:ptof4}
	\QavgW&=&\frac{5}{6}\QptW+\frac{1}{24}\QptWu+\bigO(h^{4}),\\
	\label{eq:ptof6}
	\QavgW&=&\frac{1159}{1440}\QptW+\frac{1}{20}\QptWu-\frac{17}{5760}\QptWuu+\frac{1}{576}\QptWuv+\bigO(h^{6}),\\
\label{eq:ptof8} 	\QavgW&=&\frac{47939}{60480}\QptW+\frac{52223}{967680}\QptWu-\frac{241}{48384}\QptWuu+\frac{367}{967680}\QptWuuu+\frac{47}{17280}\QptWuv-\frac{17}{138240}\QptWuuv+\bigO(h^{8}),\\
\nonumber 	\QavgW&=&\frac{91251029}{116121600}\QptW+\frac{163519}{2903040}\QptWu-\frac{186331}{29030400}\QptWuu+\frac{12011}{14515200}\QptWuuu\\
\label{eq:ptof10} &-&\frac{27859}{464486400}\QptWuuuu+\frac{193537}{58060800}\QptWuv-\frac{3667}{14515200}\QptWuuv+\frac{367}{23224320}\QptWuuuv+\frac{289}{33177600}\QptWuuvv+\bigO(h^{10}).
\end{eqnarray}

The \formulas above, for area-averages normal to $\vz$, are generalizable for the other directions in a straightforward way.

From these area-averages$\lra$point values transformations, one can deduce line-averages$\lra$point values transformations, which are needed when solving two-dimensional problems, see appendix \ref{app:transf_form_lp}. When applying an internal energy sink (section \ref{sec:Esink}) or some driving to the system (section \ref{sec:forcing}), one needs however volume-averages$\lra$point values transformations. These \formulas are given in appendix \ref{app:transf_form_vp}.

\subsubsection{Constrained-transport module \modiCT}
\label{sec:CTm}

\newcommand{\idE}{E^{ideal}}
\newcommand{\idEz}{E^{ideal}_z}

\ifx\showtikz\undefined
\else
\def\yarrDoubleCOffs{0.9375*\gscdc}
\def\xHDC{1.875*\gscdc}
\def\baseI{\xHDC+\lin}
\def\baseII{\baseI+\bsize+\lout+\xHDC+\lout+0.5+\lout}
\def\baseIII{\baseII+\bsize+\lout+\xHDC+\lin}

\def\CTplus{\baseI+\bsize+\lout+\xHDC+\lout+0.25}

\def\xminBBlockCT{\xHDC+\lin*0.5}
\def\xmaxBBlockCT{\baseIII+\bsize+\lout*0.5}
\def\yminBBlockCT{\yarrLO-\xHDC*0.5-0.15}
\def\ymaxBBlockCT{\yarrHI+0.25+\xHDC*0.5+0.75}

\def\XdcT{\xHDC*0.5}
\def\YdcT{\xHDC+0.2}
\def\XdcR{\xHDC+0.2}
\def\YdcR{\xHDC*0.5}

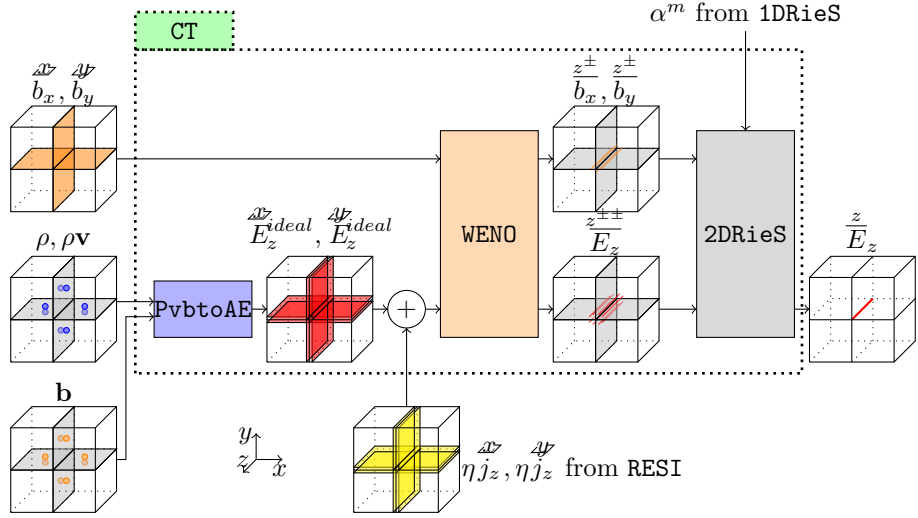
\begin{figure}
\centering
\begin{tikzpicture}
\draw[dotted, line width=1pt] (\xminBBlockCT,\yminBBlockCT) rectangle (\xmaxBBlockCT,\ymaxBBlockCT);
\draw[dotted, line width=1pt,fill=\clrBgen,fill opacity=0.3] (\xminBBlockCT,\ymaxBBlockCT) rectangle (\xminBBlockCT+\bsize,\ymaxBBlockCT+0.5);
\node at (\xminBBlockCT+\bsize*0.5,\ymaxBBlockCT+0.25) {\modCT};

	\path (0,\yarrLO-\yarrDoubleCOffs) pic {pdcube=\clrH};
	\node at (0+\XdcT,\yarrLO-\yarrDoubleCOffs+\YdcT) {$\rho,\rho\Fvel$};

	\draw[->] (\xHDC,\yarrLO+0.1) -- (\xHDC+\lin,\yarrLO+0.1);
	
	\path (0,\yarrLLO-0.25-\yarrDoubleCOffs) pic {pdcube=\clrB};
	\node at (0+\XdcT,\yarrLLO-0.25-\yarrDoubleCOffs+\YdcT) {$\Fmag$};
	
	\draw[->] (\xHDC,\yarrLLO-0.25) -- (\xHDC+\lin*0.25,\yarrLLO-0.25)-- (\xHDC+\lin*0.25,\yarrLO-0.1)-- (\xHDC+\lin,\yarrLO-0.1);
	
	\path (0,\yarrHI+0.25-\yarrDoubleCOffs) pic {Adcube=\clrB};
	\node at (0+\XdcT,\yarrHI+0.25-\yarrDoubleCOffs+\YdcT+0.1) {$\areaTbx,\areaTby$};
	\draw[->] (\xHDC,\yarrHI+0.25) -- (\baseII,\yarrHI+0.25);

\draw[fill=\clrBcolella,fill opacity=0.3] (\baseI,\yarrLO-\hSblock*0.5) rectangle (\baseI+\bsize,\yarrLO+\hSblock*0.5);
\node at (\baseI+\bsize*0.5,\yarrLO) {\modPvbtoAE};
	\draw[->] (\baseI+\bsize,\yarrLO) -- (\baseI+\bsize+\lout,\yarrLO);
	\path (\baseI+\bsize+\lout,\yarrLO-\yarrDoubleCOffs) pic {Ardcube=\clrE};
	\node at (\baseI+\bsize+\lout+\XdcT,\yarrLO-\yarrDoubleCOffs+\YdcT+0.2) {$\areaT{E}{x}{ideal}{z},\areaT{E}{y}{ideal}{z}$};

\draw[->] (\baseI+\bsize+\lout+\xHDC,\yarrLO)--(\baseI+\bsize+\lout+\xHDC+\lout,\yarrLO);
\draw (\baseI+\bsize+\lout+\xHDC+\lout+0.25,\yarrLO) circle (0.25);
\node at (\CTplus,\yarrLO) {+};

\draw[->] (\CTplus,\yarrLLO-0.25+\yarrDoubleCOffs)--(\CTplus,\yarrLO-0.25);
	\path (\CTplus-\xHDC*0.5,\yarrLLO-0.25-\yarrDoubleCOffs) pic {Ardcube=\clrJ};
	\node at (\CTplus-\xHDC*0.5+\XdcR+1.3,\yarrLLO-0.25-\yarrDoubleCOffs+\YdcR) {$\eta\areaT{j_z}{x}{}{},\eta\areaT{j_z}{y}{}{}$ from \modRESI};

\path (\CTplus-2.,\yarrLLO-0.25-\yarrDoubleCOffs+\YdcR) pic[scale=0.35] {axes};

\draw[->] (\baseI+\bsize+\lout+\xHDC+\lout+0.5,\yarrLO)--(\baseII,\yarrLO);

\draw[fill=\clrBreco,fill opacity=0.3] (\baseII,\yarrLO-\hSblock*0.5) rectangle (\baseII+\bsize,\yarrHI+0.25+\hSblock*0.5);
\node at (\baseII+\bsize*0.5,\yarrLO*0.5+\yarrHI*0.5+0.125) {\modWENO};
	\draw[->] (\baseII+\bsize,\yarrLO) -- (\baseII+\bsize+\lout,\yarrLO);
	\path (\baseII+\bsize+\lout,\yarrLO-\yarrDoubleCOffs) pic {lldcube=\clrE};
	\node at (\baseII+\bsize+\lout+\XdcT,\yarrLO-\yarrDoubleCOffs+\YdcT+0.1) {$\EfLApq{z^{\pm\pm}}{z}$};

	\draw[->] (\baseII+\bsize,\yarrHI+0.25) -- (\baseII+\bsize+\lout,\yarrHI+0.25);
	\path (\baseII+\bsize+\lout,\yarrHI+0.25-\yarrDoubleCOffs) pic {bldcube=\clrB};
	\node at (\baseII+\bsize+\lout+\XdcT,\yarrHI+0.25-\yarrDoubleCOffs+\YdcT+0.1) {$\bfLApq{z^\pm}{x},\bfLApq{z^\pm}{y}$};

\draw[->] (\baseIII+\bsize*0.5,\yarrHI+1.95)--(\baseIII+\bsize*0.5,\yarrHI+0.25+\hSblock*0.5);
\node at (\baseIII+\bsize*0.5,\yarrHI+2.2) {$\alpha^m$ from \modLLF};
\draw[->] (\baseII+\bsize+\lout+\xHDC,\yarrLO)--(\baseIII,\yarrLO);
\draw[->] (\baseII+\bsize+\lout+\xHDC,\yarrHI+0.25)--(\baseIII,\yarrHI+0.25);
\draw[fill=\clrBriemann,fill opacity=0.3] (\baseIII,\yarrLO-\hSblock*0.5) rectangle (\baseIII+\bsize,\yarrHI+0.25+\hSblock*0.5);
\node at (\baseIII+\bsize*0.5,\yarrLO*0.5+\yarrHI*0.5+0.125) {\modMLLF};
	\draw[->] (\baseIII+\bsize,\yarrLO) -- (\baseIII+\bsize+\lout,\yarrLO);
	\path (\baseIII+\bsize+\lout,\yarrLO-\yarrDoubleCOffs) pic {Edcube=\clrE};
	\node at (\baseIII+\bsize+\lout+\XdcT,\yarrLO-\yarrDoubleCOffs+\YdcT+0.1) {$\EfLA{z}{z}$};
	
\end{tikzpicture}
\caption{\label{fig:code_CT}Constrained-transport module.}
\end{figure}
\fi

The key element in the \modCT module (\fig{fig:code_CT}) is the two-dimensional Riemann solver (\modMLLF), which computes the line-averaged electric field required to update the magnetic field components (\eq{eq:dtmagCT}). The method is presented here for the $z$ component: $\EfLA{z}{z}$ and is analogous for $\EfLA{x}{x}$ and $\EfLA{y}{y}$ with appropriate permutations. In this work, the line-averaged electric field is approximated by a multidimensional version of the LLF approximation \cite{BAL10}, for the sake of simplicity. Improvements are available e.g. in \cite{BNK17}. At a certain edge of a numerical grid cell, e.g. $(i-1/2,j-1/2,k)$:

\beq
\label{eq:ELflux}
\EfLA{z}{z} = \frac{1}{4}\Big[\EfLApq{z^{++}}{z}+\EfLApq{z^{+-}}{z}+\EfLApq{z^{-+}}{z}+\EfLApq{z^{--}}{z}\Big]+\frac{S}{2}\left(\bfLApq{z^+}{y}-\bfLApq{z^-}{y}\right)-\frac{S}{2}\left(\bfLApq{z^+}{x}-\bfLApq{z^-}{x}\right),
\eeq

where the positional subscripts are omitted in order to limit the amount of notation: all are at $(i-1/2,j-1/2,k)$. The first term is the average of four line-averaged electric field, with $\EfLApq{z^{mn}}{z}, m,n\in\{+,-\}$ meaning the value obtained through reconstructions at the position $x\to (x_i-\Dx/2)^m, y\to(y_j-\Dy/2)^n$. This reconstruction procedure is described below. The two other terms add artificial dissipation. The maximum speed of propagation of information, $S$, is estimated as:

\beq
\label{eq:CTmaxspeed}
S=\max(a^x_{i-1/2,j,k},a^x_{i-1/2,j-1,k},a^y_{i,j-1/2,k},a^y_{i-1,j-1/2,k}),
\eeq

where $a^x$ and $a^y$ are given by \eq{eq:maxspeed}. This choice is made for the sake of efficiency \cite{LOZ04}, as these have already been computed in the \modLLF block. One should ideally consider the maximum speed among the four reconstructed states.

The magnetic field line-averages, $\bfLApq{z^{\pm}}{x}$ and $\bfLApq{z^{\pm}}{y}$, are obtained through WENO reconstruction of the area-averages (section \ref{sec:WENO}). These area-averages are the ones used in the constrained-transport discretization, so that they do not need the passage through the blocks \modBinterp and \modWENO in the \modHBreco module. A reconstruction along $\vy$ of $\area{\FmagL}{x}{}{x}$ gives $\bfLApq{z^{\pm}}{x}$, where $\pm$ means $y\to(y_j-\Dy/2)^\pm$. Similarly, a reconstruction along $\vx$ of $\area{\FmagL}{y}{}{y}$ gives $\bfLApq{z^{\pm}}{y}$, for $x\to(x_i-\Dx/2)^\pm$. 

The four line-averaged electric field in \eq{eq:ELflux} are obtained through the following steps (here described only for the $z$-component):
\begin{enumerate}
\item Compute point values of the ideal (i.e. magnetic diffusivity $\resis=0$) electric field: $\idE=-\Fvel \times \Fmag$ from the point values of $\rho, \rho \Fvel$ and $\Fmag$ obtained from the \modAtoP block of the \modHBflux module.
\item Use a point value$\to$area-average transformation (section \ref{sec:colella}) to deduce eight area-averages of $\idEz$ (two for each face surrounding e.g. the $(i-1/2,j-1/2,k)$ edge, one on each side). This step and the previous one are included in the \modPvbtoAE block (\fig{fig:code_CT}).  
\item If needed ($\resis>0$), add the non-ideal part (see section \ref{sec:RESI}) to obtain eight area-averages of the total $E_z$.
\item Perform eight WENO reconstructions of the area-averaged $E_z$: along $\vy$ for the faces normal to $\vx$ and along $\vx$ for the faces normal to $\vy$. This gives eight line-averaged values of $\EfLA{z}{z}$: two for each cell adjacent to the considered edge.
\item Assuming that the discontinuities occur only at cell boundaries, the two reconstructed values in each cell (one from the volume-averages$\to$area-averages reconstruction along $\vx$ followed by an area-average$\to$line-average reconstruction along $\vy$, and the other from volume-averages$\to$area-averages along $\vy$ followed by area-average$\to$line-average along $\vx$) are merged into one by taking their mean. This gives one state per cell: $\EfLApq{z^{++}}{z}$ from $\Cell_{i,j,k}$, $\EfLApq{z^{+-}}{z}$ from $\Cell_{i,j-1,k}$, $\EfLApq{z^{-+}}{z}$ from $\Cell_{i-1,j,k}$ and $\EfLApq{z^{--}}{z}$ from $\Cell_{i-1,j-1,k}$. These four states are used in the two-dimensional Riemann solver \modMLLF. 
\end{enumerate}

\subsection{Deduction of the \rhs\ \moditoHBrhs}
\label{sec:toHBrhs}

\ifx\showtikz\undefined
\else
\def\xPc{1.0375}
\def\xMc{0.1625}
\def\lin{0.5}
\def\lout{0.2}
\def\lrout{0.25}
\def\baseBsumF{1.5+\lin}

\def\baseSumK{\baseBsumF+\bsize+\lout*2+0.25}
\def\baseSumM{\baseBsumF+\bsize+\lout*2+1.5+\lout*2+0.25}

\def\xminBBlocktoHBrhs{1.5+\lin*0.5}
\def\xmaxBBlocktoHBrhs{\baseSumM+0.5+\lout*2}
\def\yminBBlocktoHBrhs{-0.2}
\def\ymaxBBlocktoHBrhs{\yarrHI-\yarrcROffs+1.5+0.2}

\def\yFOR{\yarrLO-1.25}
\def\yESINK{\yarrHI+1.25}

\begin{figure}[h]
\centering
\begin{tikzpicture}

\draw[dotted, line width=1pt] (\xminBBlocktoHBrhs,\yminBBlocktoHBrhs) rectangle (\xmaxBBlocktoHBrhs,\ymaxBBlocktoHBrhs);
\draw[dotted, line width=1pt,fill=\clrBgen,fill opacity=0.3] (\xminBBlocktoHBrhs,\ymaxBBlocktoHBrhs) rectangle (\xminBBlocktoHBrhs+\bsize,\ymaxBBlocktoHBrhs+0.5);
\node at (\xminBBlocktoHBrhs+\bsize*0.5,\ymaxBBlocktoHBrhs+0.25) {\modtoHBrhs};

\path (1.5-1.0375,\yarrHI-\yarrROffs) pic {hflux=gray!50/\clrH};
\node[anchor=north] at (1.5-1.0375+\underCircleX-0.05,\yarrHI-\yarrROffs-\underCircleY) {$\areaTFH$};
\path (1.5-1.0375,\yarrLO-\yarrROffs) pic {bflux=gray!50/\clrE};
\node[anchor=north] at (1.5-1.0375+\underCircleX,\yarrLO-\yarrROffs-\underCircleYlin) {$\lineTEm$};

\draw[->] (1.5,\yarrLO)--(\baseBsumF,\yarrLO);
\draw[->] (1.5,\yarrHI)--(\baseBsumF,\yarrHI);
\draw[fill=\clrBcalc,fill opacity=0.3] (\baseBsumF,\yarrLO-0.5) rectangle (\baseBsumF+\bsize,\yarrHI+0.5);
\node at (\baseBsumF+\bsize*0.5,\yarrLO*0.5+\yarrHI*0.5) {\modSumF};

\draw[->] (\baseBsumF+\bsize,\yarrHI)--(\baseSumK-0.25,\yarrHI);
\draw[->] (\baseBsumF+\bsize,\yarrLO)--(\baseSumM-0.25,\yarrLO);

\draw (\baseSumK,\yarrHI) circle (0.25);
\node at (\baseSumK,\yarrHI) {+};

\path (\baseSumK-0.75,\yESINK) pic {volcube=\clrU!50};
\node at (\baseSumK+1.75-0.25+0.1,\yESINK+0.75) {\begin{tabular}{c}$\vol{\Usink}$ \\ from \modEsink \end{tabular}};

\draw[->] (\baseSumK,\yESINK)--(\baseSumK,\yarrHI+0.25);

\path (\baseSumK+3.,\yESINK+0.75) pic[scale=0.35] {axes};

\draw[->] (\baseSumK,\yFOR)--(\baseSumK,\yarrHI-0.25);
\path (\baseSumK-0.75,\yFOR-1.5) pic {volcube=\clrForce!50};
\node at (\baseSumK-1.75+0.1-0.15,\yFOR-0.75) {\begin{tabular}{c}$\vol{\Kforce}, \vol{\KMforceE}$ \\ from \modFor\end{tabular}};


\draw (\baseSumM,\yarrLO) circle (0.25);
\node at (\baseSumM,\yarrLO) {+};

\draw[->] (\baseSumM,\yFOR)--(\baseSumM,\yarrLO-0.25);

\path (\baseSumM-0.75,\yFOR-1.5) pic {areacube=\clrForce!50};
\node at (\baseSumM+1.65+0.15,\yFOR-0.75) {\begin{tabular}{c}$\areaT{\Mforce}{n}{}{}$ \\ from \modFor\end{tabular}};


\draw[->] (\baseSumM+0.25,\yarrLO) -- (\xmaxBBlocktoHBrhs+\lout,\yarrLO);
\draw[->] (\baseSumK+0.25,\yarrHI) -- (\xmaxBBlocktoHBrhs+\lout,\yarrHI);

\node at (\xmaxBBlocktoHBrhs+\lout+0.3,\yarrHI) {$\pat$};
\path (\xmaxBBlocktoHBrhs+\lout+0.6,\yarrHI-\yarrcROffs) pic {volcube=\clrH};
\node at (\xmaxBBlocktoHBrhs+\lout+0.6+\volnodeX,\yarrHI-\yarrcROffs+\volnodeY) {$\volH$};

\node at (\xmaxBBlocktoHBrhs+\lout+0.3,\yarrLO) {$\pat$};
\path (\xmaxBBlocktoHBrhs+\lout+0.6,\yarrLO-\yarrcROffs) pic {areacube=\clrB};
\node at (\xmaxBBlocktoHBrhs+\lout+0.6+\aXnodeX,\yarrLO-\yarrcROffs+\aXnodeY) {$\areaTbx$};
\node at (\xmaxBBlocktoHBrhs+\lout+0.6+\aYnodeX,\yarrLO-\yarrcROffs+\aYnodeY) {$\areaTby$};
\node at (\xmaxBBlocktoHBrhs+\lout+0.6+\aZnodeX,\yarrLO-\yarrcROffs+\aZnodeY) {$\areaTbz$};

\end{tikzpicture}
\caption{\label{fig:code_toHBrhs}Deducing right-hand-side.}
\end{figure}
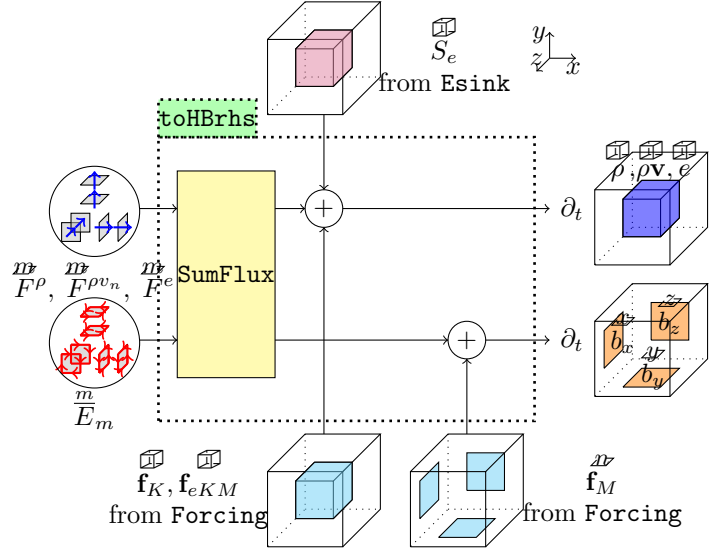
\fi

The \rhs of the MHD equations \eqsN{eq:beginMHD}{eq:dtb} is determined by (see \fig{fig:code_toHBrhs}, which is the block \modtoHBrhs in \fig{fig:code_workflow}):

\begin{itemize}
	\item plugging the interfacial fluxes (\eq{eq:fluxLLF}, after passage through the \modPtoA block, \fig{fig:code_HBflux}) into \eq{eq:patQavgT},
	\item plugging the line-averaged electric fields (\eq{eq:ELflux}) into \eq{eq:dtmagCT}. This step and the above one are represented by the \modSumF block in \fig{fig:code_toHBrhs},
	\item adding the internal energy sink term (see section \ref{sec:Esink}),
	\item adding the forcing terms for the momentum, magnetic field and total energy density (see section \ref{sec:forcing}).
\end{itemize}

\subsection{Treatment of non-ideal terms}
\label{sec:NID}

\ifx\showtikz\undefined
\else
\def\baseNIDI{1.5+\lin}
\def\baseNIDII{\baseNIDI+\bsize+\lout+1.5+\lin}
\def\baseNIDP{1.5-\xPc-\bsize-\xMc-\lrout}

\def\XacB{0.4375}
\def\YacB{-0.2-0.2}

\def\XvcB{0.75}
\def\YvcB{-0.2}

\def\xminBBlockVISCO{\baseNIDP-\lin*0.5}
\def\xmaxBBlockVISCO{\baseNIDII+\bsize+\lout*0.5}
\def\yminBBlockVISCO{\yarrVVHI-1.5}
\def\ymaxBBlockVISCO{\yarrVVHI+1.}

\def\xminBBlockRESI{\baseNIDI-\lin*0.5}
\def\xmaxBBlockRESI{\baseNIDII+\bsize+\lout*0.5}
\def\yminBBlockRESI{\yarrVHI-1.5}
\def\ymaxBBlockRESI{\yarrVHI+1.}

\def\xminBBlockESINK{\baseNIDI-\lin*0.5}
\def\xmaxBBlockESINK{\baseNIDII+\bsize+\lout*0.5}
\def\yminBBlockESINK{\yarrLO-1.3}
\def\ymaxBBlockESINK{\yarrLO+1.}

\def\offYA{-0.2}
\def\offYB{-0.4}

\begin{figure}[h]
\centering
\begin{tikzpicture}
\draw[dotted, line width=1pt] (\xminBBlockVISCO,\yminBBlockVISCO) rectangle (\xmaxBBlockVISCO,\ymaxBBlockVISCO);
\draw[dotted, line width=1pt,fill=\clrBNID,fill opacity=0.3] (\xminBBlockVISCO,\ymaxBBlockVISCO) rectangle (\xminBBlockVISCO+\bsize,\ymaxBBlockVISCO+0.5);
\node at (\xminBBlockVISCO+\bsize*0.5,\ymaxBBlockVISCO+0.25) {\modVISCO};

\draw[dotted, line width=1pt] (\xminBBlockRESI,\yminBBlockRESI+\offYA) rectangle (\xmaxBBlockRESI,\ymaxBBlockRESI+\offYA);
\draw[dotted, line width=1pt,fill=\clrBNID,fill opacity=0.3] (\xminBBlockRESI,\ymaxBBlockRESI+\offYA) rectangle (\xminBBlockRESI+\bsize,\ymaxBBlockRESI+0.5+\offYA);
\node at (\xminBBlockRESI+\bsize*0.5,\ymaxBBlockRESI+0.25+\offYA) {\modRESI};

\draw[dotted, line width=1pt] (\xminBBlockESINK,\yminBBlockESINK+\offYB) rectangle (\xmaxBBlockESINK,\ymaxBBlockESINK+\offYB);
\draw[dotted, line width=1pt,fill=\clrBNID,fill opacity=0.3] (\xminBBlockESINK,\ymaxBBlockESINK+\offYB) rectangle (\xminBBlockESINK+\bsize,\ymaxBBlockESINK+0.5+\offYB);
\node at (\xminBBlockESINK+\bsize*0.5,\ymaxBBlockESINK+0.25+\offYB) {\modEsink};

\draw[fill=\clrBcolella,fill opacity=0.3] (\baseNIDP,\yarrVVHI-\hSblock*0.5) rectangle (\baseNIDP+\bsize,\yarrVVHI+\hSblock*0.5);

\draw[->] (\baseNIDP+\bsize,\yarrVVHI)--(\baseNIDP+\bsize+\lrout,\yarrVVHI);
\node at (\baseNIDP+\bsize*0.5,\yarrVVHI) {\modPtoAv};
\draw[->] (\baseNIDP-\lin,\yarrVVHI)--(\baseNIDP,\yarrVVHI);

	\path (\baseNIDP-\lin-\xPc,\yarrVVHI-\yarrROffs) pic {pval=dotted/\clrH};
	\node[anchor=north] at (\baseNIDP-\lin-\xPc+\underCircleX,\yarrVVHI-\yarrROffs-\underCircleYp) {$\rho,\rho\Fvel$};


	\path (1.5-\xPc,\yarrVVHI-\yarrROffs) pic {reco=\clrV};
	\node[anchor=north] at (1.5-\xPc+\underCircleX,\yarrVVHI-\yarrROffs-\underCircleY) {$\areaT{v_n}{m}{}{}$};
	
	\path (1.5-\xPc,\yarrVHI-\yarrROffs+\offYA) pic {reco=\clrB};
	\node[anchor=north] at (1.5-\xPc+\underCircleX,\yarrVHI-\yarrROffs-\underCircleY+\offYA) {$\areaT{b_n}{m}{}{}$};

\path (-1.75,\yarrLO-\yarrcROffs+\offYB) pic {volcube=\clrH};
\node[anchor=north] at (-1.75+\XvcB,\yarrLO-\yarrcROffs-\underCubeY+\offYB) {$\vol{\rho},\vol{\rho\Fvel},\vol{\edens}$};
\path (0,\yarrLO-\yarrcROffs+\offYB) pic {volcube=\clrB};
\node[anchor=north] at (\XvcB,\yarrLO-\yarrcROffs-\underCubeY+\offYB) {$\vol{\Fmag}$};

\path (-1.,\yarrVHI+\offYA) pic[scale=0.35] {axes};

\draw[->] (1.5,\yarrVVHI)--(\baseNIDI,\yarrVVHI);
\draw[->] (1.5,\yarrVHI+\offYA)--(\baseNIDI,\yarrVHI+\offYA);
\draw[->] (1.5,\yarrLO+\offYB)--(\baseNIDI,\yarrLO+\offYB);

\draw[fill=\clrBDVDB,fill opacity=0.3] (\baseNIDI,\yarrVVHI-\hSblock*0.5) rectangle (\baseNIDI+\bsize,\yarrVVHI+\hSblock*0.5);
\node at (\baseNIDI+\bsize*0.5,\yarrVVHI) {\modAvtoVdv};

\draw[fill=\clrBDVDB,fill opacity=0.3] (\baseNIDI,\yarrVHI-\hSblock*0.5+\offYA) rectangle (\baseNIDI+\bsize,\yarrVHI+\hSblock*0.5+\offYA);
\node at (\baseNIDI+\bsize*0.5,\yarrVHI+\offYA) {\modAbtoVj};

\draw[fill=\clrBcolella,fill opacity=0.3] (\baseNIDI,\yarrLO-\hSblock*0.5+\offYB) rectangle (\baseNIDI+\bsize,\yarrLO+\hSblock*0.5+\offYB);
\node at (\baseNIDI+\bsize*0.5,\yarrLO+\offYB) {\modVtoP};

\draw[->] (\baseNIDI+\bsize,\yarrVVHI)--(\baseNIDI+\bsize+\lout,\yarrVVHI);
\draw[->] (\baseNIDI+\bsize,\yarrVHI+\offYA)--(\baseNIDI+\bsize+\lout,\yarrVHI+\offYA);
\draw[->] (\baseNIDI+\bsize,\yarrLO+\offYB)--(\baseNIDI+\bsize+\lout,\yarrLO+\offYB);

\path (\baseNIDI+\bsize+\lout,\yarrVVHI-\yarrcROffs) pic {volcube=\clrDV};
\node[anchor=north] at (\baseNIDI+\bsize+\lout+\XvcB,\yarrVVHI-\yarrcROffs-\underCubeY) {$\vol{\partial_{l}\FvelL_n}{}$};

\path (\baseNIDI+\bsize+\lout,\yarrVHI-\yarrcROffs+\offYA) pic {volcube=\clrJ};
\node[anchor=north] at (\baseNIDI+\bsize+\lout+\XvcB,\yarrVHI-\yarrcROffs-\underCubeY+\offYA) {$\vol{\FmagJ}$};

\path (\baseNIDI+\bsize+\lout,\yarrLO-\yarrcROffs+\offYB) pic {pvalvol=\clrU};
\node[anchor=north] at (\baseNIDI+\bsize+\lout+\XvcB,\yarrLO-\yarrcROffs-\underCubeY+\offYB) {$(\rho,\rho\Fvel,\edens,\Fmag)$};


\draw[->] (\baseNIDI+\bsize+\lout+1.5,\yarrVVHI)--(\baseNIDII,\yarrVVHI);
\draw[->] (\baseNIDI+\bsize+\lout+1.5,\yarrVHI+\offYA)--(\baseNIDII,\yarrVHI+\offYA);
\draw[->] (\baseNIDI+\bsize+\lout+1.5,\yarrLO+\offYB)--(\baseNIDII,\yarrLO+\offYB);

\draw[fill=\clrBreco,fill opacity=0.3] (\baseNIDII,\yarrVVHI-\hSblock*0.5) rectangle (\baseNIDII+\bsize,\yarrVVHI+\hSblock*0.5);
\node at (\baseNIDII+\bsize*0.5,\yarrVVHI) {\modWENO};

\draw[fill=\clrBreco,fill opacity=0.3] (\baseNIDII,\yarrVHI-\hSblock*0.5+\offYA) rectangle (\baseNIDII+\bsize,\yarrVHI+\hSblock*0.5+\offYA);
\node at (\baseNIDII+\bsize*0.5,\yarrVHI+\offYA) {\modWENO};

\draw[fill=\clrBcolella,fill opacity=0.3] (\baseNIDII,\yarrLO-\hSblock*0.5+\offYB) rectangle (\baseNIDII+\bsize,\yarrLO+\hSblock*0.5+\offYB);
\node at (\baseNIDII+\bsize*0.5,\yarrLO+\offYB) {\modPtoVSe};

\draw[->] (\baseNIDII+\bsize,\yarrVVHI)--(\baseNIDII+\bsize+\lrout,\yarrVVHI);
\draw[->] (\baseNIDII+\bsize,\yarrVHI+\offYA)--(\baseNIDII+\bsize+\lrout,\yarrVHI+\offYA);
\draw[->] (\baseNIDII+\bsize,\yarrLO+\offYB)--(\baseNIDII+\bsize+\lout,\yarrLO+\offYB);

\path (\baseNIDII+\bsize+\resXoffs,\yarrVVHI-\yarrROffs) pic {reco=\clrDV};
\node[anchor=north] at (\baseNIDII+\bsize+\resXoffs+\underCircleX,\yarrVVHI-\yarrROffs-\underCircleY) {$\mu\areaTDV$};

\path (\baseNIDII+\bsize+\resXoffs,\yarrVHI-\yarrROffs+\offYA) pic {reco=\clrJ};
\node[anchor=north] at (\baseNIDII+\bsize+\resXoffs+\underCircleX,\yarrVHI-\yarrROffs-\underCircleY+\offYA) {$\eta\areaTJ$};

\path (\baseNIDII+\bsize+\lout,\yarrLO-\yarrcROffs+\offYB) pic {volcube=\clrU};
\node[anchor=north] at (\baseNIDII+\bsize+\resXoffs+\XacB,\yarrLO-\yarrcROffs-\underCubeY+\offYB) {$\vol{\Usink}$};

\end{tikzpicture}
\caption{\label{fig:code_NID}Modules for non-ideal terms.}
\end{figure}
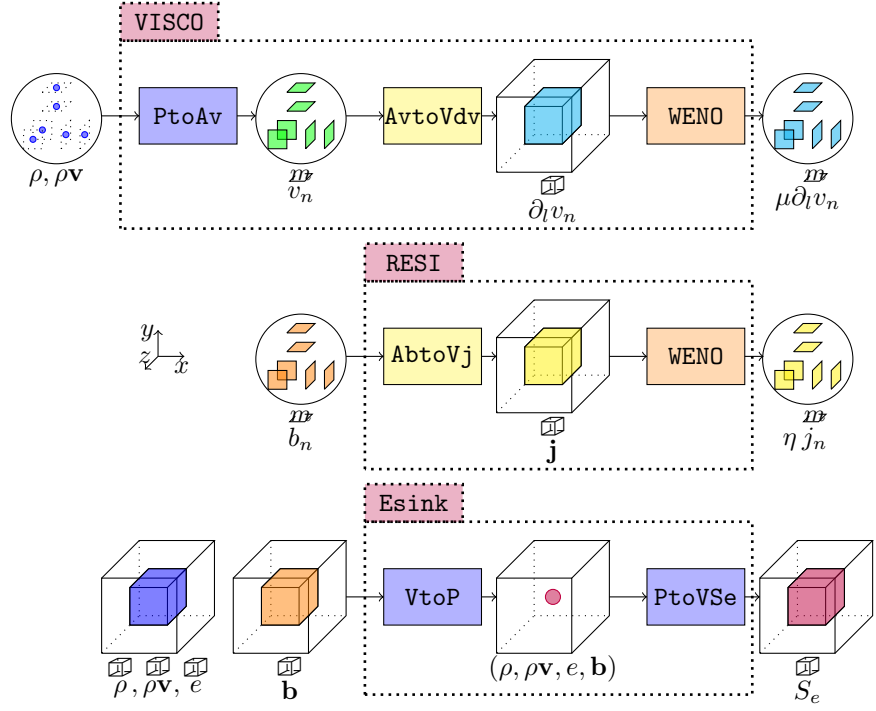
\fi

The dissipative terms are computed as sketched in \fig{fig:code_NID}. The next three subparts describe the treatment of the viscosity (\modVISCO), magnetic diffusivity (\modRESI) and internal energy sink (\modEsink), in this order.

\subsubsection{Viscosity \modiVISCO}
\label{sec:VISCO}
The viscous contribution (\eq{eq:divviscstress}) to the hydrodynamic fluxes is computed through the following steps:

\begin{enumerate}
\item From the point values in the middle of each interface between cells (after block \modAtoP, \fig{fig:code_HBflux}): compute velocity field point values and deduce its area-averages through a point value$\to$area-average transformation (section \ref{sec:colella}, block \modPtoAv in \fig{fig:code_NID}).
\item Deduce volume averages of the derivatives $\pal \FvelL_n$ by noticing that:
\beq
\label{eq:volpaxvn}(\vol{\pax v_n})_{i,j,k}=\frac{1}{\Dx}\left(\frac{1}{\Dy\Dz}\iint dydz \int \pax v_n dx\right)=\frac{\AaSx{v}{n,i+1/2,j,k}-\AaSx{v}{n,i-1/2,j,k}}{\Dx},
\eeq
and similarly for the $\vy$-and $\vz$-directions (block \modAvtoVdv in \fig{fig:code_NID}). Since \eq{eq:volpaxvn} is analytically exact, the scheme's discretization order is preserved.
\item Reconstruct area-averages of $\pal \FvelL_n$ through a WENO procedure (section \ref{sec:WENO}). Four reconstructions are needed for the fluxes in each dimension: in dimension $l$, one needs $\pal v_n, n\in\{x,y,z\}$ and $\nabla \cdot \Fvel$.
\item Multiply by the viscosity $\mu$.
\item Since the LLF flux is used in this solver (section \ref{sec:LLF}), it is enough for each interface to consider the average of the two reconstructed states surrounding it. This average is added to the computed area-averaged fluxes (see \fig{fig:code_HBflux}). When using a different Riemann solver that does not have the same coefficients in front of the ``left'' and ``right'' states, one would need to plug these area-averages in the \modAtoP block (\fig{fig:code_HBflux}) and provide them to the Riemann solver.
\end{enumerate}

\subsubsection{Magnetic diffusivity \modiRESI}
\label{sec:RESI}

The procedure for the magnetic diffusivity is very similar to that for the viscosity. The volume-averages of $\partial_l b_n$ with $l \neq n$ are computed from the area-averages (cf. \eq{eq:volpaxvn}), from which the volume averaged current $\FmagJ=\nabla \times \Fmag$ is deduced (\modAbtoVj block, \fig{fig:code_NID}), for example:

\beq
	\vol{\FmagJL_x}=\frac{\AaSy{\FmagL}{z,i,j+1/2,k}-\AaSy{\FmagL}{z,i,j-1/2,k}}{\Dy}-\frac{\AaSz{\FmagL}{y,i,j,k+1/2}-\AaSz{\FmagL}{y,i,j,k-1/2}}{\Dz}.
\eeq

After multiplication by $\eta$, a WENO procedure gives the area average of the electric field's non-ideal part. This term is added to the area-averaged electric field's ideal part (\fig{fig:code_CT}).

\subsubsection{Internal energy sink \modiEsink}
\label{sec:Esink}

As pointed out in section \ref{sec:colella}, a passage through point values allows to build schemes of order strictly greater than two. Thus, the internal energy sink is computed through the following steps:

\begin{enumerate}
\item Transform the volume averages of $(\rho,\rho\Fvel,\edens,\Fmag)$ to point values in the middle of the volume (\modVtoP in \fig{fig:code_NID}).
\item Deduce the internal energy $\intE=p/(\gadia-1)$ (cf. \eq{eq:pressure}) and thus the Stefan-Boltzman like radiative loss term $\Usink$ (\eq{eq:Usink}) as a point value.
\item Transform $\Usink$ to $\vol{\Usink}$ through a point value$\to$volume average transformation (\modPtoVSe).
\item Add this contribution to the \rhs of $\edens$ (cf. \fig{fig:code_workflow}).
\end{enumerate}

The \formulas for the point value$\lra$volume average transformations are given in appendix \ref{app:transf_form_vp}.

\subsection{Forcing terms}
\label{sec:forcing}

\newcommand{\fOUK}{{\bf f}^{OU,K}}
\newcommand{\fOUM}{{\bf f}^{OU,M}}
\newcommand{\fOUML}{f^{OU,M}}

\newcommand{\fOUkF}{\hat{\bf f}^{OU,K}}
\newcommand{\fOUbF}{\hat{\bf f}^{OU,M}}
\newcommand{\ftOUbF}{\tilde{\bf f}^{OU,M}}
\newcommand{\fOUFk}{\hat{\bf f}^{OU,m}_{\vk}}
\newcommand{\fat}{t_{auto}}
\newcommand{\famp}{F_0}
\newcommand{\specw}{\zeta} 
\newcommand{\projk}{{\cal{P}^{\specw}}} 
\newcommand{\projkij}{{{\cal{P}}^{\specw}_{ij}}} 
\newcommand{\projkT}{\underline{\projk}} 

\newcommand{\fAK}{A_K}
\newcommand{\fAM}{A_M}

\ifx\showtikz\undefined
\else
\def\xPc{1.0375}
\def\xMc{0.1625}
\def\lin{0.2}
\def\lout{0.2}
\def\lrout{0.25}
\def\baseBsumF{1.5+\lin}

\def\baseSumK{\baseBsumF+\bsize+\lout*2+0.25}
\def\baseSumM{\baseBsumF+\bsize+\lout*2+1.5+\lout*2+0.25}

\def\yFOR{\yarrLO-1.25}
\def\yESINK{\yarrHI+1.25}

\def\baseUpdateF{1.5}
\def\baseUpdateR{\baseUpdateF+\bsize+\lout}
\def\baseFEMK{\baseUpdateR+\bsize+\lout+1.5+\lout}

\def\baseVforce{\baseFEMK+\bsize+\lout+1.5+\lout}

\def\baseOUT{\baseVforce+\bsize+\lout*5+1.5}

\def\hCblocks{3.}
\def\yarrTsep{1.25}
\def\yarrDsep{0.8}

\def\yarrLfalse{\yarrLO*0.+\yarrLLO-0.25}
\def\yarrLmid{2.625*0.5-0.875*0.5}
\def\yinterL{0.1}
\def\xminBBlockForce{\baseUpdateF-\lin-0.6-0.3}
\def\xmaxBBlockForce{\baseVforce+\bsize+\lout*0.5}
\def\yminBBlockForce{\yarrHI-\yarrTsep-0.75-0.25}
\def\ymaxBBlockForce{\yarrHI+\hCblocks*0.5+1.1}

\begin{figure}
\centering
\begin{tikzpicture}

\draw[dotted, line width=1pt] (\xminBBlockForce,\yminBBlockForce) rectangle (\xmaxBBlockForce,\ymaxBBlockForce);
\draw[dotted, line width=1pt,fill=\clrBforce,fill opacity=0.3] (\xminBBlockForce,\ymaxBBlockForce) rectangle (\xminBBlockForce+\bsize,\ymaxBBlockForce+0.5);
\node at (\xminBBlockForce+\bsize*0.5,\ymaxBBlockForce+0.25) {\modFor};

\path (\baseUpdateF-\lin-0.3,\yarrHI+\yarrDsep) pic {Fdot=\clrForce};
\path (\baseUpdateF-\lin-0.3,\yarrHI-\yarrDsep) pic {Fdot=\clrForce};

\node at (\baseUpdateF-\lin-0.3-0.05,\yarrHI+\yarrDsep-0.45) {$\fOUkF$};
\node at (\baseUpdateF-\lin-0.3-0.05,\yarrHI-\yarrDsep-0.45) {$\fOUbF$};
\draw[->] (\baseUpdateF-\lin,\yarrHI+\yarrDsep)--(\baseUpdateF,\yarrHI+\yarrDsep);
\draw[->] (\baseUpdateF-\lin,\yarrHI-\yarrDsep)--(\baseUpdateF,\yarrHI-\yarrDsep);

\draw[fill=\clrBcalc,fill opacity=0.3] (\baseUpdateF,\yarrHI-\hCblocks*0.5) rectangle (\baseUpdateF+\bsize,\yarrHI+\hCblocks*0.5);
\node at (\baseUpdateF+\bsize*0.5,\yarrHI) {\modUpdateF};

\draw[->] (\baseUpdateF+\bsize,\yarrHI)--(\baseUpdateR,\yarrHI);
\draw[fill=\clrBcalc,fill opacity=0.3] (\baseUpdateR,\yarrHI-\hCblocks*0.5) rectangle (\baseUpdateR+\bsize,\yarrHI+\hCblocks*0.5);
\node at (\baseUpdateR+\bsize*0.5,\yarrHI) {\modUpdateR};

\draw[->] (\baseUpdateR+\bsize,\yarrHI-\yarrcROffs+0.75) -- (\baseUpdateR+\bsize+\lout,\yarrHI-\yarrcROffs+0.75);
\draw[->] (\baseUpdateR+\bsize,\yarrHI+\yarrTsep-\yarrcROffs+0.75) -- (\baseUpdateR+\bsize+\lout,\yarrHI+\yarrTsep-\yarrcROffs+0.75);

\path (\baseUpdateR+\bsize+\lout,\yarrHI-\yarrcROffs) pic {pvalvol=\clrForce};
\path (\baseUpdateR+\bsize+\lout,\yarrHI+\yarrTsep-\yarrcROffs) pic {pvalvol=\clrForce};

\node at (\baseUpdateR+\bsize+\lout+0.75,\yarrHI-\yarrcROffs-0.2) {$\fAM\fOUM$};
\node at (\baseUpdateR+\bsize+\lout+0.75,\yarrHI+\yarrTsep-\yarrcROffs+1.7) {$\Kforce$};

\draw[->] (\baseUpdateR+\bsize+\lout+1.5,\yarrHI) -- (\baseFEMK,\yarrHI);
\draw[->] (\baseUpdateR+\bsize+\lout+1.5+\lout*0.5,\yarrHI+\yarrTsep) -- (\baseUpdateR+\bsize+\lout+1.5+\lout*0.5,\yarrHI+0.5*\yarrTsep) -- (\baseFEMK,\yarrHI+0.5*\yarrTsep);

\draw[fill=\clrBcalc,fill opacity=0.3] (\baseFEMK,\yarrHI-\hCblocks*0.25) rectangle (\baseFEMK+\bsize,\yarrHI+\hCblocks*0.25);
\node at (\baseFEMK+\bsize*0.5,\yarrHI+\hCblocks*0.) {\modFEMK};
\draw[->] (\baseFEMK+\bsize,\yarrHI) -- (\baseFEMK+\bsize+\lout,\yarrHI);
\path (\baseFEMK+\bsize+\lout,\yarrHI-\yarrcROffs) pic {pvalvol=\clrForce};

\node at (\baseFEMK+\bsize+\lout+0.75,\yarrHI-\yarrcROffs-0.2) {$\KMforceE$};

\draw[->] (\baseFEMK+\bsize+\lout+1.5,\yarrHI) -- (\baseVforce,\yarrHI);
\draw[->] (\baseUpdateR+\bsize+\lout+1.5,\yarrHI+\yarrTsep) -- (\baseVforce,\yarrHI+\yarrTsep);

\draw[fill=\clrBcolella,fill opacity=0.3] (\baseVforce,\yarrHI-\hCblocks*0.25) rectangle (\baseVforce+\bsize,\yarrHI+\hCblocks*0.5);
\node at (\baseVforce+\bsize*0.5,\yarrHI+\hCblocks*0.125) {\modVforce};

\draw[->] (\baseVforce+\bsize,\yarrHI)--(\baseVforce+\bsize+\lout,\yarrHI);
\draw[->] (\baseVforce+\bsize,\yarrHI+\yarrTsep)--(\baseVforce+\bsize+\lout,\yarrHI+\yarrTsep);
\path (\baseVforce+\bsize+\lout,\yarrHI+\yarrTsep-\yarrcROffs) pic {volcube=\clrForce};
\path (\baseVforce+\bsize+\lout,\yarrHI-\yarrcROffs) pic {volcube=\clrForce};

\draw[->] (\baseUpdateR+\bsize,\yarrHI-\yarrTsep-\yarrcROffs+0.75) -- (\baseVforce+\bsize+\lout,\yarrHI-\yarrTsep-\yarrcROffs+0.75);
\path (\baseVforce+\bsize+\lout,\yarrHI-\yarrTsep-\yarrcROffs) pic {areacube=\clrForce};

\node at (\baseVforce+\bsize+\lout+0.75,\yarrHI+\yarrTsep+1.1) {$\vol{\Kforce}$};
\node[anchor=west] at (\baseVforce+\bsize+\lout+1.,\yarrHI) {$\vol{\KMforceE}$};
\node[anchor=north] at (\baseVforce+\bsize+\lout+0.75,\yarrHI-\yarrTsep-\yarrcROffs) {$\areaT{f}{n}{}{\!M,n}$};

\end{tikzpicture}
\caption{\label{fig:forcing}Forcing module.}
\end{figure}
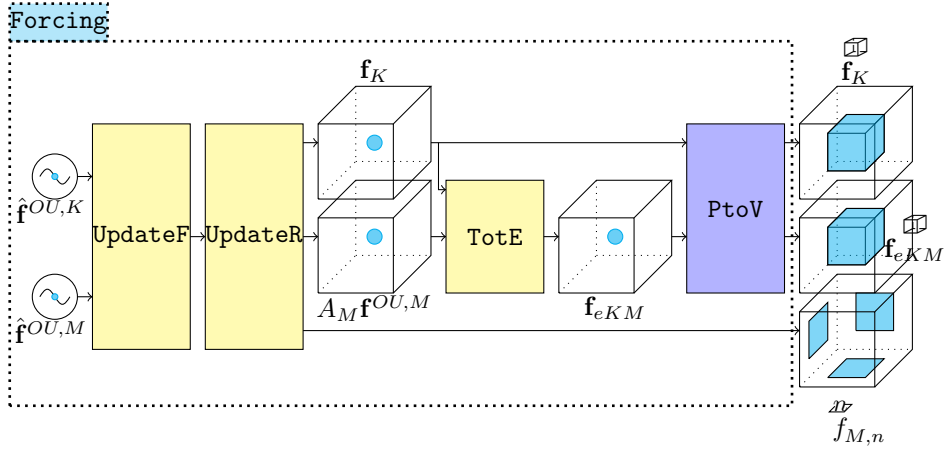
\fi

Section \ref{sec:ssturb} presents an example of application: a turbulent statistically stationary state. It results from a dynamical balance between large-scale energy injection and dissipative effects, both of physical and numerical nature. The mechanical and electromotive drivings ($\Kforce$ and $\Mforce$ in \eqsa{eq:dtrhov}{eq:dtb}) inject kinetic and magnetic energy at the largest scales, which is transported through \nonlinear effects down to the smallest scales, where dissipative effects dominate. The $\Kforce$ and $\Mforce$ forcing terms affect the total energy equation through the term $\KMforceE$ in \eq{eq:dtedens}. These three terms are computed through the following steps (see \fig{fig:forcing}), explained in the next subsections:

\begin{itemize}
\item update the mechanical and electromotive drivings in Fourier space (block \modUpdateF in \fig{fig:forcing}, section \ref{sec:forcing_Fourier_update}),
\item update and normalize them in real space to achieve a desired energy injection rate (block \modUpdateR, section \ref{sec:forcing_real_update}),
\item transform them to respect the chosen discretization: finite-volume using the constrained-transport approach (blocks \modFEMK and \modVforce, section \ref{sec:forcing_apply}).
\end{itemize}

Lastly, some remarks relevant in the context of turbulent MHD simulations are given in section \ref{sec:forcing_remarks}.

\subsubsection{Update in spectral space}
\label{sec:forcing_Fourier_update}

Both kinetic and magnetic energies are injected through an Ornstein-Uhlenbeck process. The driving fields $\fOUkF$ (kinetic) and $\fOUbF$ (magnetic), set to zero initially, are both evolved in Fourier space according to the stochastic differential equation: 

\newcommand{\OUspecprof}{\theta}

\beqa
                \label{eq:dtfOUF}
                d\fOUFk(t)=-\fOUFk(t)\frac{dt}{\fat}+\famp \left( \frac{2\OUspecprof(\vk)^2}{\fat} \right)^{1/2} \projkT \cdot d\vW(t),
\eeqa

with $\fat$ the forcing \autocorrelation time, $\famp$ an amplitude which value can be taken arbitrarily because of the normalization procedure described below and $\OUspecprof(\vk)$ a spectral profile explicited in section \ref{sec:ssturb}. A Wiener process $d\vW(t)=dt\vN(0,dt)$ modelizes a three-dimensional continuous random walk, with $\vN(0,dt)$ a 3D Gaussian distribution with zero mean and standard deviation $dt$. The random numbers drawn are different for the kinetic and magnetic forcing, so that $\fOUbF\neq\fOUkF$. Purely solenoidal forcing fields are obtained by spectral projection $\projkij(\vk)=\specw \delta_{ij}+(1-2\specw) \frac{k_ik_j}{|k|^2}$ with $\specw=1$ and the three wavenumbers $k_i$ of the corresponding wave vector.

\subsubsection{Update and normalization in configuration space}
\label{sec:forcing_real_update}
The $\fOUkF$ and $\fOUbF$ complex fields are transformed to configuration space using a fast-Fourier transform algorithm, giving respectively $\fOUK_{i,j,k}$ and $\fOUM_{i,j,k}$, the values of the driving fields at the point $(x_i,y_j,z_k)$. The point-valued electromotive driving $\fOUM_{i,j,k}$ is not used to force the magnetic field (which is defined as staggered area-averages), but to update the total energy density $\edens$, see section \ref{sec:forcing_apply}. In order to get the area-averaged electromotive driving, the Fourier coefficients $\fOUbF$ are modified by (cf. \eq{eq:defBx}):

\beq
\label{eq:ftOUbF}
\ftOUbF_{x,\vk}=\fOUbF_{x,\vk}\cdot\underbrace{e^{-ik_x\frac{\Dx}{2}}}_{\text{Staggered field}}\cdot\overbrace{\frac{(e^{ik_y\frac{\Dy}{2}}-e^{ik_y\frac{\Dy}{2}})(e^{ik_z\frac{\Dz}{2}}-e^{ik_z\frac{\Dz}{2}})}{\Dy\Dz(ik_y)(ik_z)}}^{\text{Area average}},
\eeq

and similarly for the $\vy$-and $\vz$-components. After transformation to configuration space, the corresponding $\area{f}{x}{OU,M}{x,i-1/2,j,k}$, $\area{f}{y}{OU,M}{y,i,j-1/2,k}$ and $\area{f}{z}{OU,M}{z,i,j,k-1/2}$ can be applied to the staggered magnetic field without introducing non-solenoidal contributions.

\newcommand{\Ekin}{\mathcal{E}_{K}}
\newcommand{\Emag}{\mathcal{E}_{M}}

\newcommand{\EKinj}{\epsilon^K_{inj}}
\newcommand{\EMinj}{\epsilon^M_{inj}}

In order to obtain prescribed constant kinetic and magnetic energy injection rates ($\EKinj$ and $\EMinj$, respectively), the fields in configuration space are multiplied with two normalization factors $\fAK$ and $\fAM$. These are taken as the largest root of the following second-order polynomial equations, which estimate in a first-order way the amount of energy injected during a time-interval $\Dt$, if one would be using the Euler time integrator (cf. \cite{MCL99}):

\newcommand{\Stag}{\boldsymbol{G}}

\beqa
        \label{eq:polyeinjK}\Delta \Ekin &=&\frac{1}{2} \fAK^2 \Dt^2 \sum_{i,j,k} \vol{\rho}_{i,j,k} (\fOUK_{i,j,k})^2+\fAK\Dt \sum_{i,j,k} \vol{(\rho\Fvel)}_{i,j,k} \cdot \fOUK_{i,j,k},\\
        \nonumber\Delta \Emag &=&\frac{1}{2} \fAM^2 \Dt^2 \sum_{i,j,k} \Stag(f^{OU,M})_{i,j,k}^2+\fAM\Dt \sum_{i,j,k} \Stag(\FmagL)_{i,j,k} \cdot \Stag(f^{OU,M})_{i,j,k},\\
        \label{eq:polyeinjM} &&
\eeqa

with $\Stag(g)_{i,j,k}=(\area{g}{x}{}{x,i-1/2,j,k}, \area{g}{y}{}{y,i,j-1/2,k}, \area{g}{z}{}{z,i,j,k-1/2})$ the vector containing the 3 components of the staggered vector field $\vg$, $\Dt$ the time-step used for the time integration (see section \ref{sec:timeINT}) and $\Delta \Ekin=\EKinj\Dt, \Delta \Emag=\EMinj\Dt$.

\subsubsection{Applying the forcing terms}
\label{sec:forcing_apply}

While the staggered electromotive driving $\area{f}{n}{}{M,n}=\fAM \area{f}{n}{OU,M}{n} (n \in \{x,y,z\})$ can be directly added to the \rhs (see \fig{fig:code_toHBrhs}), the mechanical driving is defined as a point-valued field at each $(x_i,y_j,z_k)$:

\beq
\Kforceijk = \rho_{i,j,k}\fAK\fOUK_{i,j,k},
\eeq

hence, a point value$\to$volume average transformation is performed before adding $\vol{\Kforce}$ to the \rhs of the momentum equation (block \modVforce in \fig{fig:forcing}). Note the multiplication by the mass density $\rho$, explained in section \ref{sec:forcing_remarks}. 

Similarly, the effect of both mechanical and electromotive drivings on the energy equation \eq{eq:dtedens} is, at each point $(x_i,y_j,z_k)$:

\beq
\KMforceEijk = \Fvel_{i,j,k}\cdot\Kforceijk+\Fmag_{i,j,k}\cdot(\fAM\fOUM_{i,j,k}).
\eeq

It is computed using the point-valued quantities from the block \modAtoP (\fig{fig:code_HBflux}) and transformed into a volume-averaged quantity (block \modFEMK followed by \modVforce in \fig{fig:forcing}), so that $\vol{\KMforceE}$ is added to the \rhs of $\edens$ (\fig{fig:code_toHBrhs}).

The formulas for point value$\leftrightarrow$volume averages transformations are given in appendix \ref{app:transf_form_vp}.

\subsubsection{Remarks for compressible MHD turbulent simulations}
\label{sec:forcing_remarks}

We conclude the description of the forcing module with two remarks relevant in the context of compressible MHD turbulence:
\begin{itemize}
\item The kinetic forcing occurs through an external acceleration field (i.e. it is of the form $\pat (\rho \Fvel)=\rho \vf$). This allows a broader inertial range in the turbulent spectra as compared to a forcing through an external \emph{force} field (of the form $\pat(\rho\Fvel)=\vg$) \cite{KWN13,ALU13}.
\item The finite-volume approach guarantees the strict conservation of \emph{momentum} up to machine precision. However, the statistical process that represents the turbulence driver can still generate a finite mean \emph{velocity} component. The weak mean velocity field is thus removed at each call to the forcing procedure. This is important in the presence of magnetic fields as those break the Galilean invariance valid in purely hydrodynamic turbulence.
\end{itemize}

\subsection{Time integration: Strong-Stability Preserving Runge-Kutta methods \modiSSPRK}
\label{sec:timeINT}

\newcommand{\textnit}{it==s}

\ifx\showtikz\undefined
\else
\newcommand{\truechar}{\textcolor{green}{\ding{51}}}
\newcommand{\falsechar}{\textcolor{red}{\ding{55}}}
\def\IFbsize{1.}
\def\scsvac{0.7}
\def\lins{0.25}
\def\baseRHS{1.5*\scsvac+\lins}
\def\baseSSPRK{\baseRHS+\bssize+\lout+1.5*\scsvac+0.3+\lins}
\def\mIF{\IFbsize*0.20711}
\def\baseIFiter{\baseSSPRK+\bssize+\lout+\mIF}
\def\baseDiag{\baseIFiter+\IFbsize+\mIF+\lout}
\def\baseTF{\baseDiag+\bDiagsize+\lout+\mIF}
\def\bDiagsize{1.3}
\def\hSDiagblock{0.75}
\def\yarrSSPRK{-0.5}
\def\yarrITloop{-0.75}
\def\yarrgenloop{-1.}
\def\baseINIT{-\bDiagsize-\lins}
\def\yIF{0.14}

\def\IFtx{\IFbsize+\mIF+0.07}
\def\IFty{0.2}
\def\IFfx{\IFbsize*0.5+0.12}
\def\IFfy{-\IFbsize*0.5-\mIF-0.1}

\def\bssize{1.1}

\begin{figure}
\centering
\begin{tikzpicture}

\draw[fill=\clrBdiag,fill opacity=0.3] (\baseINIT,\yarrLO-\hSDiagblock) rectangle (\baseINIT+\bDiagsize,\yarrLO+\hSDiagblock);
\node at (\baseINIT+\bDiagsize*0.5,\yarrLO) {\module{Init}};

\draw[->] (\baseINIT+\bDiagsize,\yarrLO)--(\baseINIT+\bDiagsize+\lins,\yarrLO);

\path (0,\yarrLO-\yarrcROffs*\scsvac) pic {svolareacube=\clrH/\clrB/\scsvac};

\draw[->] (1.5*\scsvac,\yarrLO)--(\baseRHS,\yarrLO);

\draw[fill=\clrBgen,fill opacity=0.3] (\baseRHS,\yarrLO-\hSblock) rectangle (\baseRHS+\bssize,\yarrLO+\hSblock);
\node at (\baseRHS+\bssize*0.5,\yarrLO) {\modRHS};

\draw[->] (\baseRHS+\bssize,\yarrLO)--(\baseRHS+\bssize+\lout,\yarrLO);
\node at (\baseRHS+\bssize+\lout+0.15,\yarrLO) {$\pat$};
\path (\baseRHS+\bssize+\lout+0.3,\yarrLO-\yarrcROffs*\scsvac) pic {svolareacube=\clrH/\clrB/\scsvac};

\draw[->] (\baseRHS+\bssize+\lout+0.3+1.5*\scsvac,\yarrLO)--(\baseSSPRK,\yarrLO);
\draw[fill=\clrBgen,fill opacity=0.3] (\baseSSPRK,\yarrLO-\hSblock) rectangle (\baseSSPRK+\bssize,\yarrLO+\hSblock);
\node at (\baseSSPRK+\bssize*0.5,\yarrLO) {\modSSPRK};

\draw[->] (\baseSSPRK+\bssize,\yarrLO)--(\baseIFiter-\mIF,\yarrLO);
\draw[fill=\clrBcond,fill opacity=0.3,rotate around={45:(\baseIFiter+\IFbsize*0.5,\yarrLO)}] (\baseIFiter,\yarrLO-\IFbsize*0.5) rectangle (\baseIFiter+\IFbsize,\yarrLO+\IFbsize*0.5);
\node at (\baseIFiter+\IFbsize*0.5,\yarrLO+\yIF) {\begin{tabular}{c} \module{if}\\ \module{\textnit}\end{tabular}};

\node at (\baseIFiter+\IFtx,\yarrLO+\IFty) {\truechar};
\node at (\baseIFiter+\IFfx,\yarrLO+\IFfy) {\falsechar};

\draw[->] (\baseIFiter+\IFbsize+\mIF,\yarrLO)--(\baseDiag,\yarrLO);

\draw[fill=\clrBdiag,fill opacity=0.3] (\baseDiag,\yarrLO-\hSDiagblock) rectangle (\baseDiag+\bDiagsize,\yarrLO+\hSDiagblock);
\node at (\baseDiag+\bDiagsize*0.5,\yarrLO) {\begin{tabular}{c} \module{t$\leftarrow$t+dt} \\\module{Diags}\end{tabular}};

\draw[->] (\baseDiag+\bDiagsize,\yarrLO)--(\baseTF-\mIF,\yarrLO);
\draw[fill=\clrBcond,fill opacity=0.3,rotate around={45:(\baseTF+\IFbsize*0.5,\yarrLO)}] (\baseTF,\yarrLO-\IFbsize*0.5) rectangle (\baseTF+\IFbsize,\yarrLO+\IFbsize*0.5);
\node at (\baseTF+\IFbsize*0.5,\yarrLO+\yIF) {\begin{tabular}{c}\module{if} \\ \module{t==tf}\end{tabular}};

\node at (\baseTF+\IFtx,\yarrLO+\IFty) {\truechar};
\node at (\baseTF+\IFfx,\yarrLO+\IFfy) {\falsechar};

\draw[->] (\baseTF+\IFbsize+\mIF,\yarrLO)--(\baseTF+\IFbsize+\mIF+\lout,\yarrLO);
\draw[fill=\clrBdiag,fill opacity=0.3] (\baseTF+\IFbsize+\mIF+\lout,\yarrLO-\hSDiagblock) rectangle (\baseTF+\IFbsize+\mIF+\lout+\bDiagsize,\yarrLO+\hSDiagblock);
\node at (\baseTF+\IFbsize+\mIF+\lout+\bDiagsize*0.5,\yarrLO) {\begin{tabular}{c} \module{End}\end{tabular}};

\draw[->] (\baseIFiter+\IFbsize*0.5,\yarrLO-\IFbsize*0.5-\mIF)--(\baseIFiter+\IFbsize*0.5,\yarrITloop)--(-\lins*0.5+\baseIFiter*0.5+\IFbsize*0.25,\yarrITloop);
\draw (-\lins*0.5+\baseIFiter*0.5+\IFbsize*0.25,\yarrITloop)--(-\lins*0.5,\yarrITloop) -- (-\lins*0.5,\yarrLO)--(0,\yarrLO);
\draw[->] (\baseTF+\IFbsize*0.5,\yarrLO-\IFbsize*0.5-\mIF)--(\baseTF+\IFbsize*0.5,\yarrgenloop)--(-\lins*0.5+\baseTF*0.5+\IFbsize*0.25,\yarrgenloop);
\draw[->] (-\lins*0.5+\baseTF*0.5+\IFbsize*0.25,\yarrgenloop)--(-\lins*0.5,\yarrgenloop) -- (-\lins*0.5,\yarrLO-1.25);

\end{tikzpicture}
\caption{\label{fig:code_wholescheme}Complete structure of a solver. The current iteration of the \module{s}-stages SSPRK method is noted \module{it}. The variable timestep is \module{dt}. After initialization (module \module{Init}), a main loop occurs till the final time \module{tf} is reached when the code releases the memory and exits (module \module{End}).}
\end{figure}
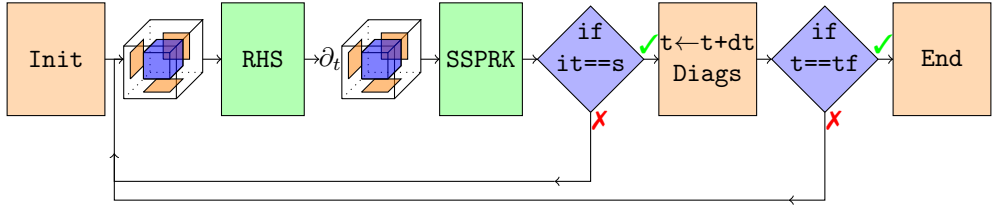
\fi

The computed \rhs (\modRHS, sections \ref{sec:HBreco}-\ref{sec:forcing}) is plugged in a Strong Stability Preserving Runge-Kutta (\modSSPRK) time integrator (see the diagram for the whole solver, \fig{fig:code_wholescheme}). The SSPRK methods require several iterations, so that the RHS is computed \module{s} times (see the loop with the ending condition \module{\textnit} in \fig{fig:code_wholescheme}). After this, the simulation time \module{t} is updated. If needed, several diagnostics are written out, e.g. complete state of the system, or only its kinetic/magnetic energy, Fourier spectra, etc. Then, the time integration procedure is repeated. This is performed until the final time \module{tf} is reached.
 
SSPRK methods prevent additional oscillations resulting from inaccuracies in the time integration process \cite{GOS98,GST01}. Although only smooth subsonic problems are considered in the present work, SSPRK methods also are very valuable in systems exhibiting discontinuities and shocks, e.g. turbulent supersonic fluids. 

The second-order and fourth-order time integration methods are taken from \cite{KET08} (pseudocode 1 with \module{s=10} and pseudocode 3 of that reference, respectively). They are ten-stages single-step SSPRK methods. Ten-stages means that \module{s=10} RHS computations are needed. Single-step means that the solution at the next instant $\solvarP$ is computed solely from the solution at the instant $t$, $\solvar$.

However, for schemes of order strictly greater than four, no explicit single-step SSPRK method exists \cite{KRA91,RSP02}. For this reason, \multistep methods, which require the solution at several points in time, have been developed \cite{KGM11,BGG17,SSPSITE}. Two-step explicit methods, which require $\solvarM$ in addition to $\solvar$ to compute $\solvarP$, have order of accuracy at most eight \cite{KGM11}. Hence, for tenth order of accuracy, three-step methods (at least), which require $\solvarMM$ as well, are needed. 

Consequently, the sixth-order method used is the two-step ten-stages one from \cite{SSPSITE}. The two-step eighth-order methods presented in \cite{KGM11} require a smaller \timestep as compared to the three-step ten-stages one from \cite{SSPSITE}, so that the latter is used. Finally, for tenth order of accuracy, a three-step twenty-stages method is used \cite{SSPSITE}. In the following, we denote the applied time integrators $SSP(p,s,n)$, with $p$ the method's order of accuracy, $s$ the number of stages and $n$ the number of steps.

\ifx\showtikz\undefined
\else
\def\endMSRK{10}
\def\unitT{1}
\def\yT{0.1}
\def\fyT{0.2}
\def\lRK{-0.125}
\def\llRK{0.35}
\def\TllRK{-0.35}

\def\arrShape{-stealth}

\def\loosellRK{0.6}
\def\ainu{-150}
\def\aoutu{-30}
\def\aino{150}
\def\aouto{30}
\def\looseF{0.5}

\def\clrArrA{black}
\def\clrArrB{orange}
\def\clrArrC{cyan}
\def\clrArrD{orange!80!black}

\def\clrArrTA{black!70}
\def\clrArrTB{black!70}
\def\clrArrTC{black!70}
\def\clrArrTD{orange!70!black}

\def\TS{-1.8}

\def\ySInt{-0.2}

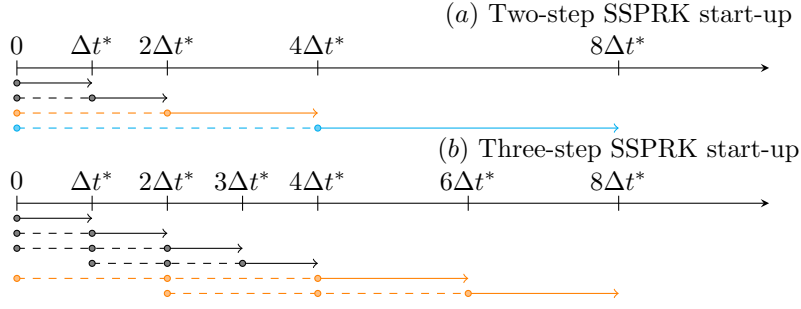
\begin{figure}
\centering
\begin{tikzpicture}

\draw[\arrShape] (0,0)--(\endMSRK,0);
\draw (0,-\yT) -- (0,\yT);
\draw (\unitT,-\yT) -- (\unitT,\yT);
\draw (\unitT*2,-\yT) -- (\unitT*2,\yT);
\draw (\unitT*4,-\yT) -- (\unitT*4,\yT);
\draw (\unitT*8,-\yT) -- (\unitT*8,\yT);

\node at (\unitT*8,0.7) {$(a)$ Two-step SSPRK start-up};
\node at (\unitT*8,\TS+0.7) {$(b)$ Three-step SSPRK start-up};

\draw[->,\clrArrA] (0,\ySInt) -- (\unitT,\ySInt);
\path (0,\ySInt) pic {timeDot=\clrArrA};

\draw[dashed,\clrArrA] (0*\unitT,2*\ySInt) -- (1*\unitT,2*\ySInt);
\draw[->,\clrArrA] (\unitT,2*\ySInt) -- (2*\unitT,2*\ySInt);
\path (0,2*\ySInt) pic {timeDot=\clrArrA};
\path (\unitT,2*\ySInt) pic {timeDot=\clrArrA};

\draw[dashed,\clrArrB] (0*\unitT,3*\ySInt) -- (2*\unitT,3*\ySInt);
\draw[->,\clrArrB] (2*\unitT,3*\ySInt) -- (4*\unitT,3*\ySInt);
\path (0,3*\ySInt) pic {timeDot=\clrArrB};
\path (2*\unitT,3*\ySInt) pic {timeDot=\clrArrB};

\draw[dashed,\clrArrC] (0*\unitT,4*\ySInt) -- (4*\unitT,4*\ySInt);
\draw[->,\clrArrC] (4*\unitT,4*\ySInt) -- (8*\unitT,4*\ySInt);
\path (0,4*\ySInt) pic {timeDot=\clrArrC};
\path (4*\unitT,4*\ySInt) pic {timeDot=\clrArrC};

\node at (0,\yT+\fyT) {$0$};
\node at (\unitT,\yT+\fyT) {$\dts$};
\node at (\unitT*2,\yT+\fyT) {$2\dts$};
\node at (\unitT*4,\yT+\fyT) {$4\dts$};
\node at (\unitT*8,\yT+\fyT) {$8\dts$};

\draw[\arrShape] (0,\TS)--(\endMSRK,\TS);
\draw (0,\TS-\yT) -- (0,\TS+\yT);
\draw (\unitT,\TS-\yT) -- (\unitT,\TS+\yT);
\draw (\unitT*2,\TS-\yT) -- (\unitT*2,\TS+\yT);
\draw (\unitT*3,\TS-\yT) -- (\unitT*3,\TS+\yT);
\draw (\unitT*4,\TS-\yT) -- (\unitT*4,\TS+\yT);
\draw (\unitT*6,\TS-\yT) -- (\unitT*6,\TS+\yT);
\draw (\unitT*8,\TS-\yT) -- (\unitT*8,\TS+\yT);

\draw[->,\clrArrA] (0,\TS+\ySInt) -- (\unitT,\TS+\ySInt);
\path (0,\TS+\ySInt) pic {timeDot=\clrArrA};

\draw[dashed,\clrArrA] (0*\unitT,\TS+2*\ySInt) -- (1*\unitT,\TS+2*\ySInt);
\draw[->,\clrArrA] (\unitT,\TS+2*\ySInt) -- (2*\unitT,\TS+2*\ySInt);
\path (0,\TS+2*\ySInt) pic {timeDot=\clrArrA};
\path (\unitT,\TS+2*\ySInt) pic {timeDot=\clrArrA};

\draw[dashed,\clrArrA] (0*\unitT,\TS+3*\ySInt) -- (2*\unitT,\TS+3*\ySInt);
\draw[->,\clrArrA] (2*\unitT,\TS+3*\ySInt) -- (3*\unitT,\TS+3*\ySInt);
\path (0,\TS+3*\ySInt) pic {timeDot=\clrArrA};
\path (\unitT,\TS+3*\ySInt) pic {timeDot=\clrArrA};
\path (2*\unitT,\TS+3*\ySInt) pic {timeDot=\clrArrA};

\draw[dashed,\clrArrA] (1*\unitT,\TS+4*\ySInt) -- (3*\unitT,\TS+4*\ySInt);
\draw[->,\clrArrA] (3*\unitT,\TS+4*\ySInt) -- (4*\unitT,\TS+4*\ySInt);
\path (\unitT,\TS+4*\ySInt) pic {timeDot=\clrArrA};
\path (2*\unitT,\TS+4*\ySInt) pic {timeDot=\clrArrA};
\path (3*\unitT,\TS+4*\ySInt) pic {timeDot=\clrArrA};

\draw[dashed,\clrArrB] (0*\unitT,\TS+5*\ySInt) -- (4*\unitT,\TS+5*\ySInt);
\draw[->,\clrArrB] (4*\unitT,\TS+5*\ySInt) -- (6*\unitT,\TS+5*\ySInt);
\path (0*\unitT,\TS+5*\ySInt) pic {timeDot=\clrArrB};
\path (2*\unitT,\TS+5*\ySInt) pic {timeDot=\clrArrB};
\path (4*\unitT,\TS+5*\ySInt) pic {timeDot=\clrArrB};

\draw[dashed,\clrArrB] (2*\unitT,\TS+6*\ySInt) -- (6*\unitT,\TS+6*\ySInt);
\draw[->,\clrArrB] (6*\unitT,\TS+6*\ySInt) -- (8*\unitT,\TS+6*\ySInt);
\path (2*\unitT,\TS+6*\ySInt) pic {timeDot=\clrArrB};
\path (4*\unitT,\TS+6*\ySInt) pic {timeDot=\clrArrB};
\path (6*\unitT,\TS+6*\ySInt) pic {timeDot=\clrArrB};

\node at (0,\TS+\yT+\fyT) {$0$};
\node at (\unitT,\TS+\yT+\fyT) {$\dts$};
\node at (\unitT*2,\TS+\yT+\fyT) {$2\dts$};
\node at (\unitT*3,\TS+\yT+\fyT) {$3\dts$};
\node at (\unitT*4,\TS+\yT+\fyT) {$4\dts$};
\node at (\unitT*6,\TS+\yT+\fyT) {$6\dts$};
\node at (\unitT*8,\TS+\yT+\fyT) {$8\dts$};

\end{tikzpicture}
\caption{\label{fig:msrk}Start-up procedure for two-and three-step SSPRK methods, until $8\dts$ is reached. The dots linked with dashed lines correspond to the instants in time required to determine the next state. A change of color means a doubling of the \timestep.}
\end{figure}
\fi

Since \multistep methods require several instants in time to determine the next one, a particular start-up procedure is needed. The start-up is performed by doubling step-by-step $\dts=\DtCFL/2^\alpha$, $\alpha\in\mathbb{N}$, chosen initially small enough so that the overall accuracy of the scheme, $p$, is preserved \cite{KGM11}:

\begin{equation}
(\dts)^{q+1}=\bigO(\DtCFL^p),
\end{equation}

with $q$ the order of the single-step procedure used for the first iteration. In this work, the first iteration is performed by the $SSP(q=4,10,1)$ method. Then, for a two-step method $SSP(p,s,2)$, $\solvarA{0}$ and $\solvarA{\dts}$ allow to deduce the solution at $2\dts$. Knowing $\solvarA{0}$ and $\solvarA{2\dts}$, the \timestep can be doubled and one can determine $\solvarA{4\dts}$, etc. (\fig{fig:msrk}.$(a)$). This process is repeated until $\DtCFL$ is reached. For a three-step method $SSP(p,s,3)$, the $SSP(6,10,2)$ method is used to deduce $\solvarA{2\dts}$ from $(\solvarA{0},\solvarA{\dts})$. Starting then, only the highest-order $SSP(p,s,3)$ method is used. It is used two times before each doubling of the \timestep: $\solvarA{3\dts}$ and $\solvarA{4\dts}$ are gained successively from $(\solvarA{2\dts},\solvarA{\dts},\solvarA{0})$ and $(\solvarA{3\dts},\solvarA{2\dts},\solvarA{\dts})$. The timestep can then be doubled: $(\solvarA{4\dts},$\\$\solvarA{2\dts},\solvarA{0})$ allows to determine $\solvarA{6\dts}$, etc. (\fig{fig:msrk}.$(b)$).

The \timestep $\Dt$ is limited by the Courant-Friedrichs-Lewy criterion:

\newcommand{\Ceff}{C_{eff}}
\newcommand{\CCFL}{C_{CFL}}
\newcommand{\KCFL}{K}
\newcommand{\CVISCO}{C_{VISCO}}
\newcommand{\CRESI}{C_{RESI}}
\newcommand{\CEsink}{C_{Esink}}
\newcommand{\CSSP}{C_{SSP}}

\beq
\label{eq:CFLcond}
\Dt \leq \DtCFL=\CCFL \underset{i,j,k}{\min}(\Delta x/a^x_{i,j,k},\Delta y/a^y_{i,j,k},\Delta z/a^z_{i,j,k}),
\eeq
with $a^m$ the maximum speed of propagation of information in the $\vm$-direction (see section \ref{sec:LLF}) and $\CCFL$ the Courant number. 

\begin{table}
\centering
\begin{tabular}{c|c|c}
\textbf{Method [source]} & ${\bf \CCFL}$  & ${\bf \CVISCO}$ \\
\hline
$SSP(2,10,1)$ \cite{KET08} & 1.9 & 5.5 \\
\hline
$SSP(4,10,1)$ \cite{KET08} & 1.9 & 1.4\\
\hline
$SSP(6,10,2)$ \cite{SSPSITE} & 1.5 & 0.6\\
\hline
$SSP(8,10,3)$ \cite{SSPSITE}& 1.2 & 0.3\\
\hline
$SSP(10,20,3)$ \cite{SSPSITE}& 2.1 & 0.5\\

\end{tabular}
\caption{\label{tab:CFL} Time integrators and Courant number used for the simulations. Note that although the 10th order method has a higher $\CCFL$, it needs twenty stages as compared to the other ones that need only ten.}
\end{table}

In non-ideal systems, the \timestep may be further reduced by the viscosity, the magnetic diffusivity or the internal energy sink:

\beqa
\label{eq:dtvisco} \Dt \leq \DtVISCO&=&\CVISCO\frac{\min(\vol{\rho}_{i,j,k})\min(\Dx,\Dy,\Dz)^2}{\visc},\\
\label{eq:dtresi} \Dt \leq \DtRESI&=&\CRESI\frac{\min(\Dx,\Dy,\Dz)^2}{\resis},\\
\Dt \leq \DtEsink&=&\CEsink\frac{1}{\UsinkF \max(\intE)^3},
\eeqa

with a priori different constants $\CVISCO, \CRESI, \CEsink$. The numerical solver is stable for:

\beq
\Dt \leq \min(\DtCFL,\DtVISCO,\DtRESI,\DtEsink).
\eeq

The Courant numbers used in this work are listed in \tab{tab:CFL}.

\subsection{Numerical method: summary}
\label{sec:numSummary}

To summarize, the key points of the presented method are the following:
\begin{itemize}
\item The use of one-dimensional WENO reconstruction methods, which are computationally more affordable than multidimensional ones, especially with increasing order of accuracy.
\item This is made possible by the transformation of area-averages into point values in the middle of each face up to arbitrarily high order of accuracy. The procedure prevents the need to compute several points on each face for the application of a quadrature rule.
\item The viscous and resistive terms are formulated and computed in a way respecting the finite-volume and constrained-transport frameworks as well as the scheme's order.
\item The internal energy sink and the forcing terms require point value$\leftrightarrow$volume average transformations, which can be derived to an arbitrarily high order of accuracy. 
\item High-order strong stability preserving explicit time integrators are available in the literature. For order strictly greater than four, these are multi-step methods. 
\end{itemize}

\section{Numerical tests}
\label{sec:numtests}

Numerical solvers of order 2, 4, 6, 8 and 10 are considered. They are denoted by S2, S4, S6, S8 and S10 and use the method of the corresponding order for the magnetic field volume interpolation (section \ref{sec:Binterp}), the passage through point values (sections \ref{sec:colella} and \ref{sec:Esink}) and the time integration (section \ref{sec:timeINT}). The employed WENO schemes (section \ref{sec:WENO}) are of the next closest odd-order, that is 3, 5, 7, 9 and 11, respectively. The Courant numbers used are given in \tab{tab:CFL}.

The approach presented in this work is validated in two steps:
\begin{itemize}
\item first, a convergence test of the ideal MHD equations ($\visc=\resis=\Usink=0$) through the \nonlinear 3D MHD vortex problem (section \ref{sec:MHVort}),
\item second, a convergence test for the viscous and resistive terms (section \ref{sec:testNID}). 
\end{itemize}

Additionally, the MHD vortex test describes how a higher-order scheme leads to a higher computational efficiency by strongly reducing numerical \nonideal effects. This is also illustrated through the inertial range width and the level of visible structure details when simulating turbulent systems (section \ref{sec:ssturb}).

\subsection{Validation of ideal MHD: 3D MHD vortex problem}

\label{sec:MHVort}

\newcommand{\dUmeanCTEST}{\delta Q}
\newcommand{\EOCCTEST}{EOC}
\newcommand{\ElossCTEST}{E_{loss}}
\sisetup{round-mode=figures,round-precision=3}
\newcommand{\numCTESTerrf}[1]{\num[round-mode=figures,round-precision=2,scientific-notation=true]{#1}}
\newcommand{\numCTESTeocf}[1]{\num[scientific-notation=fixed,fixed-exponent=0]{#1}}
\newcommand{\numCTESTelossf}[1]{\num[round-mode=figures,round-precision=2,scientific-notation=true]{#1}}
\newcommand{\numCTESTelossflarge}[1]{\num[scientific-notation=fixed,fixed-exponent=0,round-precision=2]{#1}}
\newcommand{\nce}[1]{\numCTESTerrf{#1}}
\newcommand{\ncel}[1]{\numCTESTerrflarge{#1}}
\newcommand{\ncc}[1]{\numCTESTeocf{#1}}
\newcommand{\ncl}[1]{\numCTESTelossf{#1}}
\newcommand{\ncll}[1]{\numCTESTelossflarge{#1}}
\newcommand{\nceff}[1]{\num[round-mode=figures,round-precision=2,scientific-notation=true]{#1}}

\newcommand{\errsign}{\delta}
\newcommand{\LIIerr}{\errsign_2}
\newcommand{\LIIeoc}{\EOCCTEST_2}

The 3D MHD vortex is a nonlinear test case where a smooth magnetized vortex structure in force equilibrium is advected by a mean velocity field. The MHD vortex has first been introduced in 2D \cite{BAL04} and then extended to 3D \cite{MTB10} and has the initial conditions:

\newcommand{\qMH}{\phi}
\newcommand{\pinf}{p_0}

   \begin{equation}
        \label{eq:init3Dmhvort}
        \begin{pmatrix} \rho \\ \FvelL_x \\ \FvelL_y \\ \FvelL_z \\ \press \\ \FmagL_x \\ \FmagL_y \\ \FmagL_z \end{pmatrix} = \begin{pmatrix} 1 \\ 1-{y\kappa}\exp\Big[\qMH(1-r^2)\Big] \\ 1+{x\kappa}\exp\Big[\qMH(1-r^2)\Big] \\ {2} \\ \pinf+ \frac{1}{4\qMH}\big[\mu^2\big(1-2\qMH(r^2-z^2)\big)-\kappa^2 \rho\big]\exp\Big[2\qMH(1-r^2)\Big] \\ {-y\mu}\exp\Big[\qMH(1-r^2)\Big] \\ {x\mu\exp\Big[\qMH(1-r^2)\Big]} \\ 0 \end{pmatrix},
        \end{equation}
with $r=\sqrt{x^2+y^2+z^2}$ in the triply periodic computational domain $(x,y,z) \in [-5,5]\times[-5,5]\times[-5,5]$ and the parameters $\kappa=\mu=1/(2\pi)$, $\qMH=1.75$, $\pinf=360$. The ambient pressure $\pinf$ is set to be higher than in \cite{BAL04,MTB10}, so that the advection speed is about one-tenth the speed of sound. The vortex is advected for one period of motion, until $t=10$ and then compared to the initial conditions at $t=0$. The error is measured taking the $L_2$-norm:

\beq
\LIIerr(N)=\sqrt{\frac{1}{N_qN^3}\sum_{i,j,k,q} (q_{i,j,k}(t=10)-q_{i,j,k}(t=0))^2},\\
\eeq

with $q\in\{\rho,\rho \FvelL_x,\rho \FvelL_y,\rho \FvelL_z,\edens,\FmagL_x,\FmagL_y,\FmagL_z\}$ the $N_q=8$ variables and $N=N_x=N_y=N_z$ the linear resolution in each dimension so that the system is discretized in $N^3$ grid-cells. The simulations are repeated at several resolutions $(N_i)$ and the experimental order of convergence $\LIIeoc$, expected to converge to the theoretical order of the solver for $N$ high enough, is computed for $i>1$ by:

\beq
\LIIeoc(N_i)=-\frac{\log(\LIIerr(N_i))-\log(\LIIerr(N_{i-1}))}{\log(N_i)-\log(N_{i-1})}.
\eeq

\begin{table}
\centering
 \begin{tabular}{cc|cc|c}
 & Resolution & $\LIIerr$ & $\LIIeoc$ & $\ElossCTEST$ \\
S2
 & $64^3$ & \nce{0.8467881270594523E-02} & - &  \ncll{0.9945339053261583E+00}
 \\
 & $128^3$ & \nce{0.7388070583460179E-02} & \ncc{0.1968033917768960E+00} &  \ncll{0.9380735705557063E+00}
 \\
 & $192^3$ & \nce{0.5660224292466662E-02} & \ncc{0.6570308805564111E+00} & \ncll{0.7908866112880171E+00}
 \\
 & $256^3$ & \nce{0.4063848342718727E-02} & \ncc{0.1151733661798797E+01} & \ncll{0.5980864414110383E+00}
 \\
 \hline
S4
 & $64^3$ & \nce{0.3869063248683349E-02} & - & \ncll{0.5933299292673263E+00}
 \\
 & $96^3$ & \nce{0.1272692731560324E-02} & \ncc{0.2742227343224463E+01} & \ncll{0.1963064432539060E+00}
 \\
 & $128^3$ & \nce{0.4360273434459408E-03} & \ncc{0.3723503627129372E+01} & \ncll{0.6140726499177872E-01}
 \\
 & $160^3$ & \nce{0.1684515008137889E-03} & \ncc{0.4262086315035480E+01} & \ncll{0.2186819597810624E-01}
 \\
 & $192^3$ & \nce{0.7467866685585436E-04} & \ncc{0.4461641434917484E+01} & \ncl{0.9041972250249074E-02}
 \\
 \hline
S6
 & $64^3$ & \nce{0.8712744863388463E-03} & - & \ncll{0.8173544854691293E-01}
 \\
 & $96^3$ & \nce{0.1008536051190384E-03} & \ncc{0.5318058230652527E+01} & \ncl{0.7734355623655129E-02}
 \\
 & $128^3$ & \nce{0.1580682966960341E-04} & \ncc{0.6441930465457030E+01} & \ncl{0.1179556359084977E-02}
 \\
 & $160^3$ & \nce{0.4214649896925542E-05} & \ncc{0.5923879805823349E+01} & \ncl{0.2643495843483900E-03}
 \\
 & $192^3$ & \nce{0.1380616243938451E-05} & \ncc{0.6121254081684615E+01} & \ncl{0.7743761973388546E-04}
 \\
 \hline
S8
 & $64^3$ & \nce{0.1432211613243621E-03} & - & \ncl{0.9430512801770716E-02}
 \\
 & $96^3$ & \nce{0.4852181170202610E-05} & \ncc{0.8348342732488355E+01} & \ncl{0.3940328293006451E-03}
 \\
 & $128^3$ & \nce{0.4703764970400547E-06} & \ncc{0.8111906860238960E+01} & \ncl{0.3415764026407281E-04}
 \\
 & $160^3$ & \nce{0.9030726672698769E-07} & \ncc{0.7395757087708716E+01} & \ncl{0.4840408615084994E-05}
 \\
 & $192^3$ & \nce{0.2148786911855860E-07} & \ncc{0.7874709972561968E+01} & \ncl{0.9606137429411182E-06}
 \\
 \hline
S10
 & $64^3$ & \nce{0.2892911186673655E-04} & - & \ncl{0.8422618030612288E-03}
 \\
 & $96^3$ & \nce{0.5378878753627064E-06} & \ncc{0.9828104786550686E+01} & \ncl{0.1624661204993432E-04}
 \\
 & $128^3$ & \nce{0.2975662654780163E-07} & \ncc{0.1006179584071565E+02} & \ncl{0.7974138576917064E-06}
 \\
 & $160^3$ & \nce{0.2926845323732749E-08} & \ncc{0.1039298098896949E+02} & \ncl{0.7596145634725345E-07}
 \\
 & $192^3$ & \nce{0.4953140071235243E-09} & \ncc{0.9743710796562816E+01} & \ncl{0.1107146476305681E-07}
 \\
 \hline
 \hline
S8P
 & $64^3$ & \nce{0.1438464768276557E-03} & - & \ncl{0.9441423591531101E-02}
 \\
 & $96^3$ & \nce{0.1177184582746036E-04} & \ncc{0.6173246017764224E+01} & \ncl{0.3960389776004833E-03}
 \\
 & $128^3$ & \nce{0.6422904525091104E-05} & \ncc{0.2105936928728556E+01} & \ncl{0.3490969332002516E-04}
 \\
 & $160^3$ & \nce{0.4137941718850274E-05} & \ncc{0.1970354671658245E+01} & \ncl{0.5173415657171280E-05}
 \\
 & $192^3$ & \nce{0.2880354289330197E-05} & \ncc{0.1987067231697334E+01} & \ncl{0.1127909743612236E-05}
 \\
 \hline
S8T  & $64^3$ & \nce{0.1432316274881725E-03} & - & \ncl{0.9429114218682118E-02}
 \\
 & $96^3$ & \nce{0.4959651735930570E-05} & \ncc{0.8294493167212053E+01} & \ncl{0.3935680145535107E-03}
 \\
 & $128^3$ & \nce{0.7579664656113384E-06} & \ncc{0.6529609727394829E+01} & \ncl{0.3395990136513785E-04}
 \\
 & $160^3$ & \nce{0.3926667116727380E-06} & \ncc{0.2947331213526590E+01} & \ncl{0.4739154742213177E-05}
 \\
 & $192^3$ & \nce{0.2667972646895404E-06} & \ncc{0.2119728107059896E+01} & \ncl{0.9020748972071780E-06}
 \\
\hline

&  &  &  &  \\

 \hline
S2, $t=10/64$ & $64^3$ & \nce{0.3245465895790362E-02} & -  & \ncll{0.4434554738182389E+00}
 \\
 & $128^3$ & \nce{0.1098250804669870E-02} & \ncc{0.1563218039566954E+01}  & \ncll{0.9767718393344049E-01}
 \\
 & $192^3$ & \nce{0.5409344124009545E-03} & \ncc{0.1746576866640475E+01} & \ncll{0.2923995514936233E-01}
 \\
 & $256^3$ & \nce{0.3158399144538625E-03} & \ncc{0.1870337444074101E+01} & \ncll{0.1229379508377567E-01}
 \\
 \hline
\end{tabular}
\caption{\label{tab:CTEST} 3D MHD vortex problem: measured error $\LIIerr$ after one period and corresponding experimental order of convergence $\LIIeoc$, as well as proportion of energy lost $\ElossCTEST$. The S8P and S8T schemes are modifications of the eighth-order S8 one, removing the passage through point values (S8P) or using a second-order time integrator (S8T) instead. Bottom lines: S2 scheme at an earlier time $t=10/64$, instead of one period of motion ($t=10$).}
\end{table}

For this test, the WENO weights are computed using the parameters $p=1,\epsilon=10^{-12}$ (see \eq{eq:wwweno}). The typical choice $p=2,\epsilon=10^{-6}$ leads to a worse convergence for schemes of order strictly greater than 4.

In order to measure the impact of numerical dissipation, the proportion of lost energy is computed by:

\beq
\ElossCTEST=\frac{E(t=0)-E(t=10)}{E(t=0)},\quad E=\sum_{i,j,k} \frac{1}{2} (\rho(\Fvel-\Fvel_{mean})^2+\Fmag^2)_{i,j,k},
\eeq

with $\Fvel_{mean}$ the mean velocity field, $(1,1,2)$, responsible for the advection, so that only the fluctuating field is considered.

%

The dissipation is very large for the S2 scheme. This is because of the employed LLF Riemann solver. At low Mach numbers, the dissipative term $\propto \maxspeedS (q_R-q_L)$ (cf. \eq{eq:fluxLLF}), with $(q_R-q_L)$ the jump of $q$ between the left and right states, does not scale correctly. For this solver, the dissipative term is dominated by the fast magneto-sonic speed, which leads to excessive dissipation \cite{MRE15}. When dealing with flows with a Mach number below $10^{-2}$, commonly used approximate Riemann solvers, such as the Roe solver and HLL-type Riemann solvers, are indeed too dissipative to be of practical use when combined with second-order numerics \cite{LBA22}. With appropriate techniques these Riemann solvers can be modified to have a low dissipation over a wide range of Mach numbers (e.g. \cite{MRE15} for the Roe solver and \cite{MIM21} for the HLLD scheme). An alternative solution is to increase the discretization order, as in the present work. In that case, the jump $(q_R-q_L)$ becomes small enough to prevent excessive dissipation.

To illustrate the necessity of the passage through point values, a modification of the S8 solver is tested (S8P in table \ref{tab:CTEST}). For this solver, the passage through point values is suppressed. It converges as expected at a second-order rate (see \eq{eq:simple2O}). The reduced accuracy is however not reflected in the energy dissipation level (the $\ElossCTEST$ values are very close to the S8 ones), hinting at dispersive effects. Similarly, using a second-order time integrator in place of the eighth-order one (scheme S8T, using a Courant number of 0.5 instead of 1.2) leads asymptotically to a convergence order of two. 

\begin{table}
\centering
\begin{tabular}{ccc}
\hline
\textbf{Solver} & \textbf{Average step duration} ${\bf(ASD)}$ & \textbf{Cost to reach a certain $t_f$} \\
S2 & 1 & 1 \\ 
S4 & 1.3 & 1.3 \\
S6 & 1.6 & 2.0 \\
S8 & 2.2 & 3.5 \\
S10 & 5.9 & 5.3 \\
\hline
\end{tabular}
\caption{\label{tab:teff} Computational cost for each solver: average step duration ($ASD$), as well as the cost to reach a certain instant (this corresponds to $ASD$ divided by the Courant number, see \tab{tab:CFL}). The second-order scheme S2 is taken as the reference. The values are computed from the $512^3$ runs, excluding the start-up procedure of the \multistep SSPRK methods (see section \ref{sec:timeINT}). The same number of cores, compiler options, processor types, etc. are used for all the runs.}
\end{table}

\Tab{tab:CTEST} together with \tab{tab:teff} illustrate the gain in efficiency when seeking higher-order schemes for smooth problems. Even though these are more expensive than lower-order schemes at a given resolution (\tab{tab:teff}), they allow to attain the same level of precision with regard to amplitude errors on substantially coarser grids. For example, the energy losses at resolution $192^3$ for the S4 solver are comparable to those at resolution $96^3$ for the S6 one. Even though the S6 solver is about 50\% more expensive than the S4 one, a factor 2 in linear resolution $N$ means a factor $2^4=16$ in computational costs ($2^3$ due to the increase in spatial resolution and an additional factor 2 because of the \timestep reduction through the CFL criterion, see \eq{eq:CFLcond}). Thus, when using S6 instead of S4, one can gain a factor of order 10 in computational efficiency by performing a run of similar quality at a significantly lower resolution. The same pattern can be found for all increases of discretization order: for this smooth problem with a large scale structure, the gain in computing time resulting from a lower resolution is always higher than the loss because of the cost associated with the use of a higher order scheme.


\subsection{Validation of the dissipative terms}
\label{sec:testNID}

\ifx\showpng\undefined
\else
\begin{figure}[h]
\begin{minipage}[b]{0.33\linewidth}
\centering
\includegraphics[width=\textwidth]{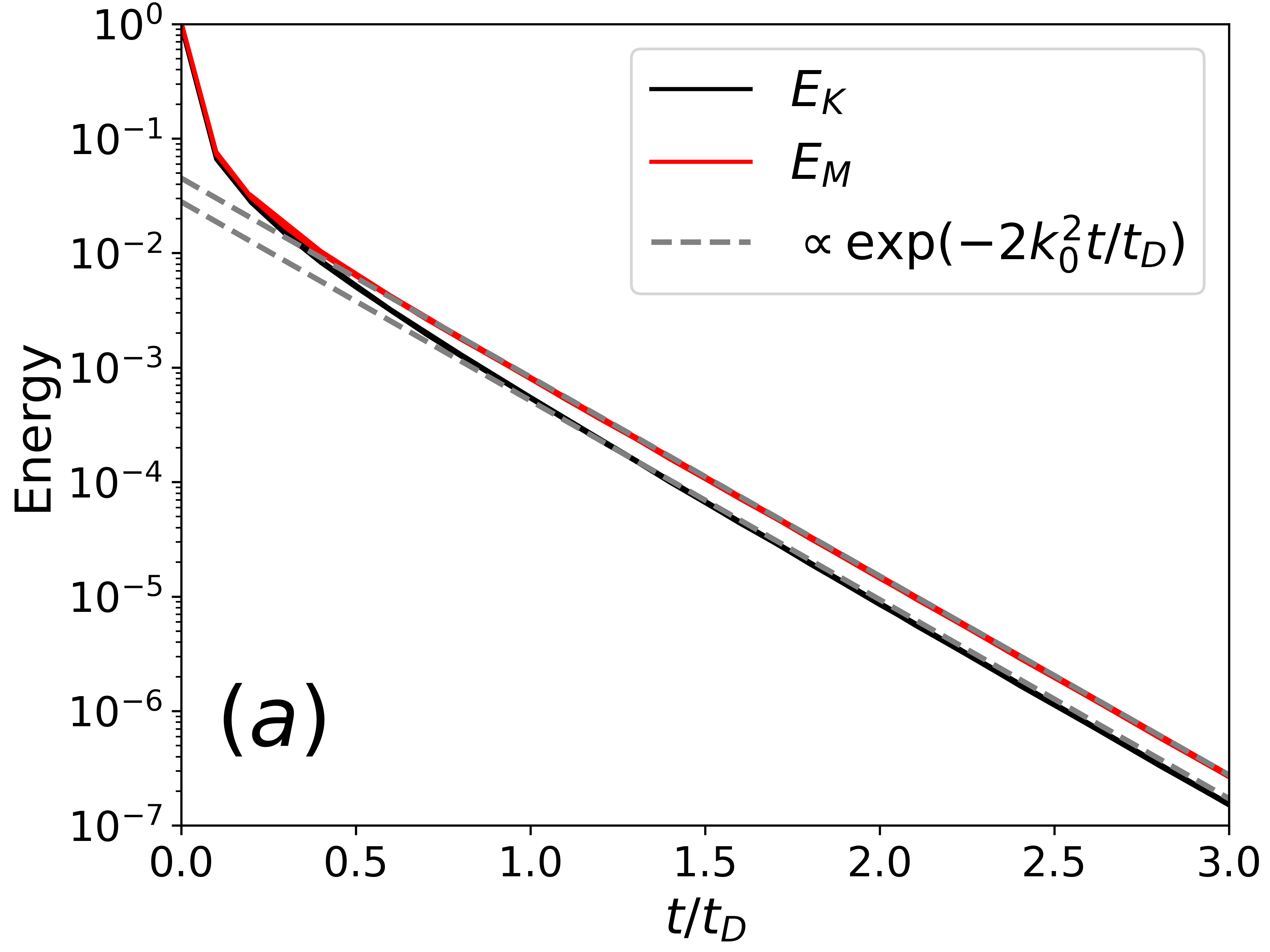}
\end{minipage}%
\begin{minipage}[b]{0.33\linewidth}
\centering
\includegraphics[width=\textwidth]{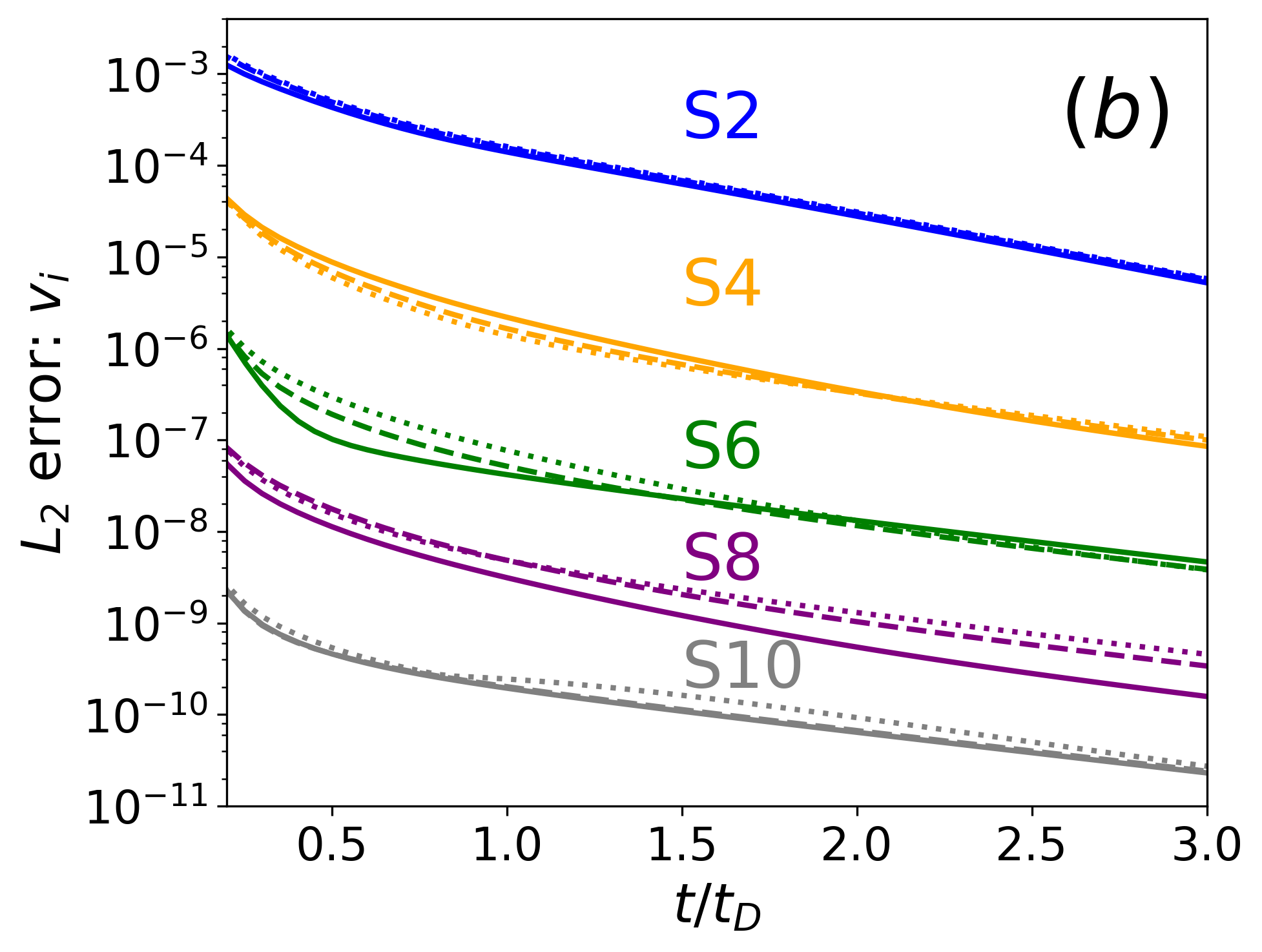}
\end{minipage}%
\begin{minipage}[b]{0.33\linewidth}
\centering
\includegraphics[width=\textwidth]{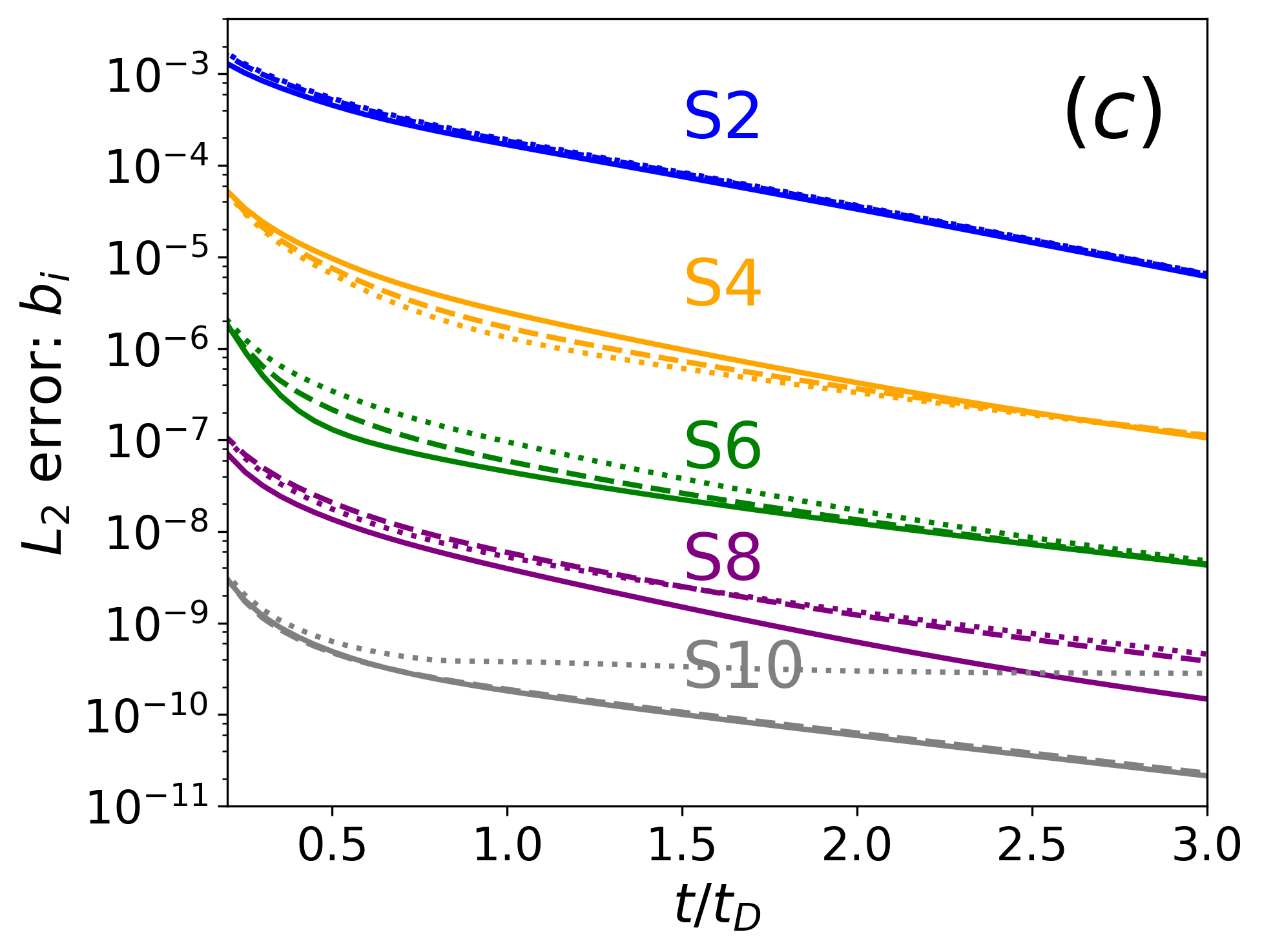}
\end{minipage}
\centering
\includegraphics[width=0.5\textwidth]{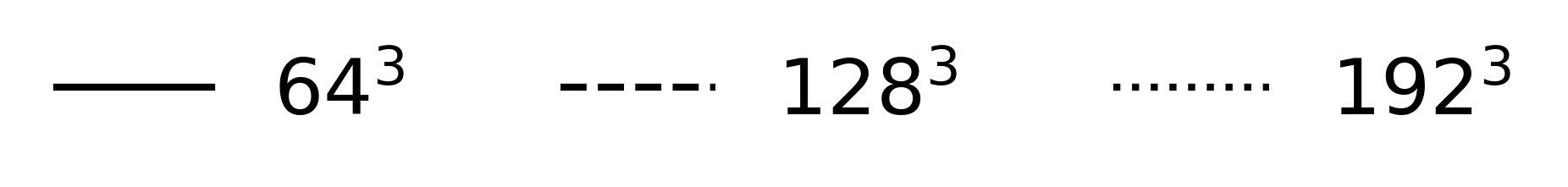}
\caption{\label{fig:viscoresi} $(a)$ Time evolution of total kinetic energy $E_K$ and magnetic energy $E_M$ (the smallest wavenumber present in the initial mixture of Fourier modes is $k_0=\sqrt{2}$). $(b)$ Time evolution of the $L_2$-error of the velocity field. The curves at linear resolution $N$ are multiplied by a factor $(N/64)^n$ for the scheme S$n$. $(c)$ Time evolution of the $L_2$-error of the magnetic field. The same rescaling as for the velocity field occurs. The time is expressed in units of the characteristic diffusion time $t_D=1/(\mu (2\pi/L)^2)$ with $L=2\pi$ the domain size.}
\end{figure}
\fi

The viscous term's implementation is validated by considering the corresponding term in the momentum equation and assuming a constant $\rho=1$:

\beq
\label{eq:visc}
\pat \Fvel = \visc (\vnabla^2 \Fvel+\frac{1}{3} \vnabla(\vnabla \cdot \Fvel)).
\eeq

\newcommand{\Nfourier}{N_F}

This equation can be solved analytically through a passage in Fourier space:

\newcommand{\ViscMat}{V}

\beq
\pat \begin{pmatrix} \FvelFx \\ \FvelFy \\ \FvelFz\end{pmatrix}=\underbrace{-\frac{\visc}{3} \begin{pmatrix} 3k^2+k_x^2 & k_xk_y & k_xk_z \\ k_xk_y & 3k^2+k_y^2 & k_yk_z \\ k_xk_z & k_yk_z & 3k^2+k_z^2  \end{pmatrix}}_{\ViscMat} \cdot \begin{pmatrix} \FvelFx \\ \FvelFy \\ \FvelFz\end{pmatrix},
\eeq
with $\FvelFm$ the Fourier coefficient of the velocity field's $\vm$-component at the corresponding wavevector $\vk$ (and $k=|\vk|$). The matrix $\ViscMat$ on the \rhs is real symmetric, and hence diagonalizable:

\beq
\ViscMat=PDP^{-1}, D=-\mu k^2 \begin{pmatrix}1 & 0 & 0 \\ 0 & 1 & 0 \\ 0 & 0 & \frac{4}{3} \end{pmatrix},P=\begin{pmatrix} -\frac{k_z}{k_x} & -\frac{k_y}{k_x} & \frac{k_x}{k_z} \\ 0 & 1 & \frac{k_y}{k_z} \\ 1 & 0 & 1 \end{pmatrix}, \quad P^{-1}=\frac{1}{k^2}\begin{pmatrix} -k_xk_z & -k_yk_z & k_x^2+k_y^2 \\ -k_xk_y & k_x^2+k_z^2 & -k_yk_z \\ k_xk_z & k_yk_z & k_z^2 \end{pmatrix}.
\eeq

The vector $\hat{\vw}_\vk=P^{-1}\hat{\vv}_\vk$ is governed by $\pat \hat{\vw}_{\vk}=D\hat{\vw}_{\vk}$, thus:
\begin{equation}
\hat{\vw}_{\vk}(t)=\begin{pmatrix}\hat{w}_{\vk,x}(0)\exp(-\mu k^2t) \\ \hat{w}_{\vk,y}(0)\exp(-\mu k^2t) \\ \hat{w}_{\vk,z}(0)\exp(-4\mu k^2t/3) \end{pmatrix}.
\end{equation}

From it, one deduces $\hat{\vv}_\vk(t)=P\hat{\vw}_\vk(t)$.

For the resistive term, the principle is the same, with the simpler equation:
\beq
\label{eq:resi}
\pat \Fmag = \resis \vnabla^2 \Fmag,
\eeq

whose solution in Fourier space is:
\beq
\hat{\vb}_\vk(t) = \hat{\vb}_\vk(0)\exp(-\resis k^2 t).
\eeq

Thus, the viscous and resistive terms are tested by considering the time evolution of a linear superposition of $\Nfourier=30$ random Fourier modes for the velocity and magnetic fields in a triply periodic domain $[0,2\pi]^3$. The modes $(k_x,k_y,k_z)$ are drawn uniformly in the integer interval $[-5..5]^3$ and are used for both the velocity and magnetic fields in order to ease the comparison of their respective time evolution. For the velocity field, the Fourier coefficients of the three components $\FvelFx, \FvelFy$ and $\FvelFz$ receive a random amplitude and phase. For the magnetic field, the coefficients are chosen with random amplitude and phase but respecting two constraints: $(i)$ they assure the magnetic field solenoidality (i.e. $k_x \FmagFx+k_y \FmagFy+k_z\FmagFz=0$) and $(ii)$ they contain the same amount of energy as the corresponding mode for the velocity field (i.e. $|\FmagFx|^2+|\FmagFy|^2+|\FmagFz|^2=|\FvelFx|^2+|\FvelFy|^2+|\FvelFz|^2$), so as to ease comparison with the latter. When all the Fourier coefficients are determined, they are rescaled such that the total kinetic energy and total magnetic energy are initially equal to unity.


The equations \eqref{eq:visc} and \eqref{eq:resi} are solved in configuration space using the finite-volume and the constrained-transport approach respectively. The timestep is limited according to the value of $\CVISCO$ given in \tab{tab:CFL}. The viscosity and magnetic diffusivity are set to $\visc=\resis=10^{-4}$. To measure the convergence of errors, the time evolutions given by the numerical schemes are compared to the analytical one at several resolutions.

The maximum speed of propagation used in the Riemann solver, $\maxspeedS$ in \eq{eq:fluxLLF}, is taken equal to the characteristic diffusion speed, e.g. in the $\vx$-direction: $\maxspeedS^x=(4/3)\visc/\Dx$. Similarly, for the 2D Riemann solver at a cell's edge along $\vz$ (see \eq{eq:CTmaxspeed}): $S=\resis/\min(\Dx,\Dy)$. Strictly speaking, the diffusive terms have a parabolic character: the characteristic manifold does not extend along the temporal dimension so that the information propagates at infinite speed. Nevertheless, the Laplacian kernel (Green's function) decreases very rapidly (exponentially with the square of the distance) so that a finite diffusion speed is a reasonable approximation. When the diffusive term is absent in the Riemann solver ($\maxspeedS=S=0$), small-scale fluctuations resulting from numerical inacurracies are amplified, leading to an unstable scheme.

Both the initialization of the fields and the comparison with the analytical solution require the transformation of the fields defined in Fourier space to configuration space. This has to be done in a way consistent with the chosen discretization. Assuming the coefficients $\hat{\vv}_\vk$ in Fourier space correspond to those at the collocation points $(x_i,y_j,z_k)$, the Fourier coefficients of the volume averages are given by:

\beq
\hat{\vol{\Fvel}}_\vk=\hat{\Fvel}_\vk \cdot \frac{(e^{ik_x\frac{\Dx}{2}}-e^{ik_x\frac{\Dx}{2}})(e^{ik_y\frac{\Dy}{2}}-e^{ik_y\frac{\Dy}{2}})(e^{ik_z\frac{\Dz}{2}}-e^{ik_z\frac{\Dz}{2}})}{\Dx\Dy\Dz(ik_x)(ik_y)(ik_z)}.
\eeq 

Similarly, if $\hat{\vb}_\vk$ corresponds to the collocation points $(x_i,y_j,z_k)$, then

\beq
\hat{\area{\FmagL}{x}{}{}}_{x,\vk}=\hat{\FmagL}_{x,\vk} \cdot \frac{e^{-ik_x\frac{\Dx}{2}}(e^{ik_y\frac{\Dy}{2}}-e^{ik_y\frac{\Dy}{2}})(e^{ik_z\frac{\Dz}{2}}-e^{ik_z\frac{\Dz}{2}})}{\Dy\Dz(ik_y)(ik_z)}
\eeq

is the Fourier coefficient of the staggered area-averaged magnetic field along the $\vx$-direction, and in an analogous manner for the $\vy$- and $\vz$-directions. This procedure is the same as the one needed to ensure that the electromotive forcing does not introduce non-solenoidal components to the magnetic field (cf. \eq{eq:ftOUbF}).

The results are shown in \fig{fig:viscoresi}. Subfigure $(a)$ presents the time evolution of both the total kinetic energy $E_K$ and magnetic energy $E_M$, for all the schemes and at all the resolutions. The shortest wavevector of the random mixture considered is $\vk_0=(1,1,0)$, with norm squared $k_0^2=2$. This is why the decay of both magnetic and kinetic energy goes asymptotically as $\exp(-2k^2_0\mu t)$ (we have $\resis=\visc$). The decay is initially faster since the larger wavevectors start with an amount of energy comparable to that at $\vk_0$. Because of the $\nabla (\nabla \cdot \Fvel)$ term in \eq{eq:visc}, the compressive part of the velocity field (along $\vk$ in Fourier space) decays faster, so that $E_K<E_M$. 

The numerical time evolution follows very closely the exact one, even at the lowest accuracy (second-order scheme at resolution $64^3$). \SFigs{fig:viscoresi}{(b)}{(c)} show the time evolution of the error (as compared to the analytical exact solution) for the velocity and magnetic fields respectively. For the numerical scheme S$n$, $n\in\{2,4,6,8,10\}$, the curve at linear resolution $N$ is multiplied by a factor $(N/64)^n$. Doing so leads to a collapse of the curves, confirming that the schemes converge at their expected order. The error is low enough for the S10 scheme to be impacted by the machine precision at resolution $192^3$ (its curves are rescaled by a factor $3^{10} \approx 6 \times 10^4$). The error in the magnetic field components stays constant at a level close to $5\times 10^{-15}$, contrary to the one of the velocity field components, which decreases with time. This difference in behaviour between the magnetic and velocity fields is due to how floating-point arithmetics are done in a computing core. Such considerations are outside the scope of the present work.


\subsection{Forced turbulence in a statistically stationary state}
\label{sec:ssturb}

\newcommand{\EKineq}{{\cal E}_{K,eq}}

\ifx\showpng\undefined
\else
\begin{figure}[h]
\includegraphics[width=\textwidth]{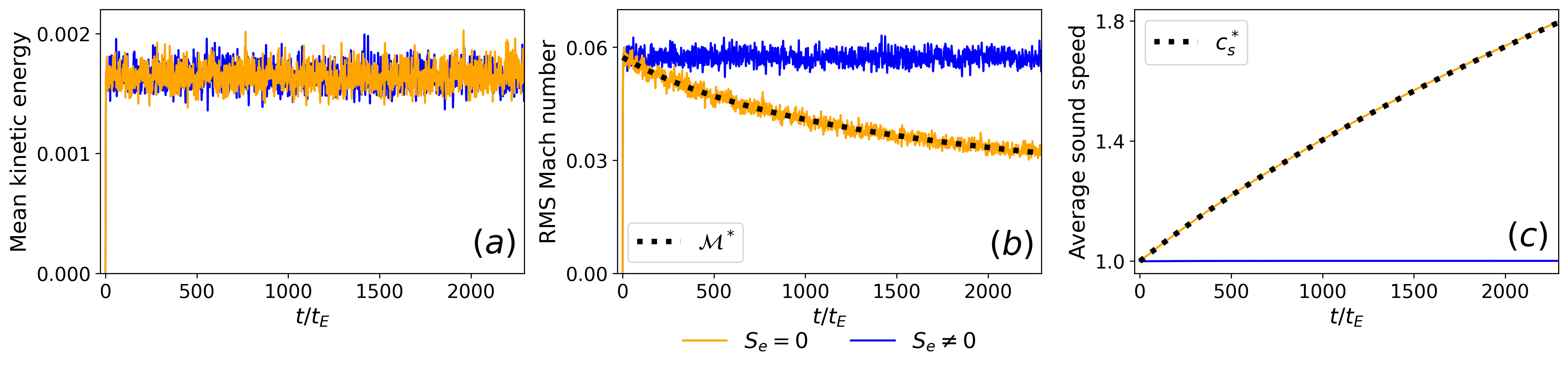}
\caption{\label{fig:impEsink} Time evolution of $(a)$ the total kinetic energy, $(b)$ the average sound speed and $(c)$ the RMS Mach number. The dotted curves are theoretical predictions given in \eqsa{eq:estimcs}{eq:estimM}. The time is given in units of the large-eddy turnover time estimated as $t_E=L/(2\sqrt{2\EKineq/\rho_0})$ with $L=1$.}
\end{figure}
\fi

As an example of application, a turbulent statistically stationary state is considered. It is obtained starting with a constant density monoatomic gas ($\rho=\rho_0=1$, $\gamma=5/3$) at rest ($\Fvel=0$) in the triply periodic computational domain $[0,1]^3$, at different resolutions. The total energy density, initially equal to the internal energy $U=p/(\gamma-1)$, is chosen so that the initial sound speed $\cs=\sqrt{\gamma p/\rho_0}$ is constant and equal to 1.

The driving terms (section \ref{sec:forcing}) inject kinetic and magnetic energies at large scales which cascade successively to smaller and smaller scales until dissipation (both of physical and numerical nature) dominates. The balance between the large scale injection and the dissipation leads to a turbulent statistically stationary state.

In the finite-volume framework, the kinetic and magnetic energies dissipated, either through the viscous and resistive terms or through numerical effects in the momentum and induction equations, are automatically transformed into heat. Indeed, the total energy density $\iiint \edens=\iiint(\frac{1}{2}\rho\Fvel^2+\frac{1}{2}\Fmag^2+U)$ is conserved down to machine precision, so that a loss of kinetic or magnetic energy means a raise in the internal energy $U$. Thus, one needs an internal energy sink $\Usink$ as well to reach a statistically stationary state. 

The importance of the internal energy sink is illustrated in section \ref{sec:importanceEsink}. Sections \ref{sec:ssturbHYDRO} and \ref{sec:ssturbMHD} present turbulent data slices and Fourier spectra during the statistically stationary state which confirm that higher-order schemes allow to resolve finer structures at a given resolution.

\subsubsection{Importance of the internal energy sink}
\label{sec:importanceEsink}

\ifx\showpng\undefined
\else
\begin{figure}[h]
\begin{minipage}[b]{0.49\linewidth}
\centering
\includegraphics[width=\textwidth]{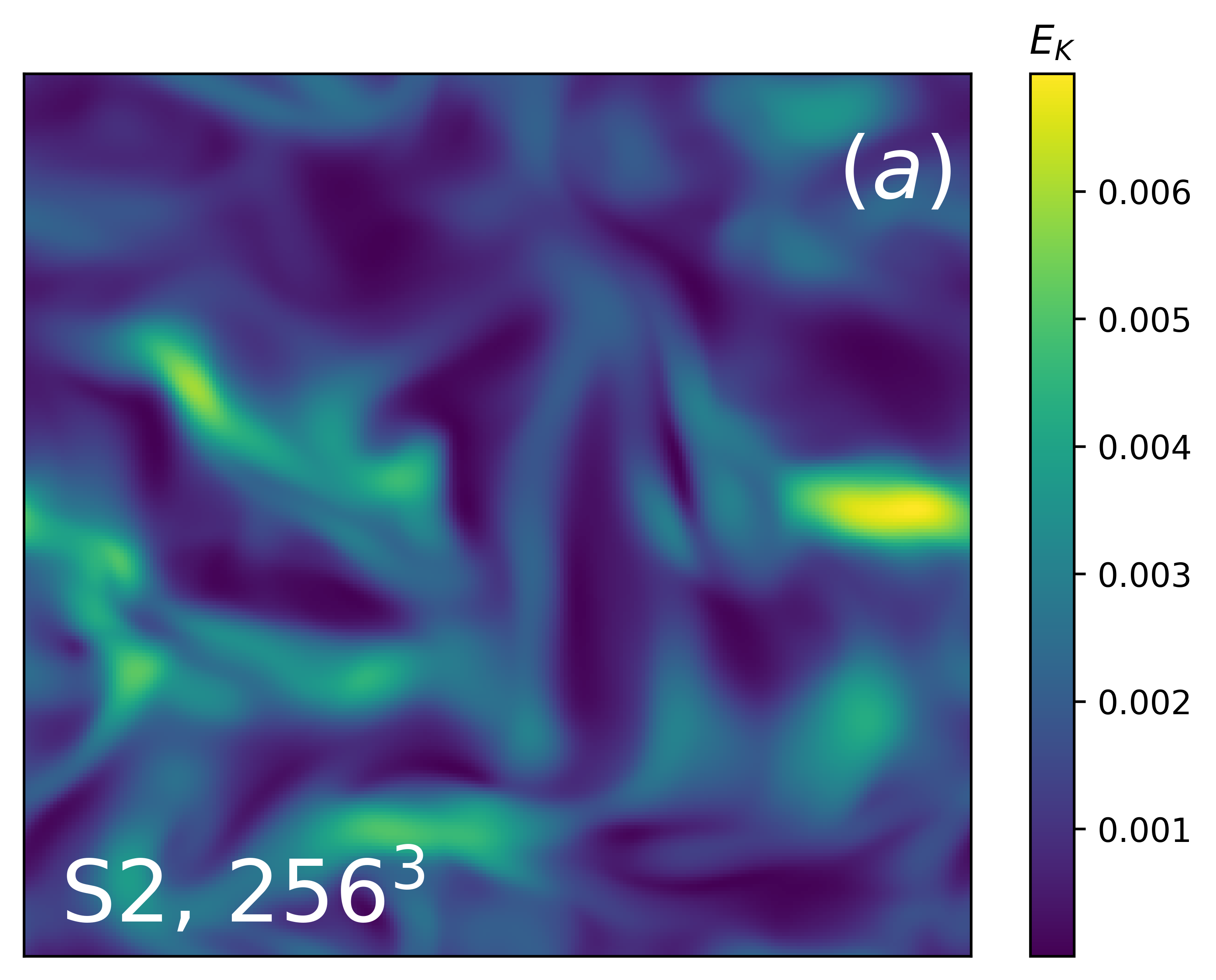}
\end{minipage}%
\begin{minipage}[b]{0.49\linewidth}
\centering
\includegraphics[width=\textwidth]{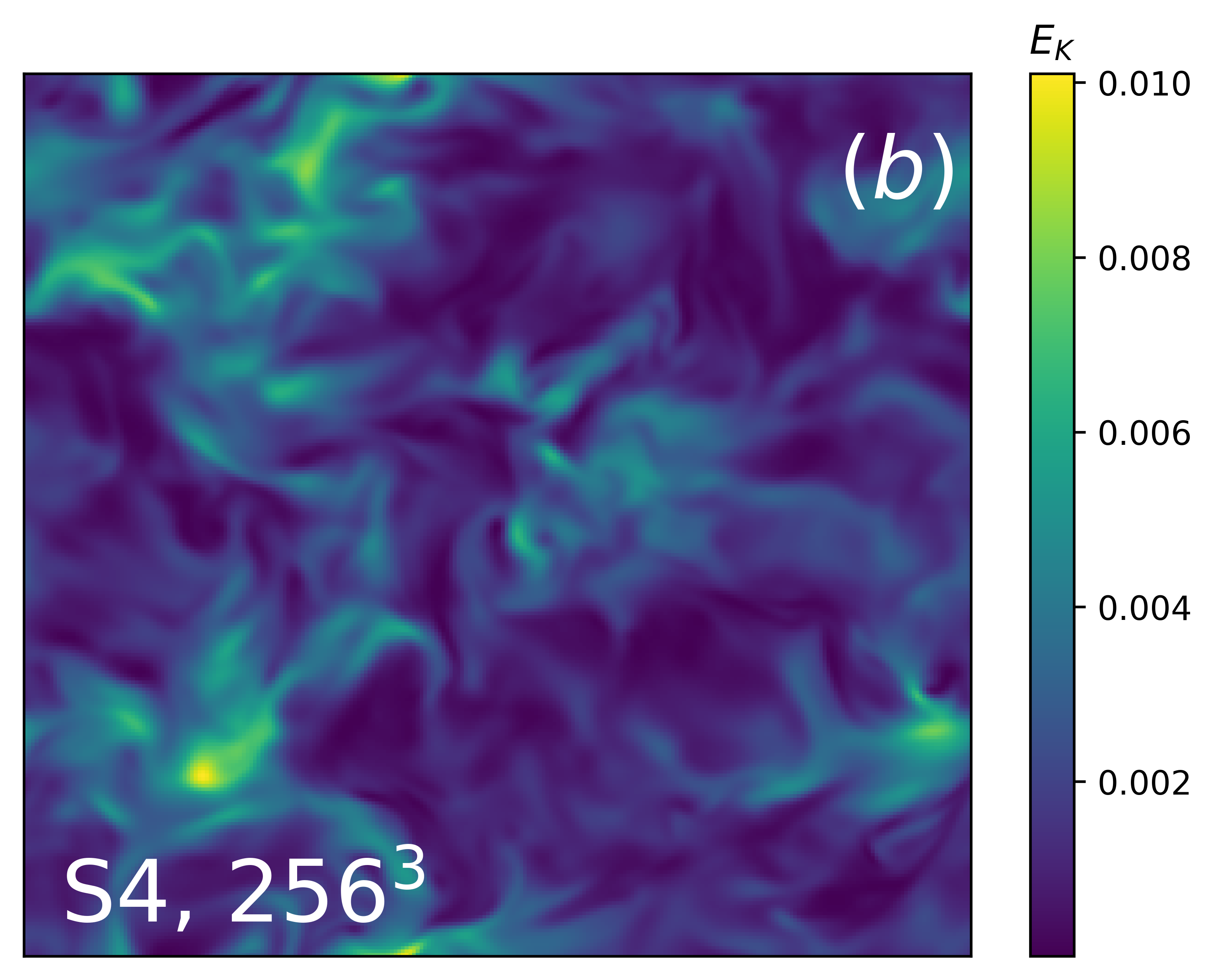}
\end{minipage}

\begin{minipage}[b]{0.49\linewidth}
\centering
\includegraphics[width=\textwidth]{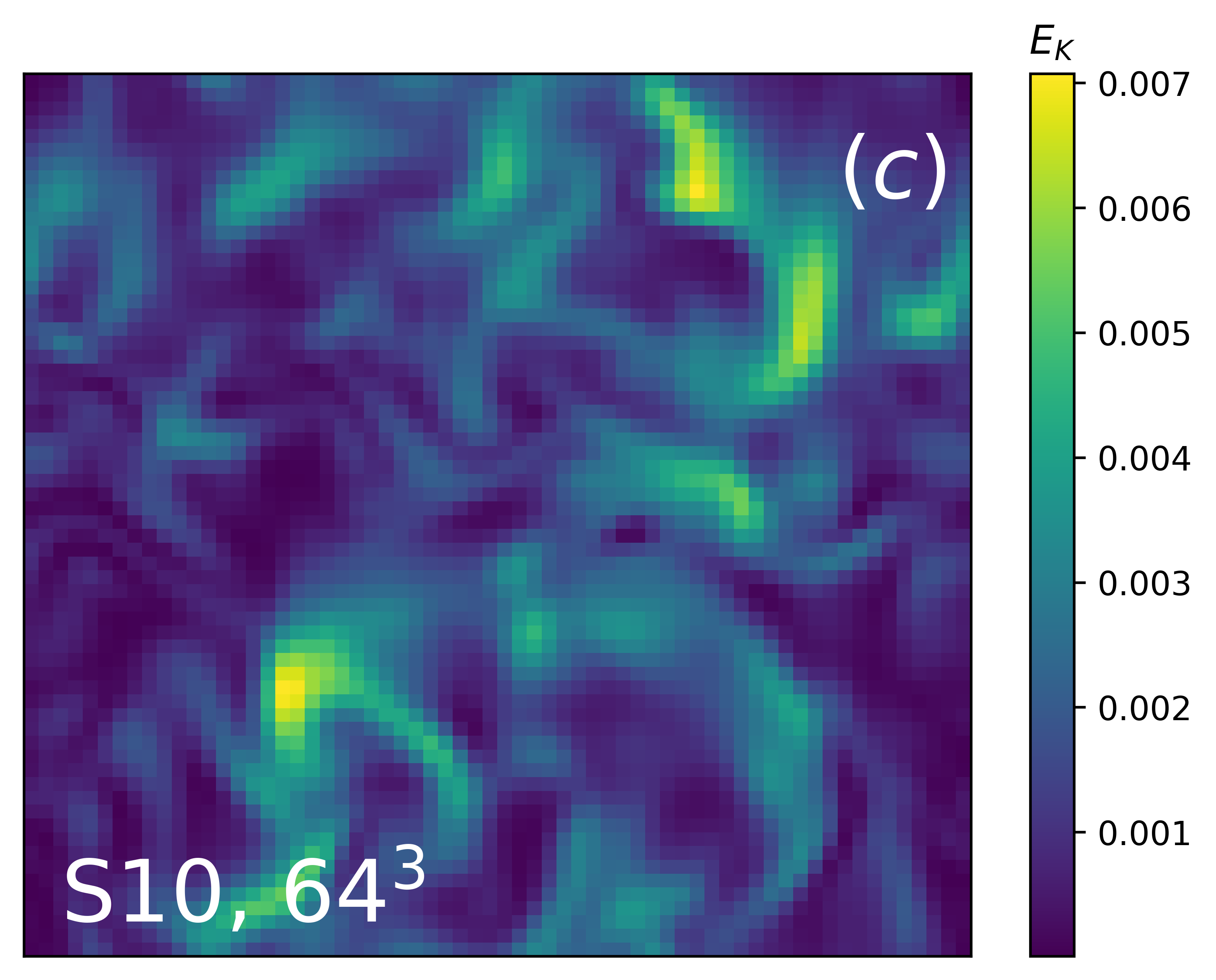}
\end{minipage}%
\begin{minipage}[b]{0.49\linewidth}
\centering
\includegraphics[width=\textwidth]{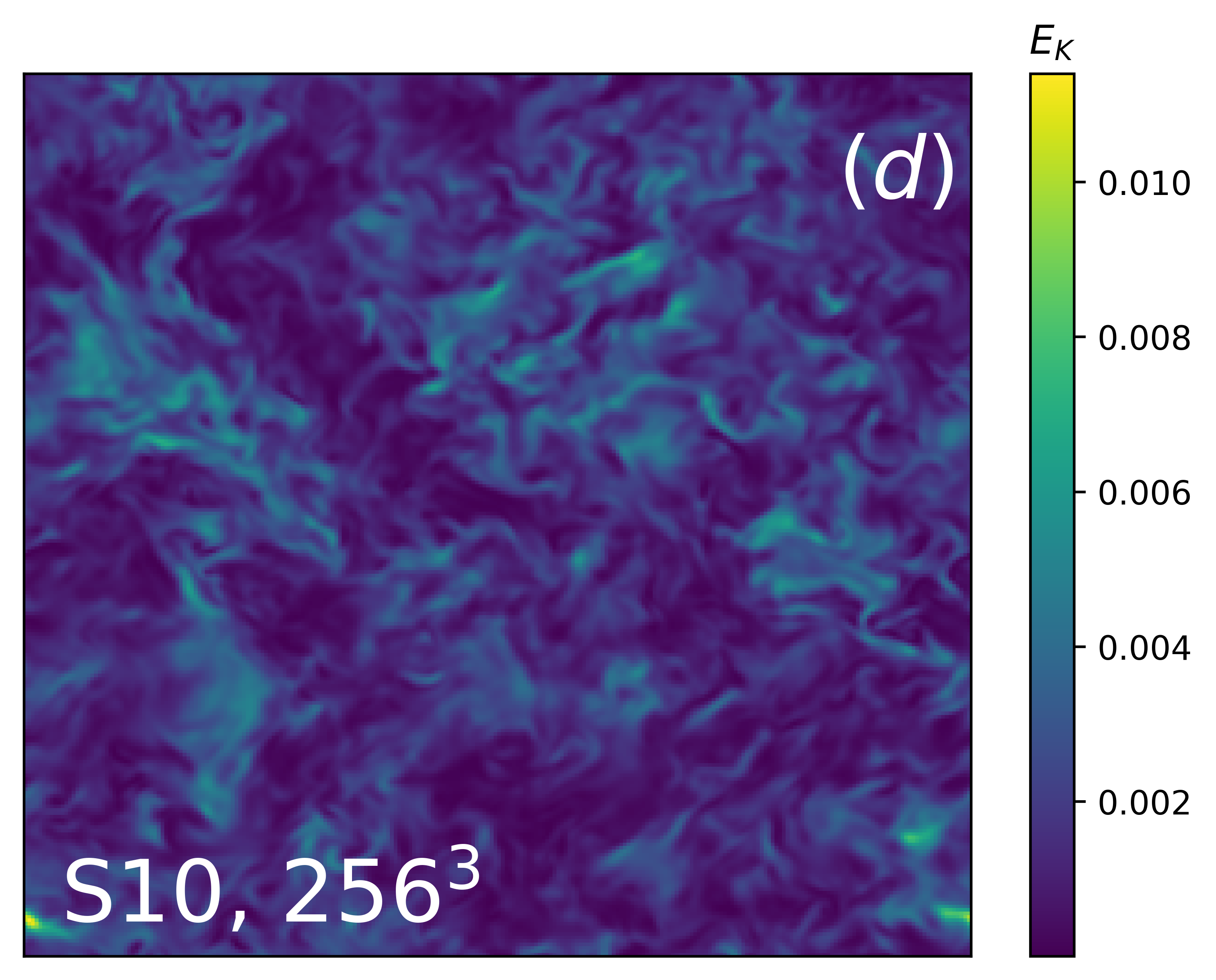}
\end{minipage}
\caption{\label{fig:slicesH} Slices of the kinetic energy during the hydrodynamic statistically stationary state.}
\end{figure}
\fi

\Fig{fig:impEsink} shows the time evolution of the total kinetic energy, the average sound speed and the root mean square (RMS) Mach number of two hydrodynamic runs executed with the S4 scheme at resolution $32^3$. One run is performed without internal energy sink $S_e$, and the other with an energy sink calibrated to keep the temperature constant. The physical viscosity is set to zero. The kinetic driving spectral profile (see \eq{eq:dtfOUF}) is taken as $\OUspecprof(\vk)=1$ for $1\leq |\vk| < 3$, zero otherwise and the energy injection rate is $\EKinj=10^{-4}$.

Because of numerical dissipation, the total kinetic energy stabilizes at a value fluctuating around $\EKineq=1.64\times 10^{-3}$ for both runs (\Sfig{fig:impEsink}{(a)}). However, the run without internal energy sink experiences an increase in temperature which translates into a decrease of the Mach number over time (\Sfig{fig:impEsink}{(c)}). The shape of the curve is expected. Indeed, the total energy in the system is governed by:

\newcommand{\totEdens}{{\cal E}_T}
\newcommand{\manRMS}{{\cal M}}

\beq
\totEdens(t)=\iiint \edens(t)=U(t=0)+\EKinj t.
\eeq

Which gives an estimate of the mean pressure in the system, once the kinetic energy begins to fluctuate around its equilibrium value $\EKineq$:

\beq
	p^*(t)=(\gamma-1)(\edens-\EKineq)\approx (\gamma-1)(U_0+\EKinj t-\EKineq).
\eeq

At low Mach numbers, the mass density remains approximately constant $\rho\approx \rho_0$ so that one can estimate the mean sound speed as:

\beq
\label{eq:estimcs}
\cs^*(t)=\sqrt{\gamma p^*(t)/\rho_0},
\eeq

and hence the RMS Mach number should go as:

\beq
\label{eq:estimM}
\manRMS^*(t)=\sqrt{2\EKineq/\rho_0}/\cs^*(t).
\eeq

\begin{figure}[h]
\includegraphics[width=\textwidth]{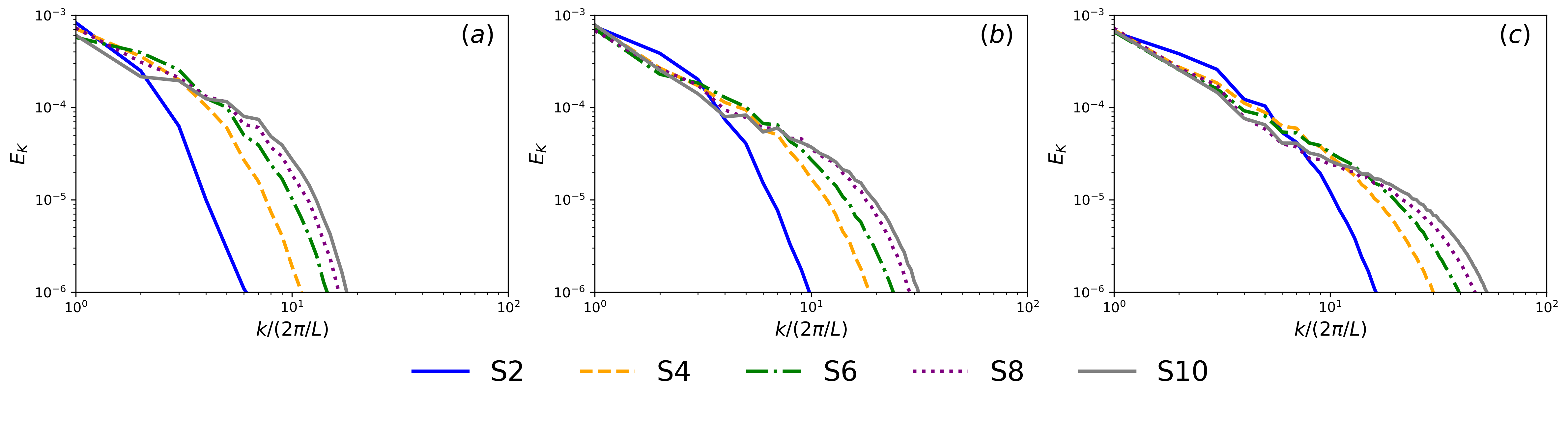}
\caption{\label{fig:specsH} Hydrodynamic turbulent statistically stationary state: velocity power spectra at resolutions $(a)$ $64^3$, $(b)$ $128^3$ and $(c)$ $256^3$.}
\end{figure}

The curves \Sfigs{fig:impEsink}{(b)}{(c)} follow these estimates very well.

In order to calibrate the internal energy sink $S_e=\UsinkF U^4$ so that the temperature stays constant (see \Sfig{fig:impEsink}{(b)}), one can estimate the internal energy variation through:

\beq
\label{eq:patUestim}
\pat U \approx \EKinj-S_e,
\eeq

which is zero for $\UsinkF=\EKinj/U^4(t=0)$.

\subsubsection{Hydrodynamic statistically stationary state}
\label{sec:ssturbHYDRO}


In this section, the same hydrodynamic setting as in section \ref{sec:importanceEsink} is considered, using the internal energy sink. The turbulent statistically stationary state obtained by schemes of different orders and at different resolutions are compared with one another.

\Fig{fig:slicesH} compares slices of the kinetic energy, obtained by the schemes S2, S4 and S10 at resolution $256^3$. The fourth-order scheme displays significantly finer structures as compared to the second-order one, but not nearly as fine as the tenth-order one. \SFig{fig:slicesH}{(c)} shows that the level of detail for the S10 scheme at the coarse $64^3$ resolution is similar to the S2 scheme's one at the higher $256^3$ resolution.

This visual impression is confirmed through the velocity power spectra (\fig{fig:specsH}). The extent of the inertial range, where the power spectrum goes as $k^m$ with $m\approx -5/3$ consistent with Kolmogorov's phenomenology, is significantly broader for the fourth and tenth-order schemes. The inertial range for the S10 scheme at resolution $64^3$ displays an extent between the one of the S2 and the one of the S4 schemes at resolution $256^3$.

This example of application illustrates that higher-order schemes can indeed be very beneficial when they allow to reach a certain results' accuracy at a significantly lower resolution as compared to, e.g., second-order schemes. 


\subsubsection{MHD statistically stationary state}
\label{sec:ssturbMHD}

\ifx\showpng\undefined
\else
\begin{figure}[h]
\includegraphics[width=\textwidth]{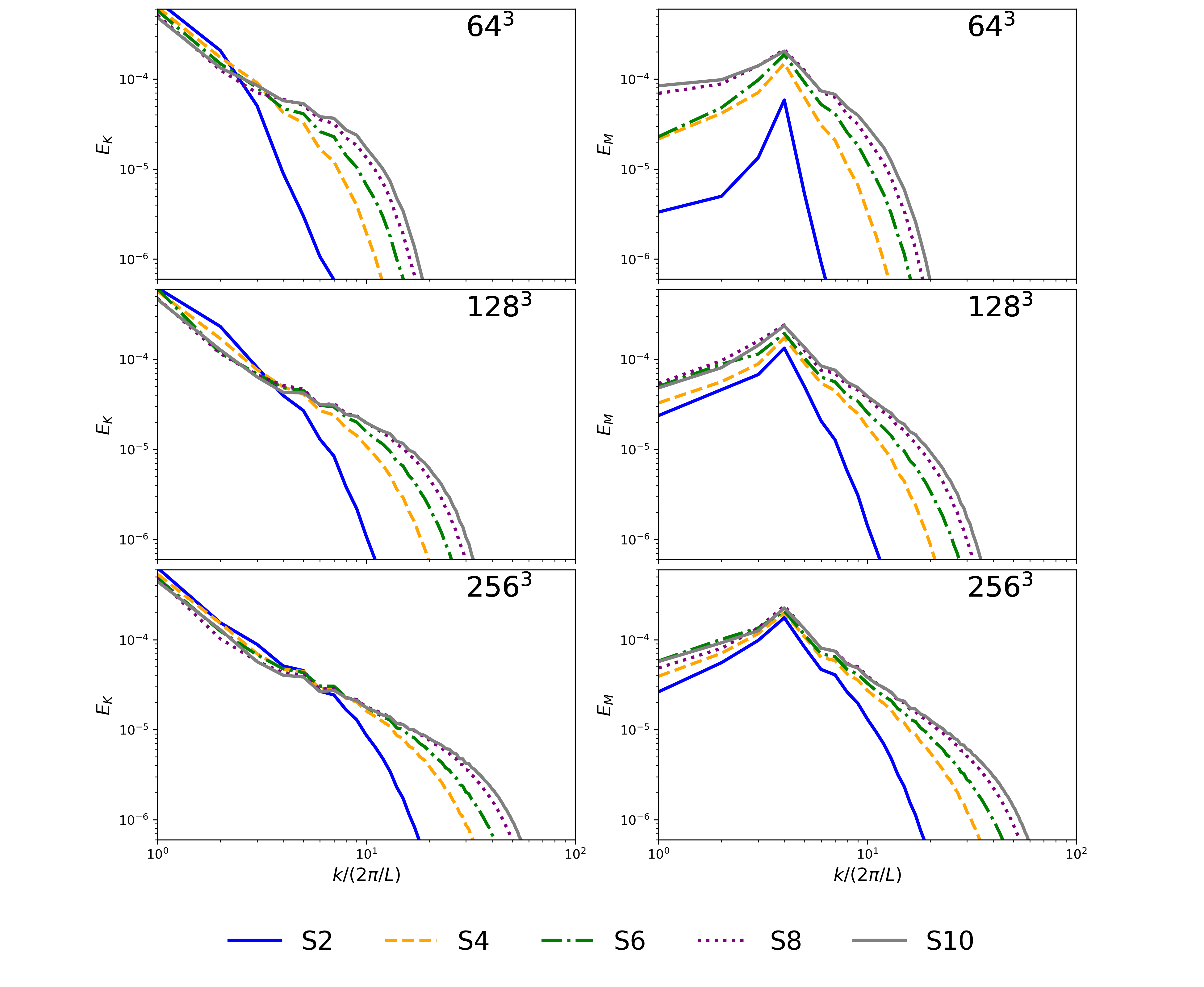}
\caption{\label{fig:MHDspecs} Fourier power spectra of the velocity (left) and magnetic fields (right) during the MHD turbulent statistically stationary state at different resolutions.}
\end{figure}
\fi

The MHD statistically stationary state is obtained using the same setting as in the hydrodynamic case, with the addition of an electromotive driving. It injects magnetic energy at the same rate as the mechanical forcing: $\EMinj=\EKinj=10^{-4}$. The spectral profile (see \eq{eq:dtfOUF}) for the electromotive forcing is $\OUspecprof_B(\vk)=1$ for $4\leq |\vk| < 5$, and zero otherwise. The electromotive driving occurs at smaller spatial scales as compared to the mechanical driving so that enough modes are excited. If this is not the case, the lack of degrees of freedom can lead to an accumulation of cross-helicity (alignment between $\Fvel$ and $\Fmag$) which affects the dynamics.

In order to maintain the temperature at its initial value, the factor of the internal energy sink is set to $\UsinkF=(\EMinj+\EKinj)/U^4(t=0)$ (cf. \eq{eq:patUestim}, considering the magnetic energy injection as well).

As in the hydrodynamic case, a higher-order scheme reveals significantly finer structures at a given resolution (slices not shown). One may notice that even though the range with an approximate power-law scaling steadily grows with increasing order of accuracy, the growth is not as significant when moving from eighth to tenth order, as compared to the changes S2$\to$S4, S4$\to$S6 and S6$\to$S8.


\section{Conclusion}
\label{sec:conclusion}

This work presents a computationally efficient method to implement high-order finite-volume compressible MHD solvers. The solvers use the constrained-transport approach to keep the magnetic field solenoidal. The key ingredient is a passage through point values, as proposed in \cite{BUH14,COC11}: the area-averages computed by a dimension-by-dimension WENO algorithm are transformed to point values in the middle of the faces, which are used to compute point-valued interfacial fluxes, finally transformed back as area-averages. We supply the required explicit \formulas up to tenth order of accuracy. This method requires only one reconstruction per face for any order of accuracy, contrary to other quadrature methods which may require an increasing number of values at each face with increasing order of accuracy.

The consistent inclusion of viscous and resistive terms respecting the finite-volume and constrained-transport formalism and preserving the order of accuracy has been shown. The usage of volume-average$\leftrightarrow$point value transformations allows to handle cooling functions and forcing terms in a high-order manner as well.

The numerical method is validated by several tests, including the advection of a 3D MHD vortex. The numerical dissipation decreases strongly when using higher-order numerics. Even though high-order schemes are more expensive at a given resolution, the results' accuracy is comparable to that of lower order schemes at a higher resolution. This trend is also visible in the fine structures resolved by higher-order schemes in turbulent systems.

In conclusion, this work shows that extending existing codes to higher-order can be beneficial and reduce computing time for given accuracy requirements. For applications in the case of non-smooth supersonic flows, common in astrophysical contexts, the oscillating behaviour of higher-order polynomials near strong gradients leads to stability issues. A solution may be the use of some ``flattening'' (local reduction of the reconstruction order in the vicinity of strong gradients \cite{VTH19,COW84,BAL12}), possibly associated with positivity preserving schemes \cite{WSH19}.



\section*{Acknowledgments}
The authors gratefully acknowledge G. Leidi for illuminating discussions. They gratefully acknowledge as well the computing time made available to them principally by the Max Planck Computing and Data Facility (MPCDF), but also on the high-performance computers HLRN-IV at GWDG at the NHR Centers NHR@G\"ottingen. These NHR Centers are jointly supported by the Federal Ministry of Education and Research and the state governments participating in the NHR (www.nhr-verein.de/unsere-partner).

\section*{Data Availability}
The data that support the findings of this study are available from the corresponding author upon reasonable request.

\section*{Conflicts of interest}
The authors declare that they have no known competing financial interests or personal relationships that could have appeared to influence the work reported in this paper.

\section*{Author Contributions}
JMT: Writing - original draft, Conceptualization, Implementation and Numerical tests. WCM: Conceptualization, Writing-review and editing.

\appendix

\newcommand{\lqmV}{\linea{q}{x}{i-5}}
\newcommand{\lqmIV}{\linea{q}{x}{i-4}}
\newcommand{\lqmIII}{\linea{q}{x}{i-3}}
\newcommand{\lqmII}{\linea{q}{x}{i-2}}
\newcommand{\lqmI}{\linea{q}{x}{i-1}}
\newcommand{\lqz}{\linea{q}{x}{i}}
\newcommand{\lqpI}{\linea{q}{x}{i+1}}
\newcommand{\lqpII}{\linea{q}{x}{i+2}}
\newcommand{\lqpIII}{\linea{q}{x}{i+3}}
\newcommand{\lqpIV}{\linea{q}{x}{i+4}}
\newcommand{\lqpV}{\linea{q}{x}{i+5}}

\newcommand{\inddx}{m_x}
\newcommand{\inddy}{m_y}
\newcommand{\inddz}{m_z}

\newcommand{\Di}{\alpha}
\newcommand{\Dj}{\beta}
\newcommand{\Dk}{\gamma}
\newcommand{\sDi}{s_\alpha}
\newcommand{\sDj}{s_\beta}
\newcommand{\sDk}{s_\gamma}

\newcommand{\LinA}{A}
\newcommand{\supAz}{0}
\newcommand{\supWx}{\presup x1}
\newcommand{\supWy}{\presup y1}
\newcommand{\supWxx}{\presup x2}
\newcommand{\supWyy}{\presup y2}
\newcommand{\supWxy}{\presup 11}

\newcommand{\apresup}[1]{\Sigma^{\pm}_{#1}}
\newcommand{\cpresup}[1]{\Sigma^{\circ}_{#1}}
\newcommand{\LinC}[1]{C^{\Di\Dj\Dk,\circ}_{#1}}

\newcommand{\Diz}{\Di_0}
\newcommand{\Djz}{\Dj_0}
\newcommand{\Dkz}{\Dk_0}
\newcommand{\ensDiDjDk}{\Omega_{\Di\Dj\Dk}}
\newcommand{\ensDizDjzDkz}{\Omega_{\Diz\Djz\Dkz}}
\newcommand{\ensDiDj}{\Omega_{\Di\Dj}}

\section{Passage through point values}
\label{app:transf_form_deriv}

In this appendix, a general method to find volume-average$\lra$point value transformation \formulas which use point symmetric stencils around the considered cell is presented. This method has been sketched in \cite{BUH14,COC11}. It is explicited in more details in sections \ref{app:pform_deriv_theor_example} and \ref{app:pform_deriv_theor_general}. The obtained \formulas are provided in sections \ref{app:transf_form_vp} and \ref{app:transf_form_lp}.

\subsection{Derivation of the \formulas: an example}
\label{app:pform_deriv_theor_example}

\newcommand{\Qlini}{\linea{q}{x}{i}}
\newcommand{\QlinpDi}{\linea{q}{x}{i+\Di}}
\newcommand{\QlinmDi}{\linea{q}{x}{i-\Di}}

\newcommand{\QlinOa}{\lineasup{q}{x}{i}{\Sigma\Di}}
\newcommand{\QlinOu}{\lineasup{q}{x}{i}{\Sigma 1}}
\newcommand{\QlinOuu}{\lineasup{q}{x}{i}{\Sigma 2}}
\newcommand{\QlinOuuu}{\lineasup{q}{x}{i}{\Sigma 3}}
\newcommand{\QlinOuuuu}{\lineasup{q}{x}{i}{\Sigma 4}}

As an instructive example, this section presents the simpler cases of fourth-order and sixth-order line-average$\leftrightarrow$point value transformations. Writing $\Qlini$ the line-average of a quantity $q$ over the cell $[x_i-\Dx/2,x_i+\Dx/2]$, a Taylor expansion up to sixth-order of accuracy gives (cf. \eq{eq:simple2O}):

\beq
	\label{eq:ltoaz}	\Qlini=\int_{-\Dx/2}^{\Dx/2} q(x_i+\eps)\deps=\QptO+\frac{1}{24}\Dx^2\pax^2\QptO+\frac{1}{1920}\Dx^4\pax^4\QptO+\bigO(\Dx^6),
\eeq

with $q_i=q(x_i)$. In order to have a sixth-order approximation of $q_i$ as a function of the known line-averages (that is, find a formula for a line-average$\to$point value transformation), one needs a way to eliminate the terms in $\pax^2\QptO$ and $\pax^4\QptO$. This can be achieved by considering $\QlinOa=\QlinpDi+\QlinmDi$, which gives, after expansion of the point-value derivatives $\pax^{2\inddx} \Qpt_{i\pm\Di}$ in Taylor series around $x_i$:

\beqa
\nonumber	\QlinOa&=&\Big(2\Qpt_{i}+(\Di\Dx)^2\pax^2 \QptO+\frac{1}{12}(\Di\Dx)^4\pax^4\QptO+\bigO(\Dx^6)\Big)+\frac{1}{24}\Dx^2\pax^2\Big(2\Qpt_i+(\Di\Dx)^2\pax^2\Qpt_i+\bigO(\Dx^4)\Big)\\
\label{eq:TexpQavgOa}	&+&\frac{1}{1920}\Dx^4\pax^4\Big(2\Qpt_i+\bigO(\Dx^2)\Big)+\bigO(\Dx^6).
\eeqa

Hence, for $\Di \in \{1,2\}$:

\beqa
	\label{eq:ltoau}
	\QlinOu&=&2\Qpt_i+\frac{13}{12}\Dx^2\pax^2\Qpt_i+\frac{121}{960}\Dx^4\pax^4\Qpt_i+\bigO(\Dx^6),\\
	\label{eq:ltoauu}
	\QlinOuu&=&2\Qpt_i+\frac{49}{12}\Dx^2\pax^2\Qpt_i+\frac{1441}{960}\Dx^4\pax^4\Qpt_i+\bigO(\Dx^6).
\eeqa

	One can find three constants $A_0, A_1$ and $A_2$ such that $A_0\Qlini+A_1 \QlinOu+A_2\QlinOuu=\QptO+\bigO(\Dx^6)$ by solving the system (see \eqsaa{eq:ltoaz}{eq:ltoau}{eq:ltoauu}):

\beq
\begin{pmatrix}
1 & 2 & 2 \\
\frac{1}{24} & \frac{13}{12} & \frac{49}{12} \\
\frac{1}{1920} & \frac{121}{960} & \frac{1441}{960}
\end{pmatrix}
\begin{pmatrix}
A_0\\
A_1\\
A_2
\end{pmatrix}
=
\begin{pmatrix}
1\\
0\\
0
\end{pmatrix},
\eeq

which has the solution $(\frac{1067}{960},-\frac{29}{480},\frac{3}{640})$ (cf. eq. 19 in \cite{BUH14}).

For a fourth-order line-average$\to$point value transformation, it is enough to remove the $\pax^2q_i$ term. This can be done by using \eqsa{eq:ltoaz}{eq:ltoau} to find $B_0, B_1$ such that $B_0\Qlini+B_1 \QlinOu=\QptO+\bigO(\Dx^4)$. The linear system to solve consists of a submatrix of the previous one:

\beq
\begin{pmatrix}
1 & 2 \\
\frac{1}{24} & \frac{13}{12}
\end{pmatrix}
\begin{pmatrix}
B_0\\
B_1
\end{pmatrix}
=
\begin{pmatrix}
1\\
0
\end{pmatrix},
\eeq

which has the solution $(\frac{13}{12},-\frac{1}{24})$ (cf. eq. 17 in \cite{BUH14}).

For point value$\to$line-average transformations, one can first rewrite the point values as a function of the line-averages only. Starting from \eq{eq:ltoaz}:

\beq
	\QptO=\Qlini-A\Dx^2\pax^2\QptO-B\Dx^4\pax^4\QptO+\bigO(\Dx^6),
\eeq

with $A=\frac{1}{24}$ and $B=\frac{1}{1920}$, one can deduce recursively:

\beqa
	\nonumber \QptO&=&\Qlini-A\Dx^2\pax^2(\Qlini-A\Dx^2\pax^2\QptO+\bigO(\Dx^4))-B\Dx^4\pax^4(\Qlini+\bigO(\Dx^2))+\bigO(\Dx^6),\\
	\nonumber &=&\Qlini-A\Dx^2\pax^2\Qlini+(A^2-B)\Dx^4\pax^4\Qlini+\bigO(\Dx^6),\\
	\label{eq:recursivQpt}&=&\Qlini-\frac{1}{24}\Dx^2\pax^2\Qlini+\frac{7}{5760}\Dx^4\pax^4\Qlini+\bigO(\Dx^6).
\eeqa

The principle is then the same as above: a Taylor expansion of $q_i^{\Sigma\alpha}=q_{i+\alpha}+q_{i-\alpha}$ in terms of $\Qlini$, $\pax^2\Qlini$ and $\pax^4\Qlini$ for $\alpha\in\{1,2\}$ gives three expressions. Combining them appropriately, one can find three constants $C_0, C_1$ and $C_2$ such that $C_0q_i+C_1q_i^{\Sigma 1}+C_2q_i^{\Sigma 2}=\Qlini+\bigO(\Dx^6)$.

\newcommand{\intseg}[2]{\int_{#1-#2/2}^{#1+#2/2}}
\newcommand{\intsym}[1]{\int_{-#1}^{#1}}
\newcommand{\epsx}{\epsilon_x}
\newcommand{\epsy}{\epsilon_y}
\newcommand{\epsz}{\epsilon_z}
\newcommand{\depsx}{d\epsx}
\newcommand{\depsy}{d\epsy}
\newcommand{\depsz}{d\epsz}
\newcommand{\dvoleps}{\depsx \depsy \depsz}
\newcommand{\indx}{n_x}
\newcommand{\indy}{n_y}
\newcommand{\indz}{n_z}
\newcommand{\suminf}[1]{\sum_{#1=0}^\infty}

\subsection{Derivation of the \formulas: general case}
\label{app:pform_deriv_theor_general}
For the general three-dimensional case, a Taylor expansion of the volume-average of $q$ over the cell $\Cellijk=[x_i-\Dx/2,x_i+\Dx/2]\times[y_j-\Dy/2,y_j+\Dy/2]\times[z_k-\Dz/2,z_k+\Dz/2]$ gives:

\beq
	\QvolT=\frac{1}{\Dx\Dy\Dz}\suminf{\indx}\suminf{\indy}\suminf{\indz}\pax^{\indx}\pay^{\indy}\paz^{\indz} \QptT \frac{1}{\indx!\indy!\indz!}\iiint_{\Cellijk} \epsx^{\indx}\epsy^{\indy}\epsz^{\indz} \dvoleps.
\eeq

The integral is non-zero only for even $\indx, \indy, \indz$. Thus, for $\indx=2\inddx, \indy=2\inddy, \indz=2\inddz$ with $\inddx, \inddy, \inddz \in \mathbb{N}$:

\beq
	\label{eq:QavgFctQpt}
	\QvolT=\suminf{\inddx}\suminf{\inddy}\suminf{\inddz} \frac{\Dx^{2\inddx}\Dy^{2\inddy}\Dz^{2\inddz}\pax^{2\inddx}\pay^{2\inddy}\paz^{2\inddz} \QptT}{(2\inddx+1)!(2\inddy+1)!(2\inddz+1)!2^{2(\inddx+\inddy+\inddz)}}.
\eeq

\newcommand{\accOrd}{p}

An expression up to an order of accuracy $2\accOrd$ contains thus all derivatives $\pax^{2\inddx}\pay^{2\inddy}\paz^{2\inddz} \Qpt$ with $0\leq 2(\inddx+\inddy+\inddz)\leq 2(\accOrd-1)$. In order to eliminate these terms as in the example explicited in section \ref{app:pform_deriv_theor_example}, one can consider a Taylor expansion (in terms of derivatives of $\Qpt$) of all the terms $\Qvol^{\Sigma \Di\Dj\Dk}_{ijk}$, which consist of the sum of all volume-averages where one offset with respect to $(i,j,k)$ is $\pm\Di$, another $\pm\Dj$ and the third one $\pm\Dk$, considering each combination of offset only once (generalizing the definition of $\area{q}{z}{\Sigma mn}{ij}$ to 3D, see section \ref{sec:explsigmamn}). A concrete expression for these terms is given in section \ref{app:transf_form_vp}. With the offsets $\Di,\Dj,\Dk\geq 0$ and $0\leq \Di+\Dj+\Dk\leq\accOrd-1$, there are as many such terms as there are different derivative terms up to order $2\accOrd$, and their Taylor expansion contain all the derivative terms that appear in \eq{eq:QavgFctQpt}. Since the offsets considered are symmetric with respect to $(i,j,k)$, the odd derivative terms in the expansion cancel out. 

In order to find a volume-average$\to$ point value transformation formula of order $2\accOrd$, one can find coefficients $A_{\Di\Dj\Dk}$ such that:

\beq
\sum_{\Di,\Dj,\Dk} A_{\Di\Dj\Dk}\Qvol^{\Sigma\Di\Dj\Dk}_{ijk}=\Qpt_{ijk}+\bigO(\Dx^{2\accOrd}+\Dy^{2\accOrd}+\Dz^{2\accOrd}),
\eeq

and hence solve a system of linear equations with as many equations as unknowns.

For the reverse transformation, from point values to volume averages, one can similarly to the 1D example of appendix \ref{app:pform_deriv_theor_example} rewrite \eq{eq:QavgFctQpt} as:

\beq
	\QptT=\QvolT-\sum_{\substack{\inddx,\inddy,\inddz=0 \\ \inddx+\inddy+\inddz>0}}^{\infty}\frac{\Dx^{2\inddx}\Dy^{2\inddy}\Dz^{2\inddz}\pax^{2\inddx}\pay^{2\inddy}\paz^{2\inddz} \QptT}{(2\inddx+1)!(2\inddy+1)!(2\inddz+1)!2^{2(\inddx+\inddy+\inddz)}}.
\eeq

And recursively (cf. \eq{eq:recursivQpt}) deduce an expression of the form:

\beq
	\QptT=\QvolT-\sum_{\substack{\inddx,\inddy,\inddz=0 \\ \inddx+\inddy+\inddz>0}}^{\inddx+\inddy+\inddz=p-1} K_{\inddx,\inddy,\inddz}\pax^{2\inddx}\pay^{2\inddy}\paz^{2\inddz}\QvolT+\bigO(\Dx^{2p}+\Dy^{2p}+\Dz^{2p}).
\eeq

A Taylor expansion of the terms $\Qpt^{\Sigma \Di\Dj\Dk}_{ijk}$ (in terms of derivatives of the volume-averages $\Qvol$) gives then a linear system, whose solution delivers the coefficients for the point$\to$volume-averages transformation.

\subsection{Formulas for volume-averages$\leftrightarrow$point values transformations}

\label{app:transf_form_vp}

Applying the method described above, one obtains the following \formulas for volume-averages$\to$point values transformations up to tenth order of accuracy:

\begin{eqnarray}
	\QptT&=&\QvolT+\bigO(h^{2}),\\
	\QptT&=&\frac{5}{4}\QvolT-\frac{1}{24}\QvolTu+\bigO(h^{4}),\\
	\QptT&=&\frac{1301}{960}\QvolT-\frac{97}{1440}\QvolTu+\frac{3}{640}\QvolTuu+\frac{1}{576}\QvolTuv+\bigO(h^{6}),\\
\nonumber 	\QptT&=&\frac{341519}{241920}\QvolT-\frac{80881}{967680}\QvolTu+\frac{173}{17920}\QvolTuu-\frac{5}{7168}\QvolTuuu\\
&+&\frac{119}{34560}\QvolTuv-\frac{1}{5120}\QvolTuuv-\frac{1}{13824}\QvolTuvw+\bigO(h^{8}),\\
\nonumber 	\QptT&=&\frac{15997789}{11059200}\QvolT-\frac{1829207}{19353600}\QvolTu+\frac{25751}{1843200}\QvolTuu-\frac{65}{36864}\QvolTuuu\\
\nonumber &+&\frac{35}{294912}\QvolTuuuu+\frac{11723}{2419200}\QvolTuv-\frac{1019}{2150400}\QvolTuuv+\frac{5}{172032}\QvolTuuuv\\
&+&\frac{9}{409600}\QvolTuuvv-\frac{47}{276480}\QvolTuvw+\frac{1}{122880}\QvolTuuvw+\bigO(h^{10}),
\end{eqnarray}

with $h=\Dx+\Dy+\Dz$. The notation $\Qvol^{\Sigma \Di\Dj\Dk}_{ijk}$ is explained in section \ref{app:pform_deriv_theor_general}. Up to tenth-order of accuracy, one has either $\Dk=\Dj$ or $\Dk=0$, resulting in:

\beqa
\nonumber \Qvol^{\Sigma \Di\Dj\Dj}_{ijk}&\underset{\substack{\Di,\Dj>0 \\ \Di\neq\Dj}}{=}&\Qvol_{i+\Di,j+\Dj,k+\Dj}+\Qvol_{i+\Di,j+\Dj,k-\Dj}+\Qvol_{i+\Di,j-\Dj,k+\Dj}+\Qvol_{i+\Di,j-\Dj,k-\Dj}\\
\nonumber 	&+&\Qvol_{i-\Di,j+\Dj,k+\Dj}+\Qvol_{i-\Di,j+\Dj,k-\Dj}+\Qvol_{i-\Di,j-\Dj,k+\Dj}+\Qvol_{i-\Di,j-\Dj,k-\Dj}\\
\nonumber 	&+&\Qvol_{i+\Dj,j+\Di,k+\Dj}+\Qvol_{i+\Dj,j+\Di,k-\Dj}+\Qvol_{i+\Dj,j-\Di,k+\Dj}+\Qvol_{i+\Dj,j-\Di,k-\Dj}\\
\nonumber 	&+&\Qvol_{i-\Dj,j+\Di,k+\Dj}+\Qvol_{i-\Dj,j+\Di,k-\Dj}+\Qvol_{i-\Dj,j-\Di,k+\Dj}+\Qvol_{i-\Dj,j-\Di,k-\Dj}\\
\nonumber 	&+&\Qvol_{i+\Dj,j+\Dj,k+\Di}+\Qvol_{i+\Dj,j+\Dj,k-\Di}+\Qvol_{i+\Dj,j-\Dj,k+\Di}+\Qvol_{i+\Dj,j-\Dj,k-\Di}\\
\label{eq:explainSumDiDjDkbegin}	&+&\Qvol_{i-\Dj,j+\Dj,k+\Di}+\Qvol_{i-\Dj,j+\Dj,k-\Di}+\Qvol_{i-\Dj,j-\Dj,k+\Di}+\Qvol_{i-\Dj,j-\Dj,k-\Di},\\
\nonumber \Qvol^{\Sigma \Dj\Dj\Dj}_{ijk}&\underset{\Dj>0}{=}&\Qvol_{i+\Dj,j+\Dj,k+\Dj}+\Qvol_{i+\Dj,j+\Dj,k-\Dj}+\Qvol_{i+\Dj,j-\Dj,k+\Dj}+\Qvol_{i+\Dj,j-\Dj,k-\Dj}\\
 	&+&\Qvol_{i-\Dj,j+\Dj,k+\Dj}+\Qvol_{i-\Dj,j+\Dj,k-\Dj}+\Qvol_{i-\Dj,j-\Dj,k+\Dj}+\Qvol_{i-\Dj,j-\Dj,k-\Dj},\\
\nonumber \Qvol^{\Sigma \Di\Dj 0}_{ijk}&\underset{\substack{\Di,\Dj>0 \\ \Di\neq\Dj}}{=}&\Qvol_{i+\Di,j+\Dj,k}+\Qvol_{i+\Di,j-\Dj,k}+\Qvol_{i-\Di,j+\Dj,k}+\Qvol_{i-\Di,j-\Dj,k}\\
\nonumber	&+&\Qvol_{i+\Dj,j+\Di,k}+\Qvol_{i+\Dj,j-\Di,k}+\Qvol_{i-\Dj,j+\Di,k}+\Qvol_{i-\Dj,j-\Di,k}\\
\nonumber	&+&\Qvol_{i+\Di,j,k+\Dj}+\Qvol_{i+\Di,j,k-\Dj}+\Qvol_{i-\Di,j,k+\Dj}+\Qvol_{i-\Di,j,k-\Dj}\\
\nonumber	&+&\Qvol_{i+\Dj,j,k+\Di}+\Qvol_{i+\Dj,j,k-\Di}+\Qvol_{i-\Dj,j,k+\Di}+\Qvol_{i-\Dj,j,k-\Di}\\
\nonumber	&+&\Qvol_{i,j+\Di,k+\Dj}+\Qvol_{i,j+\Di,k-\Dj}+\Qvol_{i,j-\Di,k+\Dj}+\Qvol_{i,j-\Di,k-\Dj}\\
		&+&\Qvol_{i,j+\Dj,k+\Di}+\Qvol_{i,j+\Dj,k-\Di}+\Qvol_{i,j-\Dj,k+\Di}+\Qvol_{i,j-\Dj,k-\Di},\\
\nonumber \Qvol^{\Sigma \Dj\Dj 0}_{ijk}&\underset{\Dj>0}{=}&\Qvol_{i+\Dj,j+\Dj,k}+\Qvol_{i+\Dj,j-\Dj,k}+\Qvol_{i-\Dj,j+\Dj,k}+\Qvol_{i-\Dj,j-\Dj,k}\\
\nonumber	&+&\Qvol_{i+\Dj,j,k+\Dj}+\Qvol_{i+\Dj,j,k-\Dj}+\Qvol_{i-\Dj,j,k+\Dj}+\Qvol_{i-\Dj,j,k-\Dj}\\
		&+&\Qvol_{i,j+\Dj,k+\Dj}+\Qvol_{i,j+\Dj,k-\Dj}+\Qvol_{i,j-\Dj,k+\Dj}+\Qvol_{i,j-\Dj,k-\Dj},\\
\label{eq:explainSumDiDjDkend} \Qvol^{\Sigma \Dj 00}_{ijk}&\underset{\Dj>0}{=}&\Qvol_{i+\Dj,j,k}+\Qvol_{i-\Dj,j,k}+\Qvol_{i,j+\Dj,k}+\Qvol_{i,j-\Dj,k}+\Qvol_{i,j,k+\Dj}+\Qvol_{i,j,k-\Dj}.
\eeqa

Using the same notation for the sum of point-averages $\Qpt^{\Sigma \Di\Dj\Dk}_{ijk}$ with offsets $\pm\Di$,$\pm\Dj$ and $\pm\Dk$, the reverse point values$\to$volume-averages transformations are:

\begin{eqnarray}
	\QvolT&=&\QptT+\bigO(h^{2}),\\
	\QvolT&=&\frac{3}{4}\QptT+\frac{1}{24}\QptTu+\bigO(h^{4}),\\
	\QvolT&=&\frac{689}{960}\QptT+\frac{67}{1440}\QptTu-\frac{17}{5760}\QptTuu+\frac{1}{576}\QptTuv+\bigO(h^{6}),\\
\nonumber 	\QvolT&=&\frac{34025}{48384}\QptT+\frac{47477}{967680}\QptTu-\frac{2291}{483840}\QptTuu+\frac{367}{967680}\QptTuuu\\
&+&\frac{89}{34560}\QptTuv-\frac{17}{138240}\QptTuuv+\frac{1}{13824}\QptTuvw+\bigO(h^{8}),\\
\nonumber 	\QvolT&=&\frac{53802803}{77414400}\QptT+\frac{2939507}{58060800}\QptTu-\frac{138211}{23224320}\QptTuu+\frac{15403}{19353600}\QptTuuu\\
\nonumber &-&\frac{27859}{464486400}\QptTuuuu+\frac{5581}{1814400}\QptTuv-\frac{4691}{19353600}\QptTuuv+\frac{367}{23224320}\QptTuuuv\\
&+&\frac{289}{33177600}\QptTuuvv+\frac{37}{276480}\QptTuvw-\frac{17}{3317760}\QptTuuvw+\bigO(h^{10}).
\end{eqnarray}

\subsection{Formulas for area-averages/line-averages$\leftrightarrow$point values transformations}

The 2D transformation \formulas (area-average$\leftrightarrow$point values) are obtained by projecting the 3D formulas on a plane. All the terms with a subscript $(i\pm\Di,j\pm\Dj,k\pm\Dk)$ are replaced by terms with a subscript $(i\pm\Di,j\pm\Dj)$. After rearrangement of the terms this gives the \formulas presented in section \ref{sec:colella}, which are not repeated here.

Projecting these \formulas again on a line (that is, the subscripts $(i\pm\Di,j\pm\Dj)$ become $(i\pm\Di)$, one obtains line-average$\leftrightarrow$point value transformations, useful e.g. when solving 2D problems:

\label{app:transf_form_lp}

\begin{eqnarray}
	\QptO&=&\Qlini+\bigO(\Dx^{2})\\
	\QptO&=&\frac{13}{12}\Qlini-\frac{1}{24}\QlinOu+\bigO(\Dx^{4})\\
	\QptO&=&\frac{1067}{960}\Qlini-\frac{29}{480}\QlinOu+\frac{3}{640}\QlinOuu+\bigO(\Dx^{6})\\
	\QptO&=&\frac{30251}{26880}\Qlini-\frac{7621}{107520}\QlinOu+\frac{159}{17920}\QlinOuu-\frac{5}{7168}\QlinOuuu+\bigO(\Dx^{8})\\
 	\QptO&=&\frac{5851067}{5160960}\Qlini-\frac{100027}{1290240}\QlinOu+\frac{31471}{2580480}\QlinOuu-\frac{425}{258048}\QlinOuuu+\frac{35}{294912}\QlinOuuuu+\bigO(\Dx^{10})
\end{eqnarray}

\begin{eqnarray}
	\Qlini&=&\QptO+\bigO(\Dx^{2})\\
	\Qlini&=&\frac{11}{12}\QptO+\frac{1}{24}\QptOu+\bigO(\Dx^{4})\\
	\Qlini&=&\frac{863}{960}\QptO+\frac{77}{1440}\QptOu-\frac{17}{5760}\QptOuu+\bigO(\Dx^{6})\\
	\Qlini&=&\frac{215641}{241920}\QptO+\frac{6361}{107520}\QptOu-\frac{281}{53760}\QptOuu+\frac{367}{967680}\QptOuuu+\bigO(\Dx^{8})\\
	\Qlini&=&\frac{41208059}{46448640}\QptO+\frac{3629953}{58060800}\QptOu-\frac{801973}{116121600}\QptOuu+\frac{49879}{58060800}\QptOuuu-\frac{27859}{464486400}\QptOuuuu+\bigO(\Dx^{10})
\end{eqnarray}

	\bibliographystyle{elsarticle-num} 
	\bibliography{biblio}

\begin{thebibliography}{10}
\expandafter\ifx\csname url\endcsname\relax
  \def\url#1{\texttt{#1}}\fi
\expandafter\ifx\csname urlprefix\endcsname\relax\def\urlprefix{URL }\fi
\expandafter\ifx\csname href\endcsname\relax
  \def\href#1#2{#2} \def\path#1{#1}\fi

\bibitem{ROE81}
P.~L. Roe, Approximate {R}iemann solvers, parameter vectors, and difference
  schemes, Journal of Computational Physics 43 (1981) 357--372.
\newblock \href {https://doi.org/10.1016/0021-9991(81)90128-5}
  {\path{doi:10.1016/0021-9991(81)90128-5}}.

\bibitem{TSS94}
E.~F. Toro, M.~Spruce, W.~Speares, Restoration of the contact surface in the
  {HLL}-{R}iemann solver, Shock Waves 4 (1994) 25--34.
\newblock \href {https://doi.org/10.1007/BF01414629}
  {\path{doi:10.1007/BF01414629}}.

\bibitem{POW97}
K.~G. Powell, An approximate {R}iemann solver for magnetohydrodynamics, Upwind
  and High-Resolution Schemes (1997) 570--583\href
  {https://doi.org/10.1007/978-3-642-60543-7_23}
  {\path{doi:10.1007/978-3-642-60543-7_23}}.

\bibitem{TOR09}
E.~F. Toro, {R}iemann Solvers and Numerical Methods for Fluid Dynamics,
  Springer Berlin Heidelberg, 2009.
\newblock \href {https://doi.org/10.1007/b79761} {\path{doi:10.1007/b79761}}.

\bibitem{MRE15}
F.~Miczek, F.~K. Röpke, P.~V. Edelmann, Astronomy and Astrophysics\href
  {https://doi.org/10.1051/0004-6361/201425059}
  {\path{doi:10.1051/0004-6361/201425059}}.

\bibitem{MIM21}
T.~Minoshima, T.~Miyoshi, A low-dissipation {HLLD} approximate {R}iemann solver
  for a very wide range of {M}ach numbers, Journal of Computational Physics 446
  (2021) 110639.
\newblock \href {https://doi.org/10.1016/J.JCP.2021.110639}
  {\path{doi:10.1016/J.JCP.2021.110639}}.

\bibitem{VTH19}
P.~S. Verma, J.-M. Teissier, O.~Henze, W.-C. Müller, Fourth-order accurate
  finite-volume {CWENO} scheme for astrophysical {MHD} problems, Monthly
  Notices of the Royal Astronomical Society 482 (2019) 416--437.
\newblock \href {https://doi.org/10.1093/mnras/sty2641}
  {\path{doi:10.1093/mnras/sty2641}}.

\bibitem{COC11}
P.~McCorquodale, P.~Colella, A high-order finite-volume method for conservation
  laws on locally refined grids, Communications in Applied Mathematics and
  Computational Science 6 (2011) 1--25.
\newblock \href {https://doi.org/10.2140/camcos.2011.6.1}
  {\path{doi:10.2140/camcos.2011.6.1}}.

\bibitem{BUH14}
P.~Buchmüller, C.~Helzel, Improved accuracy of high-order {WENO} finite volume
  methods on cartesian grids, Journal of Scientific Computing 61 (2014)
  343--368.
\newblock \href {https://doi.org/10.1007/S10915-014-9825-1}
  {\path{doi:10.1007/S10915-014-9825-1}}.

\bibitem{ROM16}
J.~N. de~la Rosa, C.-D. Munz, {XTROEM-FV}: a new code for computational
  astrophysics based on very high order finite-volume methods – {I}.
  magnetohydrodynamics, Monthly Notices of the Royal Astronomical Society 455
  (2016) 3458--3479.
\newblock \href {https://doi.org/10.1093/mnras/stv2531}
  {\path{doi:10.1093/mnras/stv2531}}.

\bibitem{BAG16}
D.~S. Balsara, T.~Amano, S.~Garain, J.~Kim, A high-order relativistic two-fluid
  electrodynamic scheme with consistent reconstruction of electromagnetic
  fields and a multidimensional {R}iemann solver for electromagnetism, Journal
  of Computational Physics 318 (2016) 169--200.
\newblock \href {https://doi.org/10.1016/J.JCP.2016.05.006}
  {\path{doi:10.1016/J.JCP.2016.05.006}}.

\bibitem{WSH19}
K.~Wu, C.~W. Shu, Provably positive high-order schemes for ideal
  magnetohydrodynamics: analysis on general meshes, Numerische Mathematik 142
  (2019) 995--1047.
\newblock \href {https://doi.org/10.1007/S00211-019-01042-W}
  {\path{doi:10.1007/S00211-019-01042-W}}.

\bibitem{EVH88}
C.~R. Evans, J.~F. Hawley, Simulation of magnetohydrodynamic flows - a
  constrained transport method, The Astrophysical Journal 332 (1988) 659.
\newblock \href {https://doi.org/10.1086/166684} {\path{doi:10.1086/166684}}.

\bibitem{ZIE04}
U.~Ziegler, A central-constrained transport scheme for ideal
  magnetohydrodynamics, Journal of Computational Physics 196 (2004) 393--416.
\newblock \href {https://doi.org/10.1016/J.JCP.2003.11.003}
  {\path{doi:10.1016/J.JCP.2003.11.003}}.

\bibitem{BMD13}
D.~S. Balsara, C.~Meyer, M.~Dumbser, H.~Du, Z.~Xu, Efficient implementation of
  {ADER} schemes for {E}uler and magnetohydrodynamical flows on structured
  meshes – speed comparisons with {R}unge–{K}utta methods, Journal of
  Computational Physics 235 (2013) 934--969.
\newblock \href {https://doi.org/10.1016/J.JCP.2012.04.051}
  {\path{doi:10.1016/J.JCP.2012.04.051}}.

\bibitem{LOZ00}
P.~Londrillo, L.~D. Zanna, High‐order upwind schemes for multidimensional
  magnetohydrodynamics, The Astrophysical Journal 530 (2000) 508--524.
\newblock \href {https://doi.org/10.1086/308344} {\path{doi:10.1086/308344}}.

\bibitem{HEO87}
A.~Harten, B.~Engquist, S.~Osher, S.~R. Chakravarthy, Uniformly high order
  accurate essentially non-oscillatory schemes, {III}, Journal of Computational
  Physics 71 (1987) 231--303.
\newblock \href {https://doi.org/10.1016/0021-9991(87)90031-3}
  {\path{doi:10.1016/0021-9991(87)90031-3}}.

\bibitem{LOC94}
X.~D. Liu, O.~S, C.~T, Weighted essentially non-oscillatory schemes, Journal of
  Computational Physics 115 (1994) 200--212.
\newblock \href {https://doi.org/10.1006/JCPH.1994.1187}
  {\path{doi:10.1006/JCPH.1994.1187}}.

\bibitem{JIS96}
G.~S. Jiang, C.~W. Shu, Efficient implementation of weighted {ENO} schemes,
  Journal of Computational Physics 126 (1996) 202--228.
\newblock \href {https://doi.org/10.1006/JCPH.1996.0130}
  {\path{doi:10.1006/JCPH.1996.0130}}.

\bibitem{SHU97}
C.-W. Shu, \href{https://ntrs.nasa.gov/citations/19980007543}{Essentially
  non-oscillatory and weighted essentially non-oscillatory schemes for
  hyperbolic conservation laws}, ICASE Report 97-65 (1997).
\newline\urlprefix\url{https://ntrs.nasa.gov/citations/19980007543}

\bibitem{BAS00}
D.~S. Balsara, C.~W. Shu, Monotonicity preserving weighted essentially
  non-oscillatory schemes with increasingly high order of accuracy, Journal of
  Computational Physics 160 (2000) 405--452.
\newblock \href {https://doi.org/10.1006/JCPH.2000.6443}
  {\path{doi:10.1006/JCPH.2000.6443}}.

\bibitem{GSV09}
G.~A. Gerolymos, D.~Sénéchal, I.~Vallet, Very-high-order {WENO} schemes,
  Journal of Computational Physics 228 (2009) 8481--8524.
\newblock \href {https://doi.org/10.1016/J.JCP.2009.07.039}
  {\path{doi:10.1016/J.JCP.2009.07.039}}.

\bibitem{LPR99}
D.~Levy, G.~Puppo, G.~Russo, Central {WENO} schemes for hyperbolic systems of
  conservation laws, ESAIM: Mathematical Modelling and Numerical Analysis 33
  (1999) 547--571.
\newblock \href {https://doi.org/10.1051/M2AN:1999152}
  {\path{doi:10.1051/M2AN:1999152}}.

\bibitem{RUS61}
V.~V. Rusanov, The calculation of the interaction of non-stationary shock waves
  and obstacles, USSR Computational Mathematics and Mathematical Physics 1
  (1962) 304--320.
\newblock \href {https://doi.org/10.1016/0041-5553(62)90062-9}
  {\path{doi:10.1016/0041-5553(62)90062-9}}.

\bibitem{SHO89}
C.~W. Shu, S.~Osher, Efficient implementation of essentially non-oscillatory
  shock-capturing schemes, {II}, Journal of Computational Physics 83 (1989)
  32--78.
\newblock \href {https://doi.org/10.1016/0021-9991(89)90222-2}
  {\path{doi:10.1016/0021-9991(89)90222-2}}.

\bibitem{BAL10}
D.~S. Balsara, Multidimensional {HLLE} {R}iemann solver: Application to {E}uler
  and magnetohydrodynamic flows, Journal of Computational Physics 229 (2010)
  1970--1993.
\newblock \href {https://doi.org/10.1016/J.JCP.2009.11.018}
  {\path{doi:10.1016/J.JCP.2009.11.018}}.

\bibitem{BNK17}
D.~S. Balsara, B.~Nkonga, Multidimensional riemann problem with self-similar
  internal structure – part {III} – a multidimensional analogue of the
  {HLLI} {R}iemann solver for conservative hyperbolic systems, Journal of
  Computational Physics 346 (2017) 25--48.
\newblock \href {https://doi.org/10.1016/J.JCP.2017.05.038}
  {\path{doi:10.1016/J.JCP.2017.05.038}}.

\bibitem{LOZ04}
P.~Londrillo, L.~D. Zanna, On the divergence-free condition in {G}odunov-type
  schemes for ideal magnetohydrodynamics: the upwind constrained transport
  method, Journal of Computational Physics 195 (2004) 17--48.
\newblock \href {https://doi.org/10.1016/J.JCP.2003.09.016}
  {\path{doi:10.1016/J.JCP.2003.09.016}}.

\bibitem{MCL99}
M.~M. Low, The energy dissipation rate of supersonic, magnetohydrodynamic
  turbulence in molecular clouds, The Astrophysical Journal 524 (1999)
  169--178.
\newblock \href {https://doi.org/10.1086/307784} {\path{doi:10.1086/307784}}.

\bibitem{KWN13}
A.~G. Kritsuk, R.~Wagner, M.~L. Norman, Energy cascade and scaling in
  supersonic isothermal turbulence, Journal of Fluid Mechanics 729 (2013) R1.
\newblock \href {https://doi.org/10.1017/JFM.2013.342}
  {\path{doi:10.1017/JFM.2013.342}}.

\bibitem{ALU13}
H.~Aluie, Scale decomposition in compressible turbulence, Physica D: Nonlinear
  Phenomena 247 (2013) 54--65.
\newblock \href {https://doi.org/10.1016/j.physd.2012.12.009}
  {\path{doi:10.1016/j.physd.2012.12.009}}.

\bibitem{GOS98}
S.~Gottlieb, C.-W. Shu, Total variation diminishing {R}unge-{K}utta schemes,
  Mathematics of Computation 67 (1998) 73--85.
\newblock \href {https://doi.org/10.1090/S0025-5718-98-00913-2}
  {\path{doi:10.1090/S0025-5718-98-00913-2}}.

\bibitem{GST01}
S.~Gottlieb, C.~W. Shu, E.~Tadmor, Strong stability-preserving high-order time
  discretization methods, SIAM Review 43 (2001) 89--112.
\newblock \href {https://doi.org/10.1137/S003614450036757X}
  {\path{doi:10.1137/S003614450036757X}}.

\bibitem{KET08}
D.~I. Ketcheson, Highly efficient strong stability-preserving {R}unge–{K}utta
  methods with low-storage implementations, SIAM Journal on Scientific
  Computing 30 (2008) 2113--2136.
\newblock \href {https://doi.org/10.1137/07070485X}
  {\path{doi:10.1137/07070485X}}.

\bibitem{KRA91}
J.~F. Kraaijevanger, Contractivity of {R}unge-{K}utta methods, BIT Numerical
  Mathematics 31 (1991) 482--528.
\newblock \href {https://doi.org/10.1007/BF01933264}
  {\path{doi:10.1007/BF01933264}}.

\bibitem{RSP02}
S.~J. Ruuth, R.~J. Spiteri, Two barriers on strong-stability-preserving time
  discretization methods, Journal of Scientific Computing 17 (2002) 211--220.
\newblock \href {https://doi.org/10.1023/A:1015156832269}
  {\path{doi:10.1023/A:1015156832269}}.

\bibitem{KGM11}
D.~I. Ketcheson, S.~Gottlieb, C.~B. Macdonald, Strong stability preserving
  two-step {R}unge–{K}utta methods, SIAM Journal on Numerical Analysis 49
  (2011) 2618--2639.
\newblock \href {https://doi.org/10.1137/10080960X}
  {\path{doi:10.1137/10080960X}}.

\bibitem{BGG17}
C.~Bresten, S.~Gottlieb, Z.~Grant, D.~Higgs, D.~I. Ketcheson, A.~Németh,
  Explicit strong stability preserving multistep {R}unge–{K}utta methods,
  Mathematics of Computation 86 (2017) 747--769.
\newblock \href {https://doi.org/10.1090/mcom/3115}
  {\path{doi:10.1090/mcom/3115}}.

\bibitem{SSPSITE}
S.~Gottlieb, D.~Higgs, D.~I. Ketcheson,
  \href{http://www.sspsite.org/msrk.html}{Multistep multistage ({MSRK})
  methods}.
\newline\urlprefix\url{http://www.sspsite.org/msrk.html}

\bibitem{BAL04}
D.~S. Balsara, Second‐order–accurate schemes for magnetohydrodynamics with
  divergence‐free reconstruction, The Astrophysical Journal Supplement Series
  151 (2004) 149--184.
\newblock \href {https://doi.org/10.1086/381377} {\path{doi:10.1086/381377}}.

\bibitem{MTB10}
A.~Mignone, P.~Tzeferacos, G.~Bodo, High-order conservative finite difference
  {GLM}–{MHD} schemes for cell-centered {MHD}, Journal of Computational
  Physics 229 (2010) 5896--5920.
\newblock \href {https://doi.org/10.1016/J.JCP.2010.04.013}
  {\path{doi:10.1016/J.JCP.2010.04.013}}.

\bibitem{LBA22}
G.~Leidi, C.~Birke, R.~Andrassy, J.~Higl, P.~V. Edelmann, G.~Wiest,
  C.~Klingenberg, F.~K. Röpke, Astronomy and Astrophysics\href
  {https://doi.org/10.1051/0004-6361/202244665}
  {\path{doi:10.1051/0004-6361/202244665}}.

\bibitem{COW84}
P.~Colella, P.~R. Woodward, The piecewise parabolic method ({PPM}) for
  gas-dynamical simulations, Journal of Computational Physics 54 (1984)
  174--201.
\newblock \href {https://doi.org/10.1016/0021-9991(84)90143-8}
  {\path{doi:10.1016/0021-9991(84)90143-8}}.

\bibitem{BAL12}
D.~S. Balsara, Self-adjusting, positivity preserving high order schemes for
  hydrodynamics and magnetohydrodynamics, Journal of Computational Physics 231
  (2012) 7504--7517.
\newblock \href {https://doi.org/10.1016/J.JCP.2012.01.032}
  {\path{doi:10.1016/J.JCP.2012.01.032}}.

\end{thebibliography}

\end{document}